# Finite functions and the necessary use of large cardinals

By Harvey M. Friedman*

## Table of Contents

*This research was partially supported by NSF. AMS Classification number 03, 04, 05

## Preface

We present a coherent collection of finite mathematical theorems some of which can only be proved by going well beyond the usual axioms for mathematics. The proofs of these theorems illustrate in clear terms how one uses the well studied higher infinities of abstract set theory called large cardinals in an essential way in order to derive results in the context of the natural numbers. The findings raise the specific issue of what consitutes a valid mathematical proof and the general issue of objectivity in mathematics in a down to earth way.

Large cardinal axioms, which go beyond the usual axioms for mathematics, have been commonly used in abstract set theory since the 1960's (e.g., see [Sc61], [MS89]). We believe that the results reported on here are the early stages of an evolutionary process in which new axioms for mathematics will be commonly used in an essential way in the more concrete parts of mathematics.

We wish to thank John Burgess, Tim Carlson, Randy Dougherty, Charles Fefferman, Ronald Graham, Leo Harrington, Carl Jockusch, Si Kochen, Rich Laver, Tony Martin, Barry Mazur, Gian-Carlo Rota, Paul Sally, Joe Shipman, Steve Simpson, Ted Slaman, Joel Spencer, Robert Stanton, and John Steel for valuable conversations and encouragement.

## Abstract

We begin by presenting a new kind of finite counting theorem which asserts that any function on vectors of natural numbers has few low values over finite sets of arbitrary size (Theorems 0.1 and 0.2 from the introduction).

This estimate on the low values (herein called regressive values) is linear in the sizes of the finite sets.

This basic counting theorem is extended to systems of finite functions, where it is shown that every system which is "decreasing" (in appropriate senses) includes functions which have even fewer low values over finite sets of arbitrary size. This improved estimate is a constant depending on the dimension and not on the sizes of the finite sets (see Proposition A from the introduction).

However, these extended results can only be proved by using additional axioms which go well beyond the currently accepted axioms for mathematics (as usually formalized by ZFC = Zermelo Frankel set theory with the axiom of choice). Furthermore, the extended results have obviously equivalent finite forms, which remove any mention of infinite sets (see Proposition B from the introduction). See Theorem 5.91 and Corollary 1.

The additional axioms used in the proof of these extended results are among the so called large cardinal axioms that have been extensively used in set theory since the 1960's (see e.g., [Sc61], [MS89]). The proofs here illustrate in clear terms how one uses large cardinals in an essential and completely natural way in the integers. A conceptual overview of the method is described at the end of the introduction.

The quest for a simple meaningful finite mathematical theorem that can only be proved by going beyond the usual axioms for mathematics has been a goal in the foundations of mathematics since Gödel's work in the 1930's [Go31].

Prior relevant developments include the introduction of the forcing method in [Co66], which was used to establish the unprovability of the continuum hypothesis from the usual axioms of mathematics (the consistency having been previously established in [Go40]). But for theoretical reasons, the forcing method cannot establish the independence of a finite statement, at least not directly (see [Fr92, p. 54] and "absoluteness" in [Je78, p. 101–108]).

Another relevant development was the first example of a simple meaningful finite mathematical theorem (in Ramsey theory) whose proof requires some weak use of infinite sets, in [PH77]. This was followed by a simple meaningful finite mathematical theorem in graph theory (related to J. B. Kruskal's theorem) whose proof requires some weak use of uncountable sets, [Si85], and a further, simply stated theorem in graph theory (related to the graph minor theorem), [FRS87], whose proof requires stronger uses of uncountable sets.

Other relevant developments include the first example of a simple meaningful mathematical theorem involving Borel measurable functions whose proof requires use of uncountably many infinities, [Fr71], and the first example of a simple meaningful mathematical theorem involving Borel measurable functions show that proof requires axioms going far beyond the usual axioms for mathematics, [Fr81], [St85].

As is clear from the body of this abstract, the results – which are provable only by going well beyond the usual axioms for mathematics – are well connected with existing topics in finite and discrete combinatorics (Ramsey theory), and form a coherent body of examples. We anticipate that a variety of simple and basic examples will be found with increasingly strong connections to diverse areas of mathematics.

**Introduction**

We let $N$ be the set of all nonnegative integers. For $x \in N^k$, we let $\min(x)$ be the minimum coordinate of $x$, and $|x|$ be the maximum coordinate of $x$ (sup norm). We also use $|\ |$ for the cardinality of finite sets.

For $n \in N$, we write $[n]$ for $\{0, 1, \ldots, n-1\}$. For $k > 0$ and $X$ a set, we write $S(X)$ for the set of all subsets of $X$, and $S_k(X)$ for the set of all $k$ element subsets of $X$.

We begin with a statement of the usual classical infinite Ramsey theorem (IRT).

Infinite Ramsey Theorem. *Let $k, r > 0$ and $F\colon S_k(N) \to [r]$. Then there exists an infinite $E \subseteq N$ such that $F$ is constant on $S_k(E)$.*

We also state the usual finite Ramsey theorem (FRT). See [Ra30] and [GRS80] for IRT and FRT (infinite and finite Ramsey theorem).

Finite Ramsey Theorem. *Let $n \gg k, r, p > 0$ and $F\colon S_k[n] \to [r]$. Then there exists $E \in S_p[n]$ such that $F$ is constant on $S_k(E)$.*

Here, as elsewhere, we use $\gg$ for "sufficiently greater than." In this use, the full expansion of the first sentence would read:

For all $k, r, p > 0$ there exists $m$ such that for all
$n \geq m$ and $F\colon S_k[n] \to [r]$, the following holds.

Of course, in this case, as well as in all subsequent cases, this can be simplified to read:

For all $k, r, p > 0$ there exists $n$ such that for all
$F\colon S_k[n] \to [r]$, the following holds.

The first simple meaningful finite mathematical theorem that is independent of Peano Arithmetic is given in [PH77], Theorem 1.2:

Paris-Harrington-Ramsey Theorem. *Let $n \gg k, r, p > 0$ and $F\colon S_k[n] \to [r]$. Then there exists $E \subseteq [n]$, $|E| \geq p$, $\min(E)$, such that $F$ is constant on $S_k(E)$.*

We now present two counting theorems of a different flavor from the PH Ramsey Theorem, but which share its metamathematical properties. We first give their infinite forms (Theorems 0.1 and 0.2), and then follow it by their obvious finite forms (Theorems 0.3 and 0.4). Theorem 0.2 leads to the simple meaningful mathematical theorems of Part 4 (Propositions A–D) which can only be proved by going well beyond the usual axioms for mathematics.

Theorem 0.1. *Let $k, r, p > 0$ and $F\colon N^k \to N^r$ obey the inequality*

$$|F(x)| \leq \min(x)\,.$$

*Then there exists $E \in S_p(N)$ such that*

$$|F[E^k]| \leq (k^k)p\,.$$

We now turn this around so that it asserts a combinatorial property of any function $F\colon N^k \to N^r$.

Let $A, B \subseteq N^k$, and $F\colon A \to N^r$. We say that $y$ is a regressive value of $F$ on $B$ if and only if there exists $x \in B$ such that $F(x) = y$ and $|y| < \min(x)$.

*Example.* Let $F: N^2 \to N$ be given by $F(x,y) = (x-y)^2$. Then every square is a regressive value of $F$ on $N^2$. But the only regressive value of $F$ on $\{1,2,4,\dots\}^2$ is 0.

THEOREM 0.2. *Let $k,r,p > 0$ and $F: N^k \to N^r$. Then $F$ has $\leq (k^k)p$ regressive values on some $E^k \subseteq N^k$, $|E| = p$.*

We now state the obvious finite forms of Theorems 0.1 and 0.2:

THEOREM 0.3. *Let $n \gg k,r,p > 0$ and $F: [n]^k \to [n]^r$ obey the inequality*

$$|F(x)| \leq \min(x).$$

*Then there exists $E \in S_p[n]$ such that*

$$|F[E^k]| \leq (k^k)p.$$

THEOREM 0.4. *Let $n \gg k,r,p > 0$ and $F: [n]^k \to [n]^r$. Then $F$ has $\leq (k^k)p$ regressive values on some $E^k \subseteq [n]^k$, $|E| = p$.*

We could have used $\leq$ instead of $<$ in the definition of regressive value. The reason we used $<$ is because of the following set theoretic result, which is a hint of things to come.

THEOREM 0.5. *Let $k,r,p > 0$ and $F: \lambda^k \to \lambda^r$, where $\lambda$ is a suitably large cardinal. Then $F$ has $\leq k^k$ regressive values on some $E^k \subseteq \lambda^k$, $|E| = p$.*

A function assignment for a set $X$ is a mapping $U$ which assigns to each finite subset $A$ of $X$, a unique function

$$U(A): A \to A.$$

The following is easily obtained from Theorem 0.4.

THEOREM 0.6. *Let $k,p > 0$ and $U$ be a function assignment for $N^k$. Then some $U(A)$ has $\leq (k^k)p$ regressive values on some $E^k \subseteq A$, $|E| = p$.*

We want to place a natural condition on function assignments for $N^k$ so that we get the improved estimate appearing in Theorem 0.5.

The main condition that we use is #-decreasing, which is stated informally as follows:

When a vector is inserted, the function extends or drops higher up.

More formally, let $U$ be a function assignment for $N^k$.

We say that $U$ is #-decreasing if and only if for all finite $A \subseteq N^k$ and $x \in N^k$,

$$\text{either } U(A) \subseteq U(A \cup \{x\}) \text{ or there exists}$$
$$|y| > |x| \text{ such that } |U(A)(y)| > |U(A \cup \{x\})(y)|.$$

#-decreasing is stronger than it appears, and is quite stable. It implies a rather complete picture of the relationship between $U(A)$ and $U(A \cup \{x\})$. For instance, $U$ is #-decreasing if and only if for all finite $A \subseteq N^k$ and $x \in N^k$, either $U(A) \subseteq U(A \cup \{x\})$, or there exists $|y| > |x|$ such that

(i) $|U(A)(y)| > |U(A \cup \{x\})(y)|$;
(ii) for all $z \in A$, if $|z| < |y|$, then $U(A)(z) = U(A \cup \{x\})(z)$;
(iii) for all $z \in A$, if $|z| = |y|$, then $U(A)(z) = U(A \cup \{x\})(z)$ or $|U(A)(z)| > |U(A \cup \{x\})(z)|$.

We also consider lex #-decreasing function assignments $U$ for $N^k$. Lex #-decreasing is defined in exactly the same way as #-decreasing, except $>$ is replaced by $>_{\text{lex}}$, where $>_{\text{lex}}$ is according to the usual lexicographical ordering of $N^k$.

The main result is that the following proposition is independent of ZFC, and can be proved used certain large cardinal axioms, but not without them:

Proposition A. *Let $k, p > 0$ and $U$ be a #-decreasing (lex #-decreasing) function assignment for $N^k$. Then some $U(A)$ has $\leq k^k$ regressive values on some $E^k \subseteq A$, $|E| = p$.*

At the end of this introduction, we indicate the general strategy for proving Proposition A from large cardinals.

We now state the obvious finite form of Proposition A. Clearly #-decreasing and lex #-decreasing make perfectly good sense for function assignments for $[n]^k$.

Proposition B. *Let $n \gg k, p > 0$ and $U$ be a #-decreasing (lex #-decreasing) function assignment for $[n]^k$. Then some $U(A)$ has $\leq k^k$ regressive values on some $E^k \subseteq A$, $|E| = p$.*

Proposition B immediately implies Proposition A. Also Proposition B can be derived from Proposition A using a standard compactness (finitely branching tree) argument.

Note that #-decreasing is based on the ordering of $N^k$ by the sup norm and lex #-decreasing is based on the lexicographic ordering. We now give an abstract form of Proposition A which uses other orderings.

A strict order on $N^k$ is a transitive irreflexive binary relation $<$ on $N^k$.

We write $\leq_c$ for the order on $N^k$ defined by: $x \leq_c y$ if and only if for all $1 \leq i \leq k$, $x_i \leq y_i$.

We say that $<$ is an upward order on $N^k$ if and only if $<$ is a strict order on $N^k$ such that $x \leq_c y$ implies $\neg(y < x)$. Observe that $|x| < |y|$ and the strict lexicographic orderings on $N^k$ are upward orderings on $N^k$.

Let $<_1$, $<_2$ be orders on $N^k$ and $U$ be a function assignment for $N^k$.

We say that $U$ is $<_1, <_2$-#-decreasing if and only if for all finite $A \subseteq N^k$ and $x \in N^k$,

$$\text{either } U(A) \subseteq U(A \cup \{x\}) \text{ or there exists}$$
$$y >_1 x \text{ such that } U(A)(y) >_2 U(A \cup \{x\})(y).$$

PROPOSITION C. *Let $k, p > 0$, $<_1$, $<_2$ be upward orders on $N^k$, and $U$ be a $<_1, <_2$-#-decreasing function assignment for $N^k$. Then some $U(A)$ has $\leq k^k$ regressive values on some $E^k \subseteq A$, $|E| = p$.*

We also have the following equivalent finite form:

PROPOSITION D. *Let $n \gg k, p > 0$, $<_1$, $<_2$ be upward orders on $[n]^k$, and $U$ be a $<_1, <_2$-#-decreasing function assignment for $[n]^k$. Then some $U(A)$ has $\leq k^k$ regressive values on some $E^k \subseteq A$, $|E| = p$.*

We conjecture that Propositions A–D are equivalent to a single metamathematical statement which asserts the 1-consistency of ZFC $+\{$there exists a $k$-subtle cardinal$\}_k$. See the discussion of large cardinals in Part 2. In this manuscript, we will prove that all forms of Propositions A–D imply the consistency of ZFC $+ \{$there exists a $k$-subtle cardinal$\}_k$, and that they are all provable from ZFC+ "for all $k$ there exists a $k$-subtle cardinal."

In Part 2 we give ten presentations of the hierarchy of cardinals used in the previous paragraph, as well as nine formulations of the existence of the hierarchy taken as a whole (Th. 2.11). In the previous paragraph, any of these can be used instead of talking about $k$-subtle cardinals. From the set theoretic point of view, a very appealing choice is $k -$ SRP (the stationary Ramsey property of order $k$), which is stated as follows:

Let $S_k(X)$ be the set of all $k$ element subsets of the set $X$. We say that an infinite cardinal $\lambda$ is $k -$ SRP if and only if for all $F\colon S_k(\lambda) \to \{0, 1\}$ there exists stationary $E \subseteq \lambda$ such that $F$ is constant on $S_k(E)$.

The $k$-critical linear orderings are natural from the point of view of general mathematics. Let $(X, \leq)$ be a linear ordering with no endpoints. We say that $F\colon X^k \to X$ is regressive if and only if for all $x_1, \ldots, x_k$, we have $F(x_1, \ldots, x_k) < \min(x_1, \ldots, x_k)$. Now $(X, \leq)$ is $k$-critical if and only if for all regressive $F\colon X^k \to X$ there exist $b_1 < \cdots < b_{k+1}$ such that $F(b_1, \ldots, b_k) = F(b_2, \ldots, b_{k+1})$.

The general strategy for using large cardinals in the integers is as follows.

We start with a discrete (or finite) structure $X$ obeying certain hypotheses H. We wish to prove that a certain kind of finite configuration occurs in $X$, assuming that H holds. We define a suitable concept of completion in the context of arbitrary linearly ordered sets. We verify that if $X$ has a completion with the desired kind of finite configuration, then $X$ already has the desired kind of finite configuration. We then show, using hypotheses H, that $X$ has completions of every well-ordered type. We now use the existence of a suitably large cardinal $\lambda$. Using large cardinal combinatorics, we show that in any completion of order type $\lambda$, the desired kind of finite configuration exists. Hence the desired kind of finite configuration already exists in $X$.

In the case of Proposition A, we can spell out the strategy in more specific terms as follows.

(1) Show that #-decreasing implies some additional conditions, the most important of which is end-preserving (Theorem 3.10).

(2) Using the classical infinite Ramsey theorem, show that every #-decreasing function assignment can be made "uniform" so that every function is order isomorphic to an original function (Theorem 4.2).

(3) Define completion of function assignments. These are linearly ordered functions (functions from a Cartesian power of a linearly ordered set into itself) such that every finite set has a finite superset on which the function is order isomorphic to an element of the function assignment (immediately following the proof of Theorem 4.8).

(4) Show that every uniform function assignment has (unique) completions of every well-ordered type (Theorem 4.9). The crucial argument by transfinite recursion is entirely natural (Theorem 4.8). This is called the ordinal completion property.

(5) Finally we use the large cardinal. Show that every function on the large cardinal has few regressive values on Cartesian powers of every finite size (Lemma 4.10). This is a straightforward argument in large cardinal combinatorics.

(6) Conclude that some element of the function assignment has the desired property by quoting the definition of completion (Theorem 4.11).

## 1. Counting values

In this part we show that Theorems 0.3 and 0.4 from the introduction are independent of Peano Arithmetic (PA) and have other metamathematical properties. Theorems 0.2 and 0.4 are extended in Part 4 to finite combinatorial results that can only be proved by using large cardinals.

We say that $x, y \in N^k$ are of the same order type if and only if for all $1 \leq i, j \leq k$, $x_i < x_j$ if and only if $y_i < y_j$. Other notation was explained in the introduction.

Let $\mathrm{ot}(k)$ be the number of order types of elements of $N^k$. It is obvious that $\mathrm{ot}(k) \leq k^k$ (every element on $N^k$ has the same order type as an element of $[k]^k$), and a straightforward inductive argument shows that $\mathrm{ot}(k) \leq 2^k(k!)$.

It is easy to see that $\text{ot}(k) = \sum_{1 \le i \le k} f(k,i)$, where $f(k,i)$ is the number of surjective maps from $[k]$ onto $[i]$. In [Gr62] it is shown that $\text{ot}(k)$ is asymptotic to $k!/2\ln^{k+1} 2$.

Let $F\colon S_k(N) \to N$. We say that $F$ is regressive if and only if for all $x \in S_k(N)$, if $\min(x) > 0$ then $F(x) < \min(x)$.

[KM87] analyzes the following result.

THEOREM I. *Let $n \gg k, p > 0$ and $F\colon S_k[n] \to [n]$ be regressive. Then there exists $E \in S_p[n]$ such that the following holds. For all $x, y \in S_k(E)$, if $\min(x) = \min(y)$, then $F(x) = F(y)$.*

Theorem I is a variant of Proposition 2.10 in [PH77], which was shown to be unprovable in PA. In [PH77], the Paris-Harrington statement (Theorem 1.2) was shown to imply Proposition 2.10, thereby obtaining the unprovability of Theorem 1.2 from PA = Peano Arithmetic.

The following theorem is a weak form of Theorem 0.4, which in turn follows from Theorem 0.3.

THEOREM II. *For all $k > 0$ there exists $c > 0$ such that the following holds. Let $n \gg k, p > 0$ and $F\colon [n]^k \to [n]^2$. Then there exists $E \in S_p[n]$ such that $F$ has $\le cp$ regressive values on $E^k$.*

The next theorem is a strong form of Theorems I and 0.3.

THEOREM III. *Let $n \gg k, r, p > 0$ and $F\colon [n]^k \to N^r$ obey the inequality*

$$|F(x)| \le \min(x)\,.$$

*Then there exists $E \in S_p[n]$ such that the following holds. For all $x, y \in E^k$ with the same order type, if $\min(x) = \min(y)$, then $F(x) = F(y)$.*

We write 1-Con$(T)$ for the 1-consistency of the formal system $T$, which asserts that every $\Sigma^0_1$ sentence provable in $T$ is true.

Define the function $H_{\text{I}}(k,p) =$ the least $n$ such that Theorem I holds for $n$, $k$, $p$. Write $H_{\text{I},k}(p) = H_{\text{I}}(k,p)$.

LEMMA 1.1 ([KM87]). *Theorem* I *is not provable in* PA. *For each fixed $k$, Theorem* I *is provable in* PA. *It is provable in* PA *that Theorem* I *is equivalent to the* 1*-consistency of* PA. *Every $< \varepsilon_0$-recursive function is eventually $<$ some $H_{\text{I},k}$. Every $H_{\text{I},k}$ is a $< \varepsilon_0$-recursive function.*

We wish to obtain analogous information for Theorems II, III, 0.3, and 0.4.

For $E \subseteq N$ let $E_i$ be the $i^{\text{th}}$ element of $E$; $E_1$ is the least element of $E$.

For $n, r \in N$, write $n^*r$ for $(n, \ldots, n)$, where there are $r$ $n$'s.

Lemma 1.2. *The following is provable in* PA. *Theorem* III *implies Theorem* I *and Theorem* 0.3. *Theorem* 0.3 *implies Theorem* 0.4 *implies Theorem* II. *In fact, Theorem* 0.3 *for* $r = 2$ *implies Theorem* 0.4 *for* $r = 2$ *implies Theorem* II.

*Proof.* It is obvious that Theorem III implies Theorem 0.3. To see that Theorem III implies Theorem I, let $k, p > 0$. Let $n$ be as given by Theorem III for $k$, 1, $p+1$. Let $F\colon S_k[n] \to N$ be regressive. Define $F'\colon [n]^k \to N$ by $F'(x) = F(\mathrm{rng}(x))$ if $x$ is strictly increasing and $\min(x) > 0$; 0 otherwise. Let $E$ be as given by Theorem III. Then $E\backslash\{E_1\}$ is as required by Theorem I.

Note that the last implication in claims 2 and 3 are immediate, taking $c = k^k$. So it remains to see that for fixed $r$, Theorem 0.3 implies Theorem 0.4. Accordingly, assume Theorem 0.3 for $r$, and let $k, r, p > 0$. Let $n$ work for $k$, $r$, $p$ in Theorem 0.3. Let $F\colon [n]^k \to [n]^r$. Define $F'\colon [n]^k \to [n]^r$ by $F'(x) = F(x)$ if $|F(x)| \leq \min(x)$; $0^*r$ otherwise. Then Theorem 0.3 applies to $F'$ and so let $E \in S_p[n]$, where $|F'[E^k]| \leq (k^k)p$. Note that every regressive value of $F$ on $E^k$ is a value of $F'$ on $E$. Hence the number of regressive values of $F$ on $E^k$ is $\leq (k^k)p$ as required. □

We now wish to prove that Theorem I implies Theorem III. For this purpose, we introduce another statement.

Theorem IV. *Let* $n \gg k, p > 0$ *and* $F\colon S_k[n] \to [n]$ *be regressive. Then there exists* $E \in S_p[n]$ *such that the following holds.*
(i) *For all* $x, y \in S_k(E)$, *if* $\min(x) = \min(y)$ *then* $F(x) = F(y)$;
(ii) *For all* $x < y$ *from* $E$, $2^x < y$.

Lemma 1.3. *The following is provable in* PA. *Theorem* I *implies Theorem* IV.

*Proof.* Let $k, p > 0$. Without loss of generality assume $k, p > 1$. Apply Theorem I for $k$, $8p + 1 + k$ to obtain $n$ such that the conclusion of Theorem I holds. To prove that $n$ works for $k, p$ in Theorem IV, let $F\colon S_k[n] \to [n]$ be regressive.

We define a regressive $G\colon S_k[n] \to [n]$ by cases, where the operative case is the first one that applies. Let $A \in S_k[n]$.

*Case* 1. $A_2 - A_1 < A_1$. Define $G(A) = A_2 - A_1$.

*Case* 2. There exists $j < A_1$ such that $2^j > A_2$. Define $G(A) =$ greatest such $j$.

*Case* 3. Otherwise. Define $G(A) = F(A)$.

Let $E$ be as given by Theorem I for $n$, $k$, $8p + 3 + k$. Let $1 \leq i \leq 8p + 1$. Then by comparing $G(E_i, E_{i+1}, \ldots, E_{i+k-1})$ and $G(E_i, E_{i+2}, \ldots, E_{i+k})$, we see that $E_{i+2} - E_i \geq E_i$, and so $E_{i+2} \geq 2E_i$. Similarly, we have $E_{i+4} \geq 2E_{i+2}$. By comparing $G(E_i, E_{i+2}, \ldots, E_{i+k})$ and $G(E_i, E_{i+4}, \ldots, E_{i+k+2})$, we see that $2^{E_{i-1}} \leq E_{i+4}$. From this it follows that $\{E_9, E_{17}, \ldots, E_{8p+1}\}$ is as required by Theorem IV. □

LEMMA 1.4. *The following is provable in* PA. *Theorem* IV *implies Theorem* III.

*Proof.* Let $k, r, p > 0$. Without loss of generality we may assume $k, r, p > 1$. Apply Theorem IV for $k$ and some appropriate exponential expression $\alpha(k, r, p)$, to obtain $n$ such that the conclusion of Theorem IV holds. Let $F\colon [n]^k \to [n]^r$ obey the inequality in Theorem III. We define $G\colon S_k[n] \to [n]$ by cases. Let $A \in S_k[n]$.

*Case* 1. For all $1 < i \leq k$, $A_i > 2^{A_{i-1}}$ and $A_1 > \alpha(k, r, p)$. For each tuple $m_1, \ldots m_k$ from $\{1, \ldots, k\}$ consider $F(\log(A_{m_1}), \ldots, \log(A_{m_k}))$, where we round down to the nearest integer. This correspondence can be coded as an integer $< A_1$. $G(A)$ is the code.

*Case* 2. Otherwise. Set $G(A) = 0$.

Let $E$ be as given by Theorem IV. Then the largest $k$ elements of $\{\log(v)\colon v \in E\}$ is as required by Theorem III. □

LEMMA 1.5. *The following is provable in* PA. *Theorem* II *implies Theorem* I.

*Proof.* Fix $k, p > 0$. We apply Theorem II with dimension $k$ to obain $c$ such that the last two lines of Theorem II hold. We fix $p' \gg k, p, c$ (according to the requirements of the application of FRT below). Let $n$ be according to Theorem II, with $k$, $p'$. Thus the following holds:

$$\text{For all } F\colon [n]^k \to N^2 \text{ there exists } E \in S_{p'}[n] \text{ such that}$$
$$F \text{ has } \leq cp' \text{ regressive values on } E^k.$$

We verify that Theorem I holds for $n$, $k$, $p$.

Let $F\colon S_k[n] \to [n]$ be regressive. We apply Theorem II to the function $G\colon [n]^k \to [n]^2$ given by $G(x) = (F(\mathrm{rng}(x)), |\min(x) - 1|)$ if $x$ is strictly increasing; $(0, 0)$ otherwise. (Here $|\ |$ is ordinary one-dimensional, absolute value.)

According to the displayed statement, let $A \subseteq [n]$, $|A| = p'$, be such that $G$ has at most $cp'$ regressive values on $A^k$. For all $b \in A$, let $X(b) = \{x \in A^{k<}\colon x_1 = b\}$. (Here $A^{k<}$ is the set of all strictly increasing elements of $A^k$.) Since the same regressive value of $F$ cannot be obtained from arguments with different minimums, the pigeonhole principle tells us that

$$|\{b \in A\colon G \text{ has } > 2c \text{ regressive values on } X(b)\}| < p'/2\,.$$

Hence $|\{b \in A\colon G$ has $\leq 2c$ regressive values on $X(b)\}| \geq p'/2$. That is, $|S| \geq p'/2$, where $S = \{b \in A\colon F$ has $\leq 2c$ regressive values on elements of $S_k(A)$ whose min is $b\}$.

Now apply the Finite Ramsey Theorem in a well known way to obtain $C \subseteq S$, $|C| = p + 3ck = q$, such that the following holds:

(i) Let $A, B, A', B' \in S_k[C]$, where $(A, B)$ is order isomorphic to $(A', B')$. Then $F(A) < F(B)$ if and only if $F(A') < F(B')$;

(ii) For all $b \in C$, $F$ has $\leq 2c$ regressive values on elements of $S_k(C)$ whose min is $b$. (This follows from $C \subseteq S$.)

Now consider $F(\{C_1, \ldots, C_k\})$ and $F(\{C_1, C_{k+1}, \ldots, C_{2k-1}\})$. If these are distinct, then by shifting the elements above $C_1$, and keeping $C_1$ fixed, and by (i), we get a strictly monotone sequence of values. The number of these values is $> 2c$, and hence we contradict (ii). Therefore, $F(\{C_1, \ldots, C_k\}) = F(\{C_1, C_{k+1}, \ldots, C_{2k-1}\}) = F(\{C_1, C_{q-k+2}, \ldots, C_q\})$.

Let $C' = \{C_1, \ldots, C_p\}$. Then by (i), for all $A \in S_k(C')$, $F(A) = F(\{A_1, C_{q-k+2}, \ldots, C_q\})$. Therefore, $C'$ is as required by Theorem I. □

THEOREM 1.6. *The following are provably equivalent in* PA:

(i) *Theorem* I;
(ii) *Theorem* II;
(iii) *Theorem* III;
(iv) *Theorem* 0.3;
(v) *Theorem* 0.4;
(vi) 1-Con(PA).

*Proof.* (i) ↔ (vi) by Lemma 1.1. (iii) → (i) by Lemma 1.2. (iii) → (iv) → (v) → (ii) by Lemma 1.2. (ii) → (i) by Lemma 1.5. (i) → (iii) by Lemmas 1.3 and 1.4.

Hence (i) → (iii) → (iv) → (v) → (ii) → (i). Therefore (i)–(v) are equivalent. Hence (i)–(vi) are equivalent. □

We now want to consider associated numerical functions. Recall the definition of $H_{\mathrm{I}}(k, p)$ as the least $n$ such that Theorem I holds for $n, k, p$. Let $H_{\mathrm{III}}(k, r, p)$ be the least $n$ such that Theorem III holds for $n, r, k, p$. Let $H_{0.3}(k, r, p)$ be the least $n$ such that Theorem 0.3 holds for $n, k, r, p$. Let $H_{0.4}(k, r, p)$ be the least $n$ such that Theorem 0.4 holds for $n, k, r, p$. Let $H_{\mathrm{I},k}$, $H_{\mathrm{III},k}$, $H_{0.3,k}$, $H_{0.4,k}$ be the respective cross sections. We use $H'$ instead of $H$ if we fix $r = 2$.

LEMMA 1.7. *For each fixed $k$, Theorems* I, III, 0.3, *and* 0.4 *are provable in* PA. *There is a primitive recursive function $h$ such that the following hold. Let $k, r, p > 1$.*

(i) $H_{\mathrm{III}}(k,r,p) \geq H_{0.3}(k,r,p) \geq H_{0.4}(k,r,p)$;
(ii) *if* $n \geq h(k,p)$, *then* $H_{0.4}(k,2,n) \geq H_{\mathrm{I}}(k,p)$.

*Proof.* This results from the proofs of the previous lemmas. The function $h$ is associated with FRT. □

THEOREM 1.8. *The* $H_{\mathrm{I},k}$, $H_{\mathrm{III},k}$, $H_{0.3,k}$, *and* $H_{0.4,k}$ *are all* $< \varepsilon_0$*-recursive functions. Each of the four hierarchies of functions* $H_{\mathrm{I},k}$, $H'_{\mathrm{III},k}$, $H'_{0.3,k}$, *and* $H'_{0.4,k}$ *are cofinal in the* $< \varepsilon_0$*-recursive functions in the sense that every* $< \varepsilon_0$*-recursive function is eventually* $\leq$ *some function in the hierarchy.*

*Proof.* The four functions are provably recursive functions of PA, and hence are $< \varepsilon_0$-recursive functions.

For the second claim, from Lemma 1.7 we see that there exists a primitive recursive function $h$ such that for all $k > 1$ and $n \geq h(k,p)$, we have $H_{\mathrm{III},k}(n)$, $H'_{0.3,k}(n)$, $H'_{0.4,k}(n) \geq H_{\mathrm{I},k}(p)$. By Lemma 1.1, the $H_{\mathrm{I},k}$ are cofinal in the $< \varepsilon_0$-recursive functions.

Without loss of generality we may assume that for each $k \geq 1$, $h(k,p)$ is a strictly increasing function of $p \geq 1$.

Now let $J\colon N \to N$ be a $< \varepsilon_0$-recursive function. Without loss of generality, assume that $J$ is strictly increasing. Let $W\colon N \to N$ be a strictly increasing primitive recursive function such that for all $k$, and all sufficiently large $p$ relative to $k$, $W(p) \geq h(k,p+1)$.

Since $J(W(p))$ is also a $< \varepsilon_0$-recursive function, we can fix $k,t > 1$ such that for all $p \geq t$, $H_{\mathrm{I},k}(p) \geq J(W(p)) \geq J(h(k,p+1))$. Now let $n \geq h(k,t)$. Fix $p \geq t$ such that $h(k,p+1) \geq n \geq h(k,p)$. Then $H'_{\mathrm{III},k}, H'_{0.3,k}(n), H'_{0.4,k}(n) \geq H_{\mathrm{I},k}(p) \geq J(W(p)) \geq J(h(k,p+1)) \geq J(n)$, as required. □

## 2. Large cardinals

The concepts of subtle, almost ineffable, and ineffable cardinals were introduced in an unpublished manuscript of R. Jensen and K. Kunen from 1971, and a number of basic facts were proved there. These concepts were extended to that of $k$-subtle, $k$-alsmost ineffable, and $k$-ineffable cardinals by [Ba75], and a number of results were proved.

These are the large cardinals that we use to prove Propositions A–D in Part 4, and that we show are required to prove Propositions A–D in Part 5. In this part, we give a self-contained presentation of the basic facts about these concepts. We also give some new characterizations of these cardinals in much less set theoretic terms. Some of these are based on linearly ordered sets instead of ordinals. The proofs of our results will appear in [Fr$\infty$].

As is usual in set theory, we treat cardinals and ordinals as von Neumann ordinals. We use $\omega$ for the first limit ordinal, which is also $N$. For sets $X$ and $k \in \omega$, we let $S_k(X)$ be the set of all $k$ element subsets of $X$. Let $S(X)$ be the set of all subsets of $X$. Let $|X|$ be the cardinal of $X$.

Let $\lambda$ be an infinite cardinal. We say that $\lambda$ is regular if and only if every function from an element of $\lambda$ into $\lambda$ is bounded in $\lambda$. We say that $\lambda$ is strongly inaccessible if and only if $\lambda$ is uncountable, regular, and for all cardinals $\kappa < \lambda$, $|S(\kappa)| < \lambda$.

We say that $C \subseteq \lambda$ is closed unbounded if and only if the sup of $C$ is $\lambda$, and for all limit ordinals $x < \lambda$, if the sup of the elements of $C$ below $x$ is $x$, then $x \in C$.

We say that $A \subseteq \lambda$ is stationary if and only if it intersects every closed unbounded subset of $\lambda$.

For any $k$ element set $A$ and $1 \leq i \leq k$, we define $A_i$ to be the $i^{\text{th}}$ least element of $A$. In particular, $\min(A) = A_1$.

We say that $f\colon S_k(\lambda) \to \lambda$ is regressive if and only if for all $A \in S_k(\lambda)$, either $f(A) < \min(A)$ or $\min(A) = 0$.

We say that $f\colon S_k(\lambda) \to S(\lambda)$ is regressive if and only if for all $A \in S_k(\lambda)$, $f(A) \subseteq \min(A)$. We say that $E$ is $f$-homogeneous if and only if $E \subseteq \lambda$ and for all $C, D \in S_k(E)$, we have

$$f(C) \cap \min(C \cup D) = f(D) \cap \min(C \cup D)\,.$$

*Note.* There is a slight abuse of notation here in that every map into $\lambda$ can be viewed as a map into $S(\lambda)$. This will not cause any difficulties.

We now give the three leading definitions from [Ba75]:

Let $k$ be an integer $> 0$. $\lambda$ is $k$-subtle if and only if

(i) $\lambda$ is an infinite cardinal;

(ii) For all closed unbounded $C \subseteq \lambda$ and regressive $f\colon S_k(\lambda) \to S(\lambda)$, there exists an $f$-homogeneous $A \in S_{k+1}(C)$.

$\lambda$ is $k$-almost ineffable if and only if

(i) $\lambda$ is an infnite cardinal;

(ii) For all regressive $f\colon S_k(\lambda) \to S(\lambda)$, there exists an $f$-homogeneous $A \subseteq \lambda$ of cardinality $\lambda$.

$\lambda$ is $k$-ineffable if and only if

(i) $\lambda$ is an infinite cardinal;

(ii) For all regressive $f\colon S_k(\lambda) \to S(\lambda)$, there exists an $f$-homogeneous stationary $A \subseteq \lambda$.

[Ba75] does not define $\lambda$ is 0-subtle, 0-almost ineffable, or 0-ineffable. The appropriate definition is that $\lambda$ is an uncountable regular cardinal.

We begin with the hierarchy theorem.

THEOREM 2.1 ([Ba75]). *Let $k \geq 0$ and $\lambda$ be an infinite cardinal. If $\lambda$ is $k+1$-subtle ($k+1$-almost ineffable, $k+1$-ineffable), then $\lambda$ is $k$-subtle ($k$-almost ineffable, $k$-ineffable). If $\lambda$ is $k+1$-subtle ($k+1$-almost ineffable,*

*$k+1$-ineffable), then there are $\lambda$ many $k$-subtle ($k$-almost ineffable, $k$-ineffable) cardinals below $\lambda$.*

We also have the following relationships between the three concepts.

THEOREM 2.2 ([Ba75]). *Let $k \geq 0$ and $\lambda$ be an infinite cardinal. If $\lambda$ is $k$-ineffable, then $\lambda$ is $k$-almost ineffable. If $\lambda$ is $k$-almost ineffable, then $\lambda$ is $k$-subtle. Let $k \geq 1$. If $\lambda$ is $k$-ineffable, then there are $\lambda$ many $k$-almost ineffable cardinals below $\lambda$. If $\lambda$ is $k$-almost ineffable, then there are $\lambda$ many $k$-subtle cardinals below $\lambda$. If $\lambda$ is $k+1$-subtle, then there are $\lambda$ many $k$-ineffable cardinals below $\lambda$.*

So the weakest of these large cardinal properties is 1-subtle, which is normally just called subtle. From Theorems 2.1 and 2.2, we see that any $k$-subtle, $k$-almost ineffable, or $k$-ineffable cardinal, $k \geq 1$, is subtle. How big is a subtle cardinal?

THEOREM 2.3 ([Ba75]). *Every subtle cardinal is strongly inaccessible.*

Strongly inaccessible cardinals are just beyond the scope of ZFC (they cannot be proved to exist within ZFC). Thus subtle cardinals are also beyond the scope of ZFC.

Actually, much more is true. Let $k > 0$ and $\lambda$ be an infinite cardinal. We say that $\lambda$ is $k$-RP if and only if $\lambda$ is uncountable and for all functions $f\colon S_k(\lambda) \to 2$, $f$ is constant on a subset of $\lambda$ of cardinality $\lambda$. (RP stands for "Ramsey property.")

THEOREM 2.4. *Let $\lambda$ be an infinite cardinal. Then $\lambda$ is* 2-RP *if and only if $\lambda$ is $k$-RP for all $k \geq 1$. If $\lambda$ is* 2-RP, *then $\lambda$ is strongly inaccessible, and there are $\lambda$ many strongly inaccessible cardinals below $\lambda$.*

For proofs and historical discussion, see [Ka94, p. 76].

THEOREM 2.5 ([Ba75]). *Let $\lambda$ be a subtle cardinal. Then $\lambda$ is* 2-RP *and there are $\lambda$ many* 2-RP *cardinals below $\lambda$.*

Thus the $k$-RP does not carry us even up to subtle cardinals. However, we now consider the $k$-SRP. Let $k > 0$ and $\lambda$ be an infinite cardinal. We say that $\lambda$ is $k$-SRP if and only if for all functions $f\colon S_k(\lambda) \to 2$, $f$ is constant on some $S_k(E)$, where $E$ is a stationary subset of $\lambda$. (SRP stands for "stationary Ramsey property.")

THEOREM 2.6 ([Ba75]). *Let $k > 0$ and $\lambda$ be an infinite cardinal. Then $\lambda$ is $(k-1)$-ineffable if and only if $\lambda$ is regular and $k$-*SRP.

We now present analogous definitions in terms of regressive functions $f\colon S_k(\lambda) \to \lambda$. These definitions are simpler, since they involve only the constancy of $f$.

Let $k > 0$. $\lambda$ is $k$-subtle$'$ if and only if

(i) $\lambda$ is an infininte cardinal;

(ii) For all closed, unbounded $C \subseteq \lambda$ and regressive $f\colon S_k(\lambda) \to \lambda$, there exists $A \in S_{k+1}(C)$ such that $f$ is constant on $S_k(A)$.

$\lambda$ is $k$-almost ineffable$'$ if and only if

(i) $\lambda$ is an infinite cardinal;

(ii) For all regressive $f\colon S_k(\lambda) \to \lambda$, there exists an $A \subseteq \lambda$ of cardinality $\lambda$ such that $f$ is constant on $S_k(A)$.

$\lambda$ is $k$-ineffable$'$ if and only if the following hold:

(i) $\lambda$ is an infininte cardinal;

(ii) For all regressive $f\colon S_k(\lambda) \to \lambda$, there exists a stationary $A \subseteq \lambda$ such that $f$ is constant on $S_k(A)$.

Note that $k$-ineffable$'$ immediately implies $k$-SRP.

Theorem 2.7 ([Ba75]). *Let $k \geq 0$ and $\lambda$ be an infinite cardinal. $\lambda$ is $k$-subtle if and only if $\lambda$ is $(k+1)$-subtle$'$; $\lambda$ is $k$-almost ineffable if and only if $\lambda$ is $(k+1)$-almost ineffable$'$; $\lambda$ is $k$-ineffable if and only if $\lambda$ is $(k+1)$-ineffable$'$.*

Actually, there is excess set theoretic baggage in even the three definitions just given as well as in $k$-SRP. We give a number of less set theoretic formulations. These results are due to the author and will appear in [Fr∞].

Let $\alpha$ be an ordinal and $k > 0$. We say that $\alpha$ is purely $k$-subtle if and only if

(i) $\alpha$ is an ordinal;

(ii) For all regressive $f\colon S_k(\alpha) \to \alpha$, there exists $A \in S_{k+1}(\alpha\backslash\{0,1\})$ such that $f$ is constant on $S_k(A)$.

Here \ denotes set theoretic difference. Note that we have already removed the use of closed unbounded sets in this definition.

Observe that this concept is upward closed in the sense that if $\alpha$ is purely $k$-subtle and $\alpha \leq \beta$, then $\beta$ is purely $k$-subtle. This is not true for any of the earlier concepts.

Why do we write $\alpha\backslash\{0,1\}$ instead of $\alpha$ or $\alpha\backslash\{0\}$? Because then $\omega + k$ would be $k$-subtle for all $k > 0$, which is hardly a large cardinal property.

Theorem 2.8. *Let $k > 1$. The least purely $k$-subtle ordinal (if it exists) $=$ the least $k$-subtle$'$ cardinal $=$ the least $(k-1)$-subtle cardinal. The least purely 1-subtle ordinal is $\omega + \omega + 1$.*

We now distill pure $k$-subtlety down even further.

Let $\alpha$ be an ordinal and $k > 0$. We say that $\alpha$ is $k$-large if and only if

(i) $\alpha$ is an ordinal;

(ii) For all regressive $f\colon \alpha^k \to \alpha$, there exist $1 < \beta_1 < \cdots < \beta_{k+1}$ such that $f(\beta_1, \ldots, \beta_k) = f(\beta_2, \ldots, \beta_{k+1})$.

Note that this concept is also upward closed.

Obviously, every purely $k$-subtle ordinal is $k$-large.

THEOREM 2.9. (i) *Let $k > 0$. An ordinal is $k$-large if and only if it is purely $k$-subtle*;

(ii) *The least $k$-large ordinal* (*if it exists*) = *the least purely $k$-subtle ordinal*;

(iii) *Let $k > 1$. The least $k$-large ordinal* (*if it exists*) = *the least purely $k$-subtle ordinal* = *the least $k$-subtle′ cardinal* = *the least* $(k-1)$-*subtle cardinal*;

(iv) *The least* 1-*large ordinal* = *the least purely* 1-*subtle ordinal* = $\omega + \omega + 1$.

Thus the least 2-large ordinal is the least subtle cardinal.

We now go further and remove even the mention of ordinals.

Let $(X, <)$ be a linear ordering with no endpoints, and $k > 0$. We say that $f\colon X^k \to X$ is regressive if and only if it obeys the inequality $f(x) < \min(x)$.

We say that $(X, <)$ is $k$-critical if and only if

$$\text{for all regressive } f\colon X^k \to X, \text{ there exist } b_1 < \cdots < b_{k+1}$$

$$\text{such that } f(b_1, \ldots, b_k) = f(b_2, \ldots, b_{k+1}).$$

THEOREM 2.10. *Let $k > 1$. The cardinality of every $k$-critical linear ordering is a $k$-large cardinal, and a purely $k$-subtle cardinal. The least cardinality of a $k$-critical linear ordering* (*if it exists*) = *the least purely $k$-subtle ordinal* = *the least $k$-subtle′ cardinal* = *the least* $(k-1)$-*subtle cardinal. The least cardinality of a* 1-*critical linear ordering is* $\omega_1$.

We now have a number of formulations of the large cardinal hypothesis used in Part 4 to prove Theorems A–D:

THEOREM 2.11. *The following are provably equivalent within* ZFC.

(i) *for all $k > 0$ there exists a subtle* (*almost ineffable, ineffable*) *cardinal*;

(ii) *for all $k > 0$ there exists a $k$-*SRP *cardinal*;

(iii) *for all $k > 0$ there exists a subtle′* (*almost ineffable′, ineffable′*) *cardinal*;

(iv) *for all $k > 0$ there exists a purely $k$-subtle ordinal*;

(v) *for all $k > 0$ there exists a $k$-large ordinal*;

(vi) *for all $k > 0$ there exists a $k$-critical linear ordering*;

(vii) *there exists an ordinal which is purely $k$-subtle for all $k > 0$*;

(viii) *there exists an ordinal which is $k$-large for all $k > 0$*;

(ix) *there exists a linear ordering which is $k$-critical for all $k > 0$.*

## 3. Function assignments

We present a thorough discussion of the conditions on function assignments first shown in the introduction. We then give some elementary implications between the different forms of Propositions A–D from the introduction.

Recall the following introductory concepts: function assignment, #-decreasing, lex #-decreasing, $\leq_c$, strict order, upward order, and $<_1, <_2$-#-decreasing.

In the statement of all lemmas and definitions in this section, we let $k > 0$, $X$ be a set, $U$ be a function assignment for $X^k$, $<$, $<_1$, $<_2$ be strict orders on $X^k$, $x \in X^k$, and $A$, $B$, $C$ be finite subsets of $X^k$, $f\colon A \to A$, $g\colon B \to B$, $h\colon C \to C$.

When we state theorems, we will reintroduce the objects needed for that theorem and not rely on the data in the previous paragraph.

We write $A \subseteq_< B$ if and only if $A \subseteq B$, and $A$ is downward closed in $B$; i.e., $A \subseteq B$ and for all $x \in A$ and $y \in B$, if $y < x$, then $y \in A$.

We say that $x$ is a $<$-minimal element of $A$ if and only if $x \in A$, and for no $y \in A$ is $y < x$.

LEMMA 3.1. *If $A$ is nonempty, then $A$ has a $<$-minimal element. $B$ can be written as $\{x_1, \ldots, x_n\}$ without repetition, where for no $i < j$ is it the case that $x_i > x_j$. The relation $\subseteq_<$ is transitive. $A \subseteq_< B$ if and only if $A \subseteq B$ and for all $x \in B$, $A \subseteq_< A \cup \{x\}$.*

*Proof.* Let $A$ be nonempty. Define a sequence $x_1, x_2, \ldots$, where $x_1 \in A$ and $x_{i+1}$ is chosen to be an element of $A$ that is $< x_i$, where the sequence terminates if it is no longer possible to continue. By the transitivity and irreflexivity of $<$, we see that no element repeats. Since $A$ is finite, we must terminate. The last element is obviously a $<$-minimal element of $A$.

Assume $B$ is nonempty and define $x_1$ to be any minimal element of $B$. In general, define $x_{i+1}$ to be any minimal element of $B \backslash \{x_1, \ldots, x_i\}$, provided this difference is nonempty. Clearly this is an enumeration of $B$ without repetition. Now let $i < j$. then $x_j \in B \backslash \{x_1, \ldots, x_{i-1}\}$. Since $x_i$ is a minimal element of $B \backslash \{x_1, \ldots, x_{i-1}\}$, we see that $x_i > x_j$ is impossible.

Suppose $A \subseteq_< B$ and $B \subseteq_< C$. Then $A \subseteq C$. Now let $x \in A$ and $y \in C$ and $y < x$. Then $x \in B$, and so $y \in B$. Hence $y \in A$.

Suppose $A \subseteq_< B$ and $x \in B$. We claim $A \subseteq_< A \cup \{x\}$, since this is an obvious weakening of $A \subseteq_< B$. On the other hand, suppose $A \subseteq B$ and for all $x \in B \backslash A$, $A \subseteq_< A \cup \{x\}$. Let $y \in A$, $z \in B$, and $z < y$. Since $A \subseteq_< A \cup \{z\}$, we have $z \in A$. This completes the proof. □

We say that $U$ is singly $<$-end-preserving if and only if for all finite $A \subseteq X^k$ and $x \in X^k$, if $A \subseteq_< A \cup \{x\}$, then $U(A) \subseteq U(A \cup \{x\})$.

We say that $U$ is $<$-end-preserving if and only if for all finite $A \subseteq_< B \subseteq X^k$, $U(A) \subseteq U(B)$.

LEMMA 3.2. *$U$ is $<$-end-preserving if and only if $U$ is singly $<$-end-preserving.*

*Proof.* Let $A \subseteq_< B \subseteq N^k$ be finite. Write $B \backslash A = \{x_1, \ldots, x_n\}$, where for all $1 \leq i \leq n$, $A \subseteq_< A \cup \{x_i\}$, and for all $1 \leq i < j \leq n$, $\neg x_i > x_j$. Then for all $1 \leq i \leq n$, $A \cup \{x_1, \ldots, x_{i-1}\} \subseteq_< A \cup \{x_1, \ldots, x_i\}$. Since $U$ is singly $<$-end-preserving, we see that for all $1 \leq i \leq n$, $U(A \cup \{x_1, \ldots, x_{i-1}\}) \subseteq U(A \cup \{x_1, \ldots, x_i\})$. Therefore $U(A) \subseteq U(B)$ as required. □

We write $f\ (<_1, <_2)^*\ g$ if and only if for all $x \in A \cap B$,

$$\text{if for all } y \in A,\ y <_1 x \text{ implies } f(y) = g(y),$$
$$\text{then } f(x) = g(x) \text{ or } f(x) >_2 g(x).$$

We say that $U$ is $<_1, <_2 -*$-decreasing if and only if for all finite $A, B \subseteq X^k$, $U(A)\ (<_1, <_2)^*\ U(B)$.

We say that $U$ is singly $<_1, <_2 -*$-decreasing if and only if for all finite $A \subseteq X^k$ and $x \in X^k$, we have $U(A)\ (<_1, <_2)^*\ U(A \cup \{x\})$.

LEMMA 3.3. *If $U$ is $<_1, <_2 -\#$-decreasing, then $U$ is $<_1$-end-preserving.*

*Proof.* Let $U$ be as given. By Lemma 3.2, it suffices to show that $U$ is singly $<_1$-end-preserving. Accordingly, assume that $A \subseteq_{<_1} A \cup \{x\}$, where $x \notin A$. Then $x$ is not $<_1$ any element of $A$. Hence clearly $U(A) \subseteq U(A \cup \{x\})$. □

LEMMA 3.4. *Suppose $U$ is $<$-end-preserving, $A \subseteq B$, $y \in A$, and no element of $B \backslash A$ is $< y$. Then $U(A)(y) = U(B)(y)$.*

*Proof.* Let $C = \{z \in B : z < y \text{ or } z = y\}$. Then $C \subseteq_< B$. Also we claim $C = \{z \in A : z < y \text{ or } z = y\}$. To see this, let $z \in C$. Now $z \in B \backslash A$ is impossible by hypothesis. Hence $z \in A$.

We therefore have $C \subseteq_< A$. By $<$-end-preserving, we have $U(C) \subseteq U(A)$ and $U(C) \subseteq U(B)$. Since $y \in A$, $B$, and $C$, we have $U(A)(y) = U(B)(y)$. □

LEMMA 3.5. *If $U$ is $<_1, <_2 -\#$-decreasing, then $U$ is singly $<_1, <_2 -*$-decreasing.*

*Proof.* Let $U$ be $<_1, <_2 -\#$-decreasing. We must show that $U(A)$ $(<_1, <_2)^*\ U(A \cup \{x\})$.

Let $y \in A$ be such that for all $z \in A$ with $z <_1 y$, we have $U(A)(z) = U(A \cup \{x\})(z)$. We must show that $U(A)(y) >_2 U(A \cup \{x\})(y)$ or $U(A)(y) = U(A \cup \{x\})(y)$. By Lemma 3.3, $U$ is $<_1$-end-preserving.

According to Lemma 3.4, we can assume, without loss of generality, that some element of $A \cup \{x\} \backslash A$ is $<_1 y$. That is, we can assume, without loss of generality, that $x <_1 y$.

Let $B = \{z \in A \colon z <_1 y \text{ or } z = y\}$. Then $B \subseteq_{<_1} A$. We claim that $B \cup \{x\} \subseteq_{<_1} A \cup \{x\}$. To see this, let $w \in A \backslash B \backslash \{x\}$, $w <_1 z$, $z \in B \cup \{x\}$. If $z \in B$, then this is a contradiction. However, if $z = x$, then $w <_1 z = x <_1 y$, and so $w \in B$, which is again a contradiction.

Since $U$ is $<_1, <_2$ $-\#$-decreasing, we obtain the disjunction $U(B) \subseteq U(B \cup \{x\})$ or there exists $z > x$ such that $U(B)(z) >_2 U(B \cup \{x\})(z)$. Since $U$ is $<_1$-end-preserving, we see that $U(B) \subseteq U(A)$ and $U(B \cup \{x\}) \subseteq U(A \cup \{x\})$. If the first disjunct holds, then $U(B)(y) = U(B \cup \{x\})(y) = U(A)(y) = U(A \cup \{x\})(y)$, and we are done. If the second disjunct holds, then fix $z$ such that $U(B)(z) >_2 U(B \cup \{x\})(z)$. If $z = y$, then we are done. Assume $z \neq y$. Then since $z \in B$, we have $z <_1 y$ and $z \in A$. By the choice of $y$, we have $U(A)(z) = U(A \cup \{x\})(z)$, which is a contradiction. This completes the proof. □

Lemma 3.6. *Suppose $A \subseteq B \subseteq C$, $f$ $(<_1, <_2)^*$ $g$ and $g$ $(<_1, <_2)^*$ $h$. Furthermore, suppose that $g$ and $h$ agree on $B \backslash A$. Then $f$ $(<_1, <_2)^*$ $h$.*

*Proof.* (See notation from the third paragraph of this section.) Let $x \in A$ be such that for all $y \in A$ with $y <_1 x$, we have $f(y) = h(y)$.

First suppose that for all $y \in A$ with $y <_1 x$, we have $f(y) = g(y)$. Then for all $y \in A$ with $y <_1 x$, we have $g(y) = h(y)$. Hence for all $y \in B$ with $y <_1 x$, we have $g(y) = h(y)$. Hence $g(x) = h(x)$ or $g(x) >_2 h(x)$. But also $f(x) = g(x)$ or $f(x) >_2 g(x)$. Hence $f(x) = h(x)$ or $f(x) >_2 h(x)$ as required.

Now assume that for some $y \in A$ with $y <_1 x$, $f(y) \neq g(y)$. Fix $z$ to be a $<_1$-minimal element of $A$ such that $z <_1 x$, $f(z) \neq g(z)$. Then $f(z) >_2 g(z)$. By minimality, for all $w \in A$ with $w <_1 z$, we have $f(w) = g(w)$. Hence for all $w \in B$ with $w <_1 z$, we have $g(w) = h(w)$. Hence $g(z) = h(z)$ or $g(z) >_2 h(z)$. Therefore, $f(z) >_2 h(z)$. But since $z <_1 x$, we have $f(z) = h(z)$. This is a contradiction. □

Lemma 3.7. *If $U$ is $<_1, <_2$ $-\#$-decreasing and $A \subseteq B$, then $U(A)$ $(<_1, <_2)^*$ $U(B)$.*

*Proof.* By Lemma 3.1, write $B \backslash A = \{x_1, \ldots, x_n\}$, where for no $i < j$ is $x_i >_1 x_j$. We show by induction on $i$ that $U(A)$ $(<_1, <_2)^* U(A \cup \{x_1, \ldots, x_i\})$. By Lemma 3.5, this is true for $i = 1$. Suppose this is true for $1 \leq i < n$. By Lemma 3.5, $U(A \cup \{x_1, \ldots, x_i\})$ $(<_1, <_2)^*$ $U(A \cup \{x_1, \ldots, x_{i+1}\})$. Now we claim that $U(A \cup \{x_1, \ldots, x_i\})$ and $U(A \cup \{x_1, \ldots, x_{i+1}\})$ agree on $\{x_1, \ldots, x_i\}$. This follows from Lemma 3.4 since $x_{i+1}$ is not $<_1$ any element of $\{x_1, \ldots, x_i\}$.

Hence we can apply Lemma 3.6 to obtain $U(A)\ (<_1,<_2)^*\ U(A\cup\{x_1,\ldots,x_{i+1}\})$ as required. ☐

LEMMA 3.8. *If $U$ is $<_1,<_2$ $-\#$-decreasing, then $U$ is $<_1,<_2$ $-*$-decreasing.*

*Proof.* Let $A, B \subseteq N^k$ be finite, $x \in A\cap B$, and for all $y \in A$, $y <_1 x$ implies $U(A)(y) = U(B)(y)$. Let $A' = \{z \in A : z <_1 x \text{ or } z = x\}$, and let $B' = \{z \in B : z <_1 x \text{ or } z = x\}$. Then $A' \subseteq B'$. By Lemma 3.7, $U(A')(<_1,<_2)^*$ $U(B')$. Note that $A' \subseteq_{<_1} A$ and $B' \subseteq_{<_1} B$. Hence by Lemma 3.3, $U(A') \subseteq U(A)$ and $U(B') \subseteq U(B)$. Therefore, $x \in A' \cap B'$ and for all $y \in A'$, $y <_1 x$ implies $U(A')(y) = U(B')(y)$. Hence $U(A')(x) = U(B')(x)$ or $U(A')(x) >_2 U(B')(x)$. Therefore, $U(A)(x) = U(B)(x)$ or $U(A)(x) >_2 U(B)(x)$. ☐

LEMMA 3.9. *If $U$ is $<_1,<_2$ $-*$-decreasing, then $U$ is $<_1$-end-preserving.*

*Proof.* Let $A \subseteq_{<_1} A\cup\{x\}$. By Lemma 3.2 it suffices to show that $U(A) \subseteq U(A\cup\{x\})$. Assume this is false, and let $y$ be a $<_1$-minimal element of $A$ such that $U(A\cup\{x\})(y) \neq U(A)(y)$. Now $x \notin A$ and $\neg x <_1 y$. Therefore,

(i) $y \in A \cap (A\cup\{x\})$ and for all $z \in A$, $z <_1 y$ implies $U(A)(z) = U(A\cup\{x\})(z)$;

(ii) $y \in (A\cup\{x\})\cap A$ and for all $z \in A\cup\{x\}$, $z <_1 y$ implies $U(A\cup\{x\})(z) = U(A)(z)$.

Hence by $<_1,<_2$ $-*$-decreasing we see that

(iii) $U(A)(y) = U(A\cup\{x\})(y)$ or $U(A)(y) >_2 U(A\cup\{x\})(y)$;

(iv) $U(A\cup\{x\})(y) = U(A)(y)$ or $U(A\cup\{x\})(y) >_2 U(A)(y)$.

Hence, $U(A)(y) = U(A\cup\{x\})(y)$, which is a contradiction. ☐

THEOREM 3.10. *Let $k > 0$, $X$ be a set, $<_1$, $<_2$ be strict orders on $X^k$, and $U$ be a function assignment for $X^k$. Then $U$ is $<_1,<_2$ $-\#$-decreasing if and only if $U$ is $<_1,<_2$ $-*$-decreasing. If $U$ is $<_1,<_2$ $-\#$-decreasing, then $U$ is $<_1$-end-preserving.*

*Proof.* The forward direction is by Lemma 3.8. Now suppose that $U$ is $<_1,<_2$ $-*$-decreasing. Let $A \subseteq X^k$ be finite and $x \in X^k$. Assume that $\neg U(A) \subseteq U(A\cup\{x\})$. We need to find $y >_1 x$ such that $U(A)(y) >_2 U(A\cup\{x\})(y)$.

Let $y$ be a $<_1$ minimal element of $A$ such that $U(A)(y) \neq U(A\cup\{x\})(y)$. Then $U(A)(y) = U(A\cup\{x\})(y)$ or $U(A)(y) >_2 U(A\cup\{x\})(y)$. Hence $U(A)(y) >_2 U(A\cup\{x\})(y)$.

Now by Lemma 3.9, $U$ is $<_1$-end-preserving. Hence by Lemma 3.4, we see that $x <_1 y$. This completes the proof of the reverse direction. The final statement is from Lemmas 3.8 and 3.9. ☐

We now consider the various forms of Propositions A–D.

LEMMA 3.11. (i) *Proposition* A *for #-decreasing is equivalent to Proposition* B *for #-decreasing*;

(ii) *Proposition* A *for* lex *#-decreasing is equivalent to Proposition* B *for* lex *#-decreasing*;

(iii) *Proposition* C *is equivalent to Proposition* D;

(iv) *Proposition* C *implies Proposition* A *for #-decreasing and for* lex *#-decreasing.*

*Proof.* The reverse implications in i), ii), iii) are immediate since if $U$ is a function assignment for $N^k$, then the restriction of $U$ to subsets of any $[n]^k$ is a function assignment for $[n]^k$.

The forward implications in (i), (ii), (iii) are all proved by a finitely branching tree argument, where we fix $k$, $p$, and assume that the antecedent is true for $k$, $p$, and that the consequent is false. We consider the infinite finitely branching tree of function assignments (and upward orderings) for the $[n]^k$ which are counterexamples. There must be an infinite path, which yields an appropriate function assignment (and upward orderings), to which we can apply the antecedent to obtain the required contradiction.

For (iv), note that Proposition A for #-decreasing is equivalent to Proposition C for $<_1=<_2=$ the ordering given by $|x| < |y|$. Also Proposition A for lex #-decreasing is equivalent to Proposition C for $<_1=<_2=<_{\text{lex}}$. ☐

We now wish to show that Proposition A for lex #-decreasing implies Proposition A for #-decreasing. Before dealing directly with function assignments, we need the following lemma:

LEMMA 3.12. *Let $c, k > 0$ and $U$ be a function assignment for $N^k$. Suppose that for all $p > 0$, some $U(A)$ has at most $c$ regressive values on some $E^k \subseteq A$, $|E| = p$. Then for all $p > 0$, some $U(A)$ has at most $k^k$ regressive values on some $E^k \subseteq A$, $|E| = p$. In fact, $k^k$ can be replaced by* $\text{ot}(k)$.

*Proof.* Let $c$, $k$, $U$ be as given. Let $p > 0$. Choose $q \gg c, k, p$. Let $U(A)$ have at most $c$ regressive values on $E^k \subseteq A$, $|E| = q$. Let $t \in N^k$ not be a regressive value of $U(A)$ on $E^k$. Define $g\colon E^k \to N$ by $g(x) = U(A)(x)$ if $|U(A)(x)| < \min(x)$; $t$ otherwise. By the Finite Ramsey Theorem, let $E' \subseteq E$, $|E'| = p$, be homogeneous for $g$. Then $U(A)$ has at most $\text{ot}(k)$ regressive values on $E'^k$, as required. ☐

Let $U$ be a function assignment for $N^k$. We say that $U$ is $*$-decreasing if and only if $U$ is $<_1, <_2 -*$-decreasing, where $<_1=<_2$ is given by $|x| < |y|$. By Theorem 3.10, $U$ is #-decreasing if and only if $U$ is $*$-decreasing.

LEMMA 3.13. *Proposition* A *for* lex #*-decreasing implies Proposition* A *for #-decreasing.*

*Proof.* Assume Proposition A for lex #-decreasing. Let $U$ be a #-decreasing function assignment for $N^k$. We now construct a lex #-decreasing function assignment $V$ for $N^{k+1}$ as follows. Let $A \subseteq N^{k+1}$ be finite. Let $A' = \{x \in N^k: (|x|, x) \in A\}$. Define $V(A): A \to A$ as follows. For $(|x|, x) \in A$, let $V(A)(|x|, x) = (|U(A')(x)|, U(A')(x))$. For the remaining $y \in A$, let $V(A)(y) = y$.

We now show that $V$ is lex #-decreasing. By Theorem 3.10, it suffices to prove that $V$ is lex $*$-decreasing; i.e., $<_{\text{lex}}, <_{\text{lex}}$ $-*$-decreasing. To this end, let $A, B \subseteq N^{k+1}$ be finite, $y \in A \cap B$, and

$$\text{for all } z \in A \text{ with } z <_{\text{lex}} y,\ V(A)(z) = V(B)(z),$$
$$\text{and } V(A)(y) <_{\text{lex}} V(B)(y).$$

Clearly $y$ is of the form $(|x|, x)$. Hence $V(A)(y) = (|U(A')(x)|, U(A')(x))$ and $V(B)(y) = (|U(B')(x)|, U(B')(x))$. Therefore, $U(A')(x) \neq U(B')(x)$.

We now claim that for all $w \in A'$ with $|w| < |x|$, $U(A')(w) = U(B')(w)$. Suppose $|w| < |x|$ and $U(A')(w) \neq U(B')(w)$. Then

$$V(A)((|w|, w)) \neq V(B)((|w|, w)),$$

which contradicts the previous paragraph. Also $x \in A' \cap B'$ (since $y \in A \cap B$).

Now since $U$ is #-decreasing, $U$ is also $*$-decreasing, and hence $|U(A')(x)| > |U(B')(x)|$. Therefore $V(A)((|x|, x)) >_{\text{lex}} V(B)((|x|, x))$, which is the desired contradiction.

We have now shown that $V$ is a lex $*$-decreasing function assignment for $N^{k+1}$. Then for all $p > 0$, some $V(A)$ has at most $(k+1)^{k+1}$ regressive values on some $E^{k+1}$, $|E| = p$. Let $p > 0$ and fix such $A, E$.

We claim that if $y$ is a regressive value of $U(A')$ on $E^k$, then $(|y|, y)$ is a regressive value of $V(A)$ on $E^{k+1}$. To see this, let $U(A')(w) = y$, $|y| < \min(w)$, $w \in E^k$. Then $V(A)((|w|, w)) = (|y|, y)$, and $(|w|, w) \in E^{k+1}$.

From this we conclude that the number of regressive values of $U(A')$ on $E^k$ is at most the number of regressive values of $V(A)$ on $E^{k+1}$. Hence $U(A')$ has at most $(k+1)^{k+1}$ regressive values on some $E^{k+1} \subseteq A'$, $|E| = p$.

Since $p$ is arbitrary, we can now apply Lemma 3.12 with $c = (k+1)^{k+1}$ to obtain Proposition A for #-decreasing (even with $k^k$ replaced by $\text{ot}(k)$). □

THEOREM 3.14. *Proposition* C *implies all forms of Propositions* A–D. *Proposition* A *for #-decreasing is implied by all forms of Propositions* A–D. *All forms of Proposition* A–D *remain equivalent if* $k^k$ *is replaced by* $\text{ot}(k)$.

*Proof.* The first claim is from Lemma 3.11. The second claim follows from Lemmas 3.11 and 3.13. The third claim for Propositions A, C follows from Lemma 3.12. Then derive the appropriate equivalences of A with B, and C with D, as in Lemma 3.11 (i), (ii). □

## 4. Proofs using large cardinals

The goal of this part is to prove all forms of Propositions A–D from ZFC + "there exists a suitably large cardinal." According to Theorem 3.14 we need only prove Proposition C.

For clarity, we first prove Proposition A for #-decreasing function assignments.

We make heavy use of Theorem 3.10. In particular, we use that if $U$ is #-decreasing, then $U$ is $*$-decreasing and end-preserving. These concepts are of course in the context of underlying set $X = (N, <)$, where $<_1 = <_2$ is the strict order on $N^k$ given by $|x| < |y|$.

In the treatment of function assignments for $X^k$ in Part 3, we did not use a linear ordering of $X$. In fact, if $U$ is a function assignment for $X^k$, then we can always recover the set $X$ from $U$ by setting $X$ to be the union of the domains of the elements of $X$.

In this section, we will make use of a linear ordering on $X$. However, in so doing, we will not be able to recover the linear ordering of $X$ from the function assignment. Consequently in this section, we use function assignments with additional structure as follows.

A linearly ordered function assignment is a triple $(X^k, \leq, U)$, where $\leq$ is a linear ordering on the set $X$, and $U$ is a function assignment for $X^k$. For $x \in X^k$, we write $|x|$ for the maximum coordinate of $x$.

We fix two linearly ordered function assignments $(X^k, \leq, U)$ and $(Y^k, \leq', V)$.

We say that a partial function $f: X \to Y$ is order-preserving if and only if for all $a, b \in X$, $a \leq b$ if and only if $f(a) \leq' f(b)$.

Let $h: A \to B$, where $A \subseteq X$ and $B \subseteq Y$. Then $h$ naturally extends to a function from $A^k$ into $B^k$ by $h(x_1, \ldots, x_k) = (h(x_1), \ldots, h(x_k))$.

Now assume that $A$, $B$ are finite subsets of $X^k, Y^k$, respectively. We write $\text{fld}(A)$ for the set of all coordinates of elements of $A$. We say that $h$ is an order isomorphism from $A$ onto $B$ if and only if $h$ is an order-preserving bijection which maps $\text{fld}(A)$ onto $\text{fld}(B)$ such that $h$ maps $A$ onto $B$. We say that $A$, $B$ are order isomorphic if and only if there exists an order isomorphism from $A$ onto $B$. If an order isomorphism exists, then it is unique.

We say that $f$ is a finite endomorphism in $X^k$ if and only if $\text{dom}(f)$ is a finite subset of $X^k$ and $\text{rng}(f) \subseteq \text{dom}(f)$. We say that a finite endormorphism $f$ in $X^k$ is order isomorphic to a finite endormorphism $g$ in $Y^k$ if and only if their graphs, as subsets of $X^{2k}$ and of $Y^{2k}$, are order isomorphic. Clearly, if $f$, $g$ are order isomorphic, then the unique order isomorphism must be the order isomorphism between their domains.

We say that $h$ is an order isomorphism from $(X^k, \leq, U)$ onto $(Y^k, \leq', V)$ if and only if

(i) $h$ is an order isomorphism from $(X, \leq)$ onto $(Y, \leq')$;

(ii) For all $A \subseteq X^k$, $U(A)$ is order isomorphic to $V(h[A])$ by $h|\mathrm{fld}\,(A)$.

We say that $(X^k, \leq, U)$ and $(Y^k, \leq', V)$ are order isomorphic if and only if there is an order isomorphism from $(X^k, \leq, U)$ onto $(Y^k, \leq', V)$.

Let $S \subseteq X$. We write $(X^k, \leq, U) \mid S$ for $(S^k, \leq', U')$, where $\leq'$ is the intersection of $\leq$ with $S^2$, and $U'$ is the restriction of $U$ to finite subsets of $S^k$.

Let $p \in N$ and $E \subseteq X$. We say that $(X^k, \leq, U)$ is $p$-uniform on $E \subseteq X$ if and only if for all $S, T \subseteq E$ of the same cardinality $\leq p$, $(X^k, \leq, U)|S$ is order isomorphic to $(X^k, \leq, U) \mid T$. We say that $(X^k, \leq, U)$ is $p$-uniform if $(X^k, \leq, U)$ is $p$-uniform on $X$.

We say that $(X^k, \leq, U)$ is uniform if and only if it is $p$-uniform for all finite $p$. We always strive to obtain uniformity.

We say that $(Y^k, \leq', V)$ is finitely conservative over $(X^k, \leq, U)$ if and only if for all finite $S \subseteq Y$ there exists finite $T \subseteq X$ such that $(Y^k, \leq', V)|S$ is order isomorphic to $(X^k, \leq, U)|T$.

We say that $(X^k, \leq, U)$ and $(Y^k, \leq', V)$ are finitely equivalent if and only if they are finitely conservative over each other.

Obviously, if $(Y^k, \leq', V)$ is finitely conservative over $(X^k, \leq, U)$, then every $V(A)$ is order isomorphic to some $U(B)$; however, the converse seems to fail. Also, if $(X^k, \leq, U)$ and $(Y^k, \leq', V)$ are finitely equivalent, then every $U(A)$ is order isomorphic to some $V(B)$ and every $V(A)$ is order isomorphic to some $U(B)$; the converse also seems to fail.

An infinite linearly ordered function assignment is a linearly ordered function assignment whose first component, $X^k$, is infinite.

LEMMA 4.1. *Let $U^* = (X^k, \leq, U)$ be an infinite linearly ordered function assignment. Let $p \in N$, and $E \subseteq X$ be infinite. Then there is an infinite $E' \subseteq E$ such that $U^*$ is $p$-uniform on $E$.*

*Proof.* This is an immediate consequence of Ramsey's theorem once we check that for all $i \leq p$, the equivalence relation "$U^*|S$ order isomorphic to $U^*|T$" for $|S| = |T| = i$ has finitely many equivalence classes. But we can give a concrete representation of the equivalence classes as follows. Given $U^*|S$, we can list the elements of $S$ in increasing order, and define the appropriate finite set of endomorphisms using indices as domain and range elements. Clearly $U^*|S$ is order isomorphic to $U^*|T$ if and only if the corresponding set of functions for $U^*|S$ is the same as the corresponding set of function for $U^*|T$. □

THEOREM 4.2. *Let $k > 0$ and $U^*$ be an infinite linearly ordered function assignment. Then there is an infinite uniform linearly ordered function assignment $V^*$ which is finitely conservative over $U^*$. In fact, $V^*$ can be taken to be of the form $(N^k, \leq, V)$, where $\leq$ is the usual linear ordering of $N$.*

*Proof.* Let $U^* = (X^k, \leq, U)$. We define infinite sets $X = E_0 \supseteq E_1 \supseteq E_2 \ldots$ as follows. For $i > 0$, $E_i \subseteq E_{i-1}$ is chosen so that $U^*$ is $i$-uniform on $E_i$.

We now define the function assignment $V^*$ for $N^k$. Let $A \subseteq N^k$ be finite. We wish to define $V(A)$. Let $S = \text{fld}\,(A)$, $|S| = p$. Let $T \subseteq E_p$, $|T| = p$. Let $W^*$ be the unique linearly ordered function assignment $(S^k, \leq, W)$ which is order isomorphic to $U^*|T$. Take $V(A) = W(A)$. Define $V^* = (N^k, \leq, V)$.

We need to verify that $V^*$ is well defined, i.e., is independent of the choice of $T$. But this is clear from the $p$-uniformity of $U^*$ on $E_p$.

We now verify that $V^*$ is uniform. Let $S, T \subseteq N$ have cardinality $p > 0$. Let $h$ be the unique order isomorphism from $S$ onto $T$. We need to verify that $h$ is an order isomorphism from $V^*|S$ onto $V^*|T$. Let $A \subseteq S^k$. It suffices to show that $h$ is an order isomorphism from $V(A)$ onto $V(h[A])$. But for this, it suffices to know that $V(A)$ and $V(h[A])$ are order isomorphic. But this is clear since the construction of $V^*$ is well defined.

Now let $S \subseteq N$, $|\text{fld}\,(S)| = p$. Let $|T| = p$, $T \subseteq E_p$. We claim that $V^*|S$ is order isomorphic to $U^*|T$. To see this, let $h\colon S \to T$ be order-preserving, and let $A \subseteq S^k$, $|\text{fld}\,(A)| = q$. We must show that $V(A)$ is order isomorphic to $U(h[A])$. Now by construction, $V(A)$ is order isomorphic to some $U(B)$, $B \subseteq E_q^k$, $|\text{fld}\,(B)| = q$. Since $h[A] \subseteq E_p^k \subseteq E_q^k$, $|\text{fld}\,(h[A])| = q$, we see that $U(h[A])$ is order isomorphic to $U(B)$ and $V(A)$. □

Let $U^* = (X^k, \leq, U)$ be a linearly ordered function assignment. We say that $U^*$ is #-decreasing if and only if for all finite $A \subseteq X^k$ and $x \in X^k$,

$$\text{either } U(A) = U(A \cup \{x\}), \text{ or there exists } |y| > |x|$$
$$\text{such that } |U(A)(y)| > |U(A \cup \{x\})(y)|.$$

We say that $U^*$ is $*$-decreasing if and only if for all finite $A, B \subseteq X^k$ and $x \in A \cap B$,

$$\text{if for all } y \in A,\ |y| < |x| \text{ implies } U(A)(y) = U(B)(y),$$
$$\text{then } U(A)(x) = U(B)(x) \text{ or } |U(A)(x)| > |U(B)(x)|.$$

In the above, $<$ refers to $\leq$ on $X$, and $|\ |$ is the sup, which also refers to $\leq$ on $X$. Note that if $X = N$ and $\leq$ is the usual linear ordering on $N$, then $U^*$ is #-decreasing ($*$-decreasing) if and only if $U$ is #-decreasing ($*$-decreasing).

LEMMA 4.3. *Let $U^*$ and $V^*$ be linearly ordered function assignments, where $V^*$ is finitely conservative over $U^*$. If $U$ is #-decreasing, then $V$ is #-decreasing.*

*Proof.* Note that $U^*$ or $V^*$ is #-decreasing if and only if every finite $U^*|S$ or $V^*|S$ is #-decreasing, and #-decreasing is an order theoretic property. □

LEMMA 4.4. *Let $U^*$ be an infinite #-decreasing linearly ordered function assignment. Then there is a uniform infinite #-decreasing linearly ordered function assignment $V^*$ which is finitely conservative over $U^*$. Furthermore, $V^*$ can be taken to be of the form $(N^k, \leq, V)$, where $\leq$ is the usual linear ordering of $N$.*

*Proof.* This is immediate from Theorem 4.2 and Lemma 4.3. □

We say that $(X^k, \leq, U)$ is order invariant if and only if for all finite order isomorphic $A, B \subseteq X^k$, $U(A)$ is order isomorphic to $U(B)$.

LEMMA 4.5. *A linearly ordered function assignment is uniform if and only if it is order invariant.*

*Proof.* Let $U^* = (X^k, \leq, U)$ be uniform. Let $A, B \subseteq X^k$ be finite and order isomorphic, with respective fields $S$, $T$. Now $U^* \mid S$ and $U^* \mid T$ are order isomorphic. The unique order isomorphism must send $A$ to $B$. Hence $U(A)$ and $U(B)$ are order isomorphic.

Now let $U^* = (X^k, \leq, U)$ be order invariant. Let $S$, $T$ be finite subsets of $X$ of the same cardinality. Let $h$ be the unique order-preserving bijection from $S$ onto $T$. Let $A \subseteq S^k$. We must show that $U(A)$ is order isomorphic to $U(h[A])$. But this is immediate since $A$ is order isomorphic to $h[A]$. □

LEMMA 4.6. *Let $U^*$ be a uniform infinite linearly ordered function assignment, and let $(Y, \leq)$ be an infinite linear ordering. Then there is a unique linearly ordered function assignment $V^* = (Y^k, \leq, V)$ which is finitely conservative over $U^*$. This unique $V^*$ is itself uniform and is finitely equivalent to $U^*$. Furthermore, $U^*$ is #-decreasing if and only if $V^*$ is #-decreasing.*

*Proof.* For uniqueness, let $A \subseteq Y^k$ be finite. Then $V(A)$ must be order isomorphic to some $U(B)$. Since all $U(B)$, where $B$ is isomorphic to $A$, are order isomorphic, $V(A)$ is uniquely determined up to order isomorphism. But then $V(A)$ is completely determined.

The previous paragraph actually defines $V$ and $V^*$ since $U$ is infinite. By construction, $V^*$ is finitely equivalent to $U^*$. It is also clear that $V^*$ is order invariant. Therefore by Lemma 4.5, $V^*$ is uniform. The remainder follows from the fact that $U^*$ $(V^*)$ is #-decreasing if and only if for all finite $S \subseteq X(Y)$, $U^*|S$ $(V^*|S)$ is #-decreasing. □

We say that $f$ is a direct completion of $U^*$ if and only if

(i) $f\colon X^k \to X^k$;

(ii) For all finite $A \subseteq X^k$ there exists $B \supseteq A$ such that $U(B) \subseteq f$.

It is interesting to note that this definition depends only on $X$ and $U$, and not on $\leq$.

A well ordered function assignment is a linearly ordered function assignment whose $\leq$ is a well ordering.

We now wish to prove the fundamental fact that every #-decreasing well ordered function assignment has a direct completion. To this end, fix a well ordered function assignment $V^* = (X^k, \leq, V)$ which is #-decreasing.

We say that $f$ is special if and only if

(i) $f$ is a partial function from $X^k$ into $X^k$ such that for all $y \in \operatorname{dom}(f)$ and $|x| < |y|$, we have $x \in \operatorname{dom}(f)$;

(ii) For all finite $A \subseteq \operatorname{dom}(f)$ there exists $B \supseteq A$ such that $U(B) \subseteq f$.

For $A \subseteq X^k$, we say that $x$ is minimal over $A$ if and only if $x \in X^k \backslash A$, and for all $|y| < |x|$, $y \in A$.

Lemma 4.7. *Let $f$ be special and $x \in X^k$ be minimal over* $\operatorname{dom}(f)$. *Then there exists $y \in X^k$ such that $f' = f \cup \{(x, y)\}$ is special.*

*Proof.* Let $y$ be any element of $\{V(A \cup \{x\})(x) \colon V(A) \subseteq f\}$ of minimal sup. We claim that $f' = f \cup \{(x, y)\}$ is special. Clearly $\operatorname{dom}(f')$ is as required. Now fix $A \subseteq \operatorname{dom}(f)$ such that $V(A) \subseteq f$ and $V(A \cup \{x\})(x) = y$. Note that $x \notin A$.

Now let $B \subseteq \operatorname{dom}(f')$ be finite. It suffices to find $D \supseteq B$ such that $V(D) \subseteq f'$.

Since $f$ is special, let $C \supseteq A \cup B \backslash \{x\}$, where $V(C) \subseteq f$. Hence $V(A) \subseteq V(C)$. Now by Theorem 3.10, $V$ is $<,< -*$-decreasing and $<$-end-preserving, where $<$ is the relation $|a| < |b|$. Hence $V(A) \subseteq V(A \cup \{x\})$ and $V(C) \subseteq V(C \cup \{x\})$. By $*$-decreasing, $V(A \cup \{x\})\ (<,<)^*\ V(C \cup \{x\})$. Hence $V(A \cup \{x\})(x) = V(C \cup \{x\})(x)$ or $|V(A \cup \{x\})(x)| > |V(C \cup \{x\})(x)|$. That is, $V(C \cup \{x\})(x) = y$ or $|V(C \cup \{x\})(x)| < |y|$. In the former, we have $V(C \cup \{x\}) \subseteq f'$, in which case we are done by setting $D = C \cup \{x\}$. The latter case is impossible by the minimality of $y$. □

Theorem 4.8. *Every #-decreasing well ordered function assignment has a unique direct completion.*

*Proof.* It suffices to prove that $V^*$ has a unique direct completion. Note that the empty function is special, and the special functions are closed under unions of chains. Since $X$ is well ordered, for all proper subsets $B$ of $X^k$ there exists $x \in X$ which is minimal over $B$. Hence by Lemma 4.7, any maximal special function is a completion of $V$ on $X$. By Zorn's lemma there exists a completion of $V^*$ on $X$.

For uniqueness, let $f, g$ be distinct direct completions of $V^*$ and let $x$ be of minimal sup such that $f(x) \neq g(x)$. Let $x \in A$, where $U(A) \subseteq f$. Let $B \supseteq A$, where $U(B) \subseteq g$. Then by $*$-decreasing, we see that $U(A)(x) = U(B)(x)$ or $|U(A)(x)| > |U(B)(x)|$. Hence $f(x) = g(x)$ or $|f(x)| > |g(x)|$. Now we can interchange $f$, $g$ in the above argument, and obtain $g(x) = f(x)$ or $|g(x)| > |f(x)|$. Hence $f(x) = g(x)$, which is the desired contradiction. □

We now consider the more general concept of completion.

Let $U^*$ be a linearly ordered function assignment and $Y$ be a linearly ordered set. We say that $f$ is a completion of $U^*$ on $Y$ if and only if

(i) $f: Y^k \to Y^k$;

(ii) For all finite $B \subseteq Y^k$ there exists $B \subseteq C \subseteq Y^k$ such that $f \mid C$ is order isomorphic to some $U(A)$.

Note that this concept of completion makes perfectly good sense for function assignments $U$ for $N^k$, since we can view $U$ as the linearly ordered function assignment $U^* = (N^k, \leq, U)$, where $\leq$ is the usual linear ordering on $N$.

We say that $U^*$ has the (unique) ordinal completion property if and only if $U^*$ has a (unique) completion on every well ordered set.

We say that a function assignment $U$ for $N^k$ has the (unique) ordinal completion property if and only if $U$ has (unique) completions on every well ordered set.

Theorem 4.9. *Every infinite #-decreasing linearly ordered function assignment has the ordinal completion property. For all $k > 0$, every #-decreasing function assignment for $N^k$ has the ordinal completion property. Furthermore, if the assignment in question is uniform, then the unique ordinal completion property holds.*

*Proof.* For the first claim, let $U^* = (X^k, \leq, U)$ be an infinite #-decreasing linearly ordered function assignment and $Y$ be a well ordered set. If $Y$ is finite, then it is obvious that there is a completion of $U^*$ on $Y$. If $Y$ is infinite, then by Lemma 4.4, there exists a #-decreasing linearly ordered function assignment $V^* = (Y^k, \leq, V)$ which is finitely conservative over $U^*$. By Theorem 4.8, $V^*$ has a direct completion. Clearly this direct completion is a completion of $U$ on $Y$. The second claim is a special case of the first claim.

For the third claim, let $f$, $g$ be completions of $U^*$ on $Y$, and let $x$ be of minimal sup norm such that $f(x) \neq g(x)$. Let $x \in S$ and $f|S$ be order isomorphic to $U(A)$. Let $S \subseteq T$ and $g|T$ be order isomorphic to $U(B)$. Let $h$ be the unique order isomorphism from $g|T$ onto $U(B)$. Let $C = h[S]$. Then clearly $C$ is order isomorphic to $A$. Hence by order invariance, $U(C)$ is order isomorphic to $U(A)$. Hence $f|S$ is order isomorphic to $U(C)$. Furthermore, this order isomorphism is $h|S$, because $h|S$ is the unique order isomorphism from $S$ onto $C$.

We now claim that for all $y \in C$, if $|y| < |h(x)|$, then $U(C)(y) = U(B)(y)$. Let $y \in C$, $|y| < |h(x)|$. Now $U(C)(y) = h(f(h^{-1}y))$ and $U(B)(y) = h(g(h^{-1}y))$. Since $h$ is order-preserving, $|h^{-1}y| < |x|$. Hence $f(h^{-1}y) = g(h^{-1}y)$. Hence $U(C)(y) = U(B)(y)$.

Therefore, by $*$-decreasing, we have $U(C)(h(x)) = U(B)(h(x))$ or $|U(C)(h(x))| > |U(B)(h(x))|$. Therefore, $f(x) = g(x)$ or $|f(x)| > |g(x)|$.

We can now interchange $f$ and $g$ in the above argument and obtain $g(x) = f(x)$ or $|g(x)| > |f(x)|$. Hence $f(x) = g(x)$, which is the desired contradiction. □

So far we have not used any axioms beyond ZFC. We now fix $k > 0$ and a $k$-ineffable$'$ cardinal $\lambda$. For $A \subseteq \lambda$, let $\mathrm{FST}(A)$ be the set of all finite subsets of $A$.

By elementary set theory, fix a one-one onto function $f\colon \mathrm{FST}(\lambda) \to \lambda$. Let $C = \{x\colon x < \lambda \text{ and } f[\mathrm{FST}(x)] \subseteq x\}$. Clearly $C$ is closed unbounded.

Recall from Part 1 that $\mathrm{ot}(k)$ is the number of order types of elements of $N^k$, and that $\mathrm{ot}(k) \le k^k$.

Lemma 4.10. *Let $g\colon \lambda^k \to \lambda^k$. Then there exists an infinite set $B \subseteq \lambda$ such that $g$ has at most* $\mathrm{ot}(k)$ *regressive values on $B^k$.*

*Proof.* We define an auxiliary regressive function $h\colon S_k(\lambda) \to \lambda$ as follows. Let $A \in S_k(C)$. Let $A^{k\prime}$ be the set of all elements of $A^k$ whose set of coordinates form an initial segment of $A$. Note that $|A^{k\prime}| = \mathrm{ot}(k)$. Define $h(A) = f\{g(x)\colon x \in A^{k\prime} \text{ and } |g(x)| < \min(A)\}$. For other $A \in S_k(\lambda)$, define $h(A) = 0$.

Note that $h\colon S_k(\lambda) \to \lambda$ is regressive. Since $\lambda$ is $k$-ineffable$'$ let $h$ be constantly $c$ on some infinite $S_k(B)$, where $B \subseteq C$. Let $f(T) = c$. By construction, $|T| \le \mathrm{ot}(k)$. We claim that every regressive value of $g$ on $B^k$ lies in $T$. To see this, let $|g(x)| < \min(x)$, $x \in B^k$. Choose $A \in S_k(B)$ so that $x \in A^{k\prime}$. Then $h(A) = c$, and so $f\{g(x)\colon x \in A^{k\prime} \text{ and } |g(x)| < \min(A)\} = f(T)$. Hence $\{g(x)\colon x \in A^{k\prime} \text{ and } |g(x)| < \min(A)\} = T$. Hence $g(x) \in T$. □

Theorem 4.11. *Assume that for all $k$ there exists a $k$-ineffable$'$ cardinal.*

(a) *Let $k, p > 0$ and $U^*$ be an infinite linearly ordered function assignment with the ordinal completion property. Then some $U(A)$ has at most* $\mathrm{ot}(k)$ *regressive values on some $E^k$, $|E| = p$.*

(b) *Let $k, p > 0$ and $U$ be a function assignment for $N^k$ with the ordinal completion property. Then some $U(A)$ has at most* $\mathrm{ot}(k)$ *regressive values on some $E^k$, $|E| = p$.*

*Proof.* Claim (b) is a special case of claim (a). To prove (a), let $f\colon \lambda^k \to \lambda^k$ be a completion of $U^*$ on $\lambda$, where $\lambda$ is a $k$-ineffable$'$ cardinal. By Lemma 4.10, let $B \subseteq \lambda$, $|B| = p$, and $f$ have at most $\mathrm{ot}\,(k)$ regressive values on $B^k$. Choose finite $B^k \subseteq D \subseteq \lambda^k$ and $A \subseteq X^k$ such that $f|D$ is order isomorphic to $U(A)$. Then $U(A)$ has at most $\mathrm{ot}\,(k)$ regressive values on some $E^k$, $|E| = p$, as required. □

THEOREM 4.12. *Assume that for all k there exists a k-ineffable′ cardinal. Then Proposition* A *holds for #-decreasing.*

*Proof.* This is an immediate consequence of Theorems 4.9 and 4.11. □

We now begin the proof of Proposition C. We first need to introduce some terminology that involves a pair of orderings.

Recall that in Part 3, we worked with function assignments, and in this part we have been working with linearly ordered function assignments. We now add some additional structure as follows.

An extended linearly ordered function assignment (elofa) is a quintuple $(X^k, \leq, <_1, <_2, U)$, where $(X^k, \leq, U)$ is a linearly ordered function assignment (lofa) and $<_1, <_2$ are strict orders on $X^k$. Without confusion, we can treat an elofa as a lofa by suppressing $<_1$ and $<_2$, if we so specify.

Order isomorphisms between elofa's are required to preserve the two strict orders, in addition to $\leq$ and $U$. The restriction of an elofa to a set $S$ restricts the strict orders to $S^2$. And $p$-uniformity and uniformity of elofa's is defined using these order isomorphisms and restrictions. Finite conservativity and equivalence are also defined in the obvious way using these isomorphisms and restrictions.

For linear orderings $(X, \leq)$, we let $\leq_c$ be the order on $X^k$ given by $x \leq_c y$ if and only if for all $i$, $x_i \leq y_i$.

We say that $<$ is upward on $X^k$ if and only if $<$ is a strict ordering on $X^k$ such that for all $x, y \in X^k$, if $x \leq_c y$ then $\neg y < x$.

We say that an elofa $(X^k, \leq, <_1, <_2, U)$ is upward if and only if $<_1$, $<_2$ are upward orders on $X^k$.

We also say that an elofa $(X^k, \leq, <_1, <_2, U)$ is #-decreasing ($*$-decreasing, end-preserving) if and only if $U$ is $<_1, <_2$ $-$#-decreasing ($<_1, <_2$ $-*$-decreasing, $<_1$-end-preserving).

We emphasize the contrast between #-decreasing lofa's and #-decreasing elofa's. In the former, #-decreasing is stated in terms of $\leq$. In the later, #-decreasing is stated in terms of $<_1$, $<_2$ and $\leq$ plays no role. The same is true of $*$-decreasing and end-preserving. However, $\leq$ is used in an essential way to define upward.

Thus there is no relationship between an elofa being #-decreasing, $*$-decreasing, or end-preserving, and it being #-decreasing, $*$-decreasing, or end-preserving as a lofa.

The concepts of uniform and order invariant for elofa's are entirely analogous to the lofa context, where $<_1$ and $<_2$ are taken into account in the obvious way. And the concepts of restriction, finitely conservative, finitely equivalent, uniform, and order invariant between elofa's is also entirely analogous to the lofa context, where $<_1$ and $<_2$ are taken into account in the obvious way.

A direct completion of an elofa is simply a direct completion of the elofa as a lofa. A well ordered elofa is an elofa whose $\leq$ is a well ordering.

We now wish to prove the fundamental fact that every upward #-decreasing well ordered elofa has a direct completion. To this end, fix a well ordered upward #-decreasing elofa $V^* = (X^k, \leq, <_1, <_2, V)$.

LEMMA 4.13. *Let $x_1, x_2, \cdots \in X^k$. Then there exists $i < j$ such that $x_i \leq_c x_j$. In fact there is an infinite subsequence such that each term is $\leq_c$ the next.*

*Proof.* This is a standard result from wqo theory (well quasi order theory). For the sake of completeness we give a proof, which depends crucially on the hypothesis that $(X, \leq)$ is a well ordering.

Let $x_1, x_2, \ldots$ be given. It is a well known consequence of $(X, \leq)$ being a well ordering that every infinite sequence from $X$ has an infinite subsequence which is increasing (i.e., each term is $\leq$ the next). Choose an infinite subsequence of the $x$'s whose first coordinates are increasing. Choose an infinite subsequence of this infinite subsequence whose second coordinates are increasing. Repeat the process $k$ times. $\square$

LEMMA 4.14. *Every nonempty $A \subseteq X^k$ has a $<_1$-minimal element*; *i.e., an $x \in A$ such that for no $y \in A$ is $y < x$. The same is true for $<_2$.*

*Proof.* This is a consequence of the hypothesis that $<_1$ is upward on $X^k$. Suppose this is false for nonempty $A$. Let $x_1 \in A$. Inductively define $x_{i+1} < x_i$. This results in an infinite sequence $x_1, x_2, \ldots$ such that each $x_i >_1 x_{i+1}$. Hence for all $i < j$, $x_j >_1 x_i$. But by Lemma 4.13, let $i < j$ be such that $x_i \leq_c x_j$. Then $\neg x_j >_1 x_i$. This is the desired contradiction. $\square$

We say that $f$ is regular if and only if

(i) $f$ is a partial function from $X^k$ into $X^k$ such that for all $y \in \mathrm{dom}(f)$ and $x <_1 y$, we have $x \in \mathrm{dom}(f)$;

(ii) For all finite $A \subseteq \mathrm{dom}(f)$ there exists $B$ containing $A$ such that $V(B) \subseteq f$.

Note that this is the same as the previous concept "$f$ is special" with $|x| < |y|$ replaced by $x <_1 y$.

For $A \subseteq X^k$, we say that $x$ is $<_1$-minimal over $A$ if and only if $x \in X^k \backslash A$, and for all $y <_1 x$, we have $y \in A$.

LEMMA 4.15. *Let $f$ be regular and $x$ be $<_1$-minimal over* $\mathrm{dom}\,(F)$. *Then there exists $y \in X^k$ such that $f' = f \cup \{(x,y)\}$ is regular.*

*Proof.* Let $y$ be any $<_2$-minimal element of $\{(V(A \cup \{x\})(x) \colon V(A) \subseteq f\}$. Note that this set is likely to be infinite, and so we need Lemma 4.14 to obtain $y$. We claim that $f' = f \cup \{(x,y)\}$ is regular. To see this, fix $A \subseteq \mathrm{dom}\,(f)$ such that $V(A) \subseteq f$ and $V(A \cup \{x\})(x) = y$.

Now let $B \subseteq \mathrm{dom}\,(f')$ be finite. It suffices to find $D \supseteq B$ such that $V(D) \subseteq f'$.

Since $f$ is regular, let $C \supseteq A \cup B \backslash \{x\}$, where $V(C) \subseteq f$. Hence $V(A) \subseteq V(C)$. Now by Theorem 3.10, $V$ is $<_1, <_2$ $-*$-decreasing and $<_1$-end-preserving. Hence $V(A) \subseteq V(A \cup \{x\})$ and $V(C) \subseteq V(C \cup \{x\})$. By $*$-decreasing, $V(A \cup \{x\})$ $(<_1, <_2)^*$ $V(C \cup \{x\})$. Hence $V(A \cup \{x\})(x) = V(C \cup \{x\})(x)$ or $V(A \cup \{x\})(x) >_2 V(C \cup \{x\})(x)$. That is, $V(C \cup \{x\})(x) = y$ or $V(C \cup \{x\})(x) <_2 y$. In the former case we have $V(C \cup \{x\}) \subseteq f'$, in which case we are done, by setting $D = C \cup \{x\}$. The latter case is impossible by the $<_2$-minimality of $y$. □

THEOREM 4.16. *Every well ordered upward #-decreasing extended linearly ordered function assignment has a unique direct completion.*

*Proof.* It suffices to prove that $V^*$ has a unique direct completion. Note that the empty function is regular, and the regular functions are closed under unions of chains. Since $X$ is well ordered, by Lemma 4.14 we see that for all $B \subset X$ there exists $x \in X$ which is $<_1$-minimal over $B$ (here $\subset$ indicates proper inclusion). Hence by Lemma 4.15, any maximal regular function is a completion of $V$ on $X$. By Zorn's lemma there exists a completion of $V^*$ on $X$.

For uniqueness, let $f, g$ be distinct direct completions of $V^*$ and let $x$ be $<_1$-minimal such that $f(x) \neq g(x)$. Let $x \in A$, where $U(A) \subseteq f$. Let $B \supseteq A$, where $g \mid B = U(B) \subseteq g$. Then by $*$-decreasing, we see that $U(A)(x) = U(B)(x)$ or $U(A)(x) >_2 U(B)(x)$. Hence $f(x) = g(x)$ or $f(x) >_2 g(x)$. Now we can interchange $f$, $g$ in the above argument, and obtain $g(x) = f(x)$ or $g(x) >_2 f(x)$. Hence $f(x) = g(x)$, which is the desired contradiction. □

We say that an elofa $V^*$ has the (unique) ordinal completion property if and only if $V^*$ has the (unique) ordinal completion property as a lofa.

THEOREM 4.17. *Every infinite #-decreasing upward extended linearly ordered function assignment has the ordinal completion property. For all $k > 0$ and upward orders $<_1$, $<_2$ on $N^k$, every $<_1, <_2$ $-\#$-decreasing function assignment for $N^k$ has the ordinal completion property. Furthermore, if the function assignment is uniform, then the unique ordinal completion property holds.*

*Proof.* For the first claim, let $U^*$ be an infinite upward #-decreasing elofa, and $(Y, \leq)$ be a well ordered set. We can assume without loss of generality that $Y$ is infinite. Analogously to Lemma 4.4, we can find a uniform infinite #-decreasing elofa $U^{*\prime}$ which is finitely conservative over $U^*$. By finite conservativity, $U^{*\prime}$ is also upward. Analogously to Lemma 4.6, we can find a #-decreasing elofa $V^* = (Y^k, \leq, <_1, <_2, V)$ which is finitely conservative over $U^{*\prime}$, and hence finitely conservative over $U^*$, and also upward. We then apply Theorem 4.16.

The second claim is a special case of the first claim.

The third claim is also proved analogously to Theorem 4.9. We simply replace all occurrences of $<$ by $<_1$, all occurrences of $>$ by $>_2$, and remove all absolute value signs. □

THEOREM 4.18. *Assume that for all $k$ there exists a $k$-subtle cardinal. Then Proposition* C *holds. In fact, all forms of Propositions* A–D *hold, even with $k^k$ replaced by* $\mathrm{ot}\,(k)$.

*Proof.* By Theorems 2.2 and 2.7, the assumption implies that for all $k$, there exists a $k$-ineffable′ cardinal. Then apply Theorem 4.17, 4.11, and 3.14. □

## 5. Independence proofs

5.1. *Existence of finite solutions of certain recursive definitions with strong indiscernibles.* We now prove that Proposition A for #-decreasing implies the consistency of ZFC, and in fact $\mathrm{Con}\,(\mathrm{ZFC} + \{\text{there is a } k\text{-subtle cardinal}\}_k)$. By Theorem 3.14, this is enough to establish the independence of all forms of Propositions A–D from the theory $\mathrm{ZFC} + \{\text{there is a } k\text{-subtle cardinal}\}_k$, assuming that this theory is consistent.

Throughout this part, we take Proposition A for #-decreasing as the working assumption. Each of the eight section headings in this part describes the final lemma of that section.

Let $f$ be a nonempty partial function from $N^k \to N^r$ and $E \subseteq N$. We say that $f$ is regressively regular over $E$ if and only if the following hold.

(i) $E^k \subseteq \mathrm{dom}\,(f)$;

(ii) Let $x, y \in E^k$ be of the same order type. If $|f(x)| < \min(x)$, then $|f(y)| < \min(y)$ and $f(x) = f(y)$.

LEMMA 5.1. *Let $k, p > 0$ and $U$ be a #-decreasing function assignment for $N^k$. Then some $U(A)$ is regressively regular over some $E$ of cardinality $p$.*

*Proof.* Let $k$, $p$, $U$ be as given. Choose $q \gg k, p$. Let $U(A)$ have at most $k^k$ regressive values on $E^k \subseteq A$, $|E| = q$. We define $G\colon E^k \to N^k$ by $G(x) = U(A)(x)$ if $|U(A)(x)| < |x|$; $-1$ otherwise.

By the choice of $q$, and the Finite Ramsey Theorem, let $E' \subseteq E$, $|E| = p$, such that the following holds. For all $x, y \in E'^k$, if $x$, $y$ have the same order type, then $G(x) = G(y)$. Then $E'$ is as required. □

For $i \geq 0$ and $j > 0$, we let $i^*j$ be the $j$-tuple all of whose coordinates are $i$.

A function system $U$ for $N^k$ is a mapping $U$ from finite $A \subseteq N^k$ into functions $U(A)\colon A \to \text{fld}(A)$. Thus function systems are a version of function assignments with codimension 1.

We wish to consider very special kinds of function systems called inductive function systems.

Let $\text{FPF}(N^k)$ be the set of all finite partial functions from $N^k$ into $N$. For $f \in \text{FPF}(N^k)$ we let $\text{fld}(f) = \text{fld}(\text{graph}(f))$.

Let $\text{DFNL}(N^k)$ be the set of all $H\colon \text{FPF}(N^k) \times N^k \to N$ such that for all $f \subseteq g$ from $\text{FPF}(N^k)$ and $x \in N^k$,

(i) $H(f, x) \geq H(g, x)$;

(ii) $H(f, x) \in \text{fld}(f) \cup \{x_1, \ldots, x_k\}$.

Here DFNL stands for "decreasing functional."

Each $H \in \text{DFNL}(N^k)$ generates a function system $U$ for $N^k$ by the following inductive process.

Let $A \subseteq N^k$ be finite. We define $\text{RCN}(A, H)$ as the unique $F\colon A \to \text{fld}(A)$ such that for all $x \in A$, $F(x) = H(F|\{y \in A\colon |y| < |x|\}, x)$. Note that for each fixed $H \in \text{DFNL}(N^k)$, $\text{RCN}(A, H)$ defines a function system. We call such a function system an inductive function system for $N^k$. Here RCN stands for "recursion."

For $x, y \in N^k$, we write $x \subseteq y$ if and only if every coordinate of $x$ is a coordinate of $y$.

We say that $A \subseteq N^k$ is closed if and only if for all $x \subseteq y$ with $y \in A$, we have $x \in A$.

Lemma 5.2. *Let $k, p > 0$ and $U$ be an inductive function system for $N^k$. Then there exists finite closed $A$ such that $U(A)$ is regressively regular over some $E$ of cardinality $p$.*

*Proof.* For finite $A \subseteq N^k$, write $A' = \{x \in A\colon \text{for all } y \subseteq x,\ y \in A\}$. We modify the RCN construction as follows.

Let $H \in \text{DFNL}(N^k)$. Define $M\text{RCN}(A, H)$ as the unique $F\colon A \to \text{fld}(A)$ such that for all $x \in A$, $F(x) = H(F|\{y \in A'\colon |y| < |x|\}, x)$. We wish to show that $M\text{RCN}(A, H)$ defines a function system for $N^k$.

For this purpose, we define $H'\colon \text{FPF}(N^k) \times N^k \to N$ by $H'(f, x) = H(f|\{y \in \text{dom}\,(f)'\colon |y| < |x|\}, x)$. Note that if $f \subseteq g$, then $\text{dom}\,(f)' \subseteq \text{dom}\,(g)'$.

Hence $H' \in \mathrm{DFNL}(N^k)$. And note that $M\mathrm{RCN}(A,H) = \mathrm{RCN}(A,H')$. So $M\mathrm{RCN}(A,H)$ does define a function system for $N^k$.

Note that for all $x \in A$, $M\mathrm{RCN}(A,H)(x)$ is either a coordinate of some $y \in A'$ with $|y| < |x|$, or a coordinate of $x$. This is proved by induction on $|x|$.

Fix $k$, $p$, $U$ as given. Fix $H \in \mathrm{DFNL}(N^k)$ so that $U$ is the inductive function system given by $\mathrm{RCN}(A,H)$, as a function of $A$. It is convenient to define $A^\wedge = \{x \in A$: for all $y \subseteq x$, if $|y| < |x|$, then $y \in A\}$.

Let $V$ be the function assignment for $N^k$, where $V(A)\colon A \to A$ is defined by cases as follows. Let $x \in A$.

(i) $x \in A^\wedge$ and $M\mathrm{RCN}(A,H)(x) < |x|$. Define

$$V(A)(x) = M\mathrm{RCN}(A,H)(x)^*k\,.$$

(ii) Otherwise. Define $V(A)(x) = x$.

To verify that $V$ is indeed a function assignment, note that $M\mathrm{RCN}(A,H)\ (x)^*k \subseteq x$ or $\subseteq$ some element of $A'$. Hence, in (i), $M\mathrm{RCN}(A,H)(x)^*k \in A$.

We now claim that $V$ is #-decreasing. To see this, fix finite $A \subseteq N^k$, $x \in N^k$, and assume $\neg V(A) \subseteq V(A \cup \{x\})$. Fix $y \in A$ of minimal sup such that $V(A)(y) \neq V(A \cup \{x\})(y)$.

Note that the value of any $M\mathrm{RCN}(B,H)(y)$ depends only on $y$ and the elements of $B$ of norm $< |y|$. From this it follows that the value of any $V(B)(y)$ also depends only on $y$ and the elements of $B$ of norm $< |y|$.

Now $A$ and $A \cup \{x\}$ have the same elements $|z| < |x|$. Hence $V(A)$ and $(V(A \cup \{x\})$ agree at all $|z| \leq |x|$. We have thus shown that $|y| > |x|$.

We now claim that $M\mathrm{RCN}(A,H)$ and $M\mathrm{RCN}(A \cup \{x\},H)$ agree at all $z \in A'$ with $|z| < |y|$. To see this, fix $z \in A'$, $|z| < |y|$. Then $V(A)(z) = V(A \cup \{x\})(z)$. But $V(A)(z) = M\mathrm{RCN}(A,H)(z)^*k$ if $MRCN(A,H)(z) < |z|$; $z$ otherwise. And $V(A \cup \{x\})(z) = M\mathrm{RCN}(A \cup \{x\},H)(z)^*k$ if $M\mathrm{RCN}(A \cup \{x\},H)(z) < |z|$; $z$ otherwise. Since $M\mathrm{RCN}(A,H)(z) \leq |z|$ and $M\mathrm{RCN}(A \cup \{x\},H)(z) \leq |z|$, the claim is established.

We now see that by the definition of $M\mathrm{RCN}(A,H)(y)$ and the decreasing condition on $H$, that $M\mathrm{RCN}(A,H)(y) \geq M\mathrm{RCN}(A \cup \{x\},y)$. We now claim that $|V(A)(y)| > |V(A \cup \{x\})(y)|$, which completes the proof of the claim that $V$ is #-decreasing. We argue by cases.

Recall that $V(A)(y) \neq V(A \cup \{x\})(y)$. Note that $V(A)(y) = y$ or $|V(A)(y)| < |y|$, and $V(A \cup \{x\})(y) = y$ or $|V(A \cup \{x\})(y)| < |y|$.

*Case* 1. $y \in A^\wedge$ and $M\mathrm{RCN}(A,H)(y) < |y|$. Then $V(A)(y) = M\mathrm{RCN}(A,H)(y)^*k$ and $V(A \cup \{x\})(y) = M\mathrm{RCN}(A \cup \{x\})(y)^*k$. Therefore, $M\mathrm{RCN}(A,H)(y) > M\mathrm{RCN}(A \cup \{x\},y)$. Hence $|V(A)(y)| > |V(A \cup \{x\})(y)|$.

*Case* 2. Otherwise. Then $V(A)(y) = y$. Hence $|V(A \cup \{x\})(y)| < |y| = |V(A)(y)|$.

We have thus established that $V$ is #-decreasing. Now fix $A$ and $|E| = p$ such that $V(A)$ is regressively regular on $E$.

We want to prove that $\mathrm{RCN}(A', H)$ is regressively regular on $E$. This will complete the proof of the lemma, since obviously $A'$ is closed.

By examining clauses (i) and (ii), we see that for all $x \in A'$,

$$M\mathrm{RCN}(A, H)(x) = |V(A)(x)|\,.$$

Also observe that $\mathrm{RCN}(A', H) \subseteq M\mathrm{RCN}(A, H)$ by induction on the sup norm of $x \in A'$. And since $E^k \subseteq A$, we see that $E^k \subseteq A'$.

Now let $x, y \in E^k$ be of the same order type, and $\mathrm{RCN}(A', H)(x) < \min(x)$. Since $x \in A'$, $M\mathrm{RCN}(A, H)(x) < \min(x)$. Hence $V(A)(x) = M\mathrm{RCN}(A, H)(x)^*k$ and $|V(A)(x)| < \min(x)$. So $V(A)(x) = V(A)(y)$ and $|V(A)(y)| < \min(y)$. Since $y \in A'$, $V(A)(y) = M\mathrm{RCN}(A, H)(y)^*k$. Thus

$$M\mathrm{RCN}(A, H)(x) = M\mathrm{RCN}(A, H)(y) < \min(y)\,.$$

Hence $\mathrm{RCN}(A', H)(x) = \mathrm{RCN}(A', H)(y) < \min(y)$ as required. □

We are going to apply Lemma 5.2 to certain inductive function systems that are defined in logical terms.

Let $QF$ be the set of all propositional combinations of atomic formulas of the form $x < y$, where $x$ and $y$ are variables representing elements of $N$. We require that elements of $QF$ are in disjunctive normal form. We use actual variables $x_1, x_2, \ldots$.

Let $f \in \mathrm{FPF}(N^r)$. If $y \in N^{rt}$, then we write $f(y)$ for

$$(f(y_1, \ldots, y_r), \ldots, f(y_{rt-r+1}, \ldots, y_{rt}))\,.$$

Here it is important to adhere to the convention that $f(y)$ is defined if and only if each of the $t$ components is defined. That is, $f(y)$ is defined if and only if $y \in \mathrm{dom}\,(f)^t$.

Let $\mathrm{BEF}(q, t, r)$ be the set of all bounded existential formulas of the following form:

$$B(x) = (\exists\, y \in \mathrm{dom}\,(F)^t)(|y| < |x|\ \&\ D(x, y, F(y)))\,,$$

where $x$ abbreviates the list of variables $x_1, \ldots, x_q$, $y$ abbreviates the list of variables $x_{q+1}, \ldots, x_{q+rt}$, and $D$ is in $QF$. Here $F$ is viewed as a function symbol representing a finite partial function from $N^r \to N$. Note that there are no nested occurrences of $F$ in $B(x)$.

We let $\mathrm{BEF}(q, r)$ be the union over $t$ of the $\mathrm{BEF}(q, t, r)$.

If we specify an actual $f \in \mathrm{FPF}(N^r)$ and $x \in N^q$, then it is clear what we mean by asserting that $B(x)$ is true in $F$.

Let $B \in \mathrm{BEF}(r+1, r)$ and $A$ be a finite subset of $N^r$. We define $Df(B; A)$ as the unique $f\colon A \to \mathrm{fld}(A)$ such that for all $x \in A$, $f(x) = \min\{j \in \mathrm{fld}(A)\colon j = |x|$ or $B(x, j)$ is true in $f\}$.

Lemma 5.3. *Let $r, p > 0$ and $B \in \mathrm{BEF}(r+1, r)$. Then there exists finite closed $A \subseteq N^r$ such that $Df(B; A)$ is regressively regular over some $E$ of cardinality $p$.*

*Proof.* Let $r$, $p$, $B$ be as given. Define $H(F, x) = \min\{j \in \mathrm{fld}(F) \cup \{x_1, \ldots, x_k\} \colon j = |x| \text{ or } B(x, j) \text{ is true in } F\}$. Then $H \in \mathrm{DFNL}(N^k)$. By Lemma 5.2, let $A$ be a finite closed subset of $N^k$ and $E \in S_p(N)$ such that $\mathrm{RCN}(A, H)$ is regressively regular on $E$.

Now for all $x \in A$, $\mathrm{RCN}(A, H)(x) = H(\mathrm{RCN}(A, H)|\{y \in A \colon |y| < |x|\}, x)$. Since $A$ is closed, $\mathrm{fld}(\mathrm{RCN}(A, H)|\{y \in A \colon |y| < |x|\}) = \{j \in \mathrm{fld}(A) \colon j < |x|\}$. Hence $\mathrm{RCN}(A, H)(x) = \min\{j \in \{j \in \mathrm{fld}(A) \colon j < |x|\} \cup \{x_1, \ldots, x_k\} \colon j = |x| \text{ or } B(x, y) \text{ is true in } F\}$. That is, $\mathrm{RCN}(A, H)(x) = \min\{j \in \mathrm{fld}(A) \colon (j = |x| \text{ or } j \in \{x_1, \ldots, x_k\}) \text{ and } (j = |x| \text{ or } B(x, j) \text{ is true in } F)\}$. Therefore, $\mathrm{RCN}(A, H)(x) = \min\{j \in \mathrm{fld}(A) \colon j = |x| \text{ or } B(x, j) \text{ is true in } F)\}$. To see this last step, note that we are claiming that the min over two related sets are equal. Note that both sets contain $|x|$. Also note that both sets agree below $|x|$.

Thus we have shown that $\mathrm{RCN}(A, H)$ obeys the defining condition for $Df(B; A)$. Hence $\mathrm{RCN}(A, H) = Df(B; A)$, and this completes the proof. □

Let $k > r > 0$, $f \in \mathrm{FPF}(N^k)$, and $A = \mathrm{dom}(f)$. We define $A/r = \{x \in N^r \colon (x, |x| * k - r) \in A\}$, and $f/r \colon A/r \to N$ by $f/r(x) = f(x, |x| * k - r)$. We think of $A/r$ and $f/r$ as the $r$-dimensional parts of $A$ and $f$.

Lemma 5.4. *Let $A \subseteq N^k$ be closed and $k > r$. Then $A/r$ is closed, and $\mathrm{fld}(A/r) = \mathrm{fld}(A)$. Furthermore, let $f \colon A \to \mathrm{fld}(A)$ be regressively regular over $E$. Then $f/r$ is regressively regular over $E$.*

*Proof.* Let $A$, $k$, $r$, $f$ be as given. To see that $A/r$ is closed, let $x \subseteq y \in A/r$. Then $(x, |x|^*k - r) \subseteq (y, |y|^*k - r) \in A$. Hence $(x, |x|^*k - r) \in A$, and so $x \in A/r$.

Now observe that $j \in \mathrm{fld}(A)$ if and only if $j^*k \in A$. Also note that $j \in \mathrm{fld}(A/r)$ if and only if $j^*r \in A/r$ if and only if $j^*k \in A$.

Finally, to show that $f/r$ is regressively regular over $E$, let $x, y \in E^r$ have the same order type and $f/r(x) < \min(x)$. Then $f(x, |x|^*k - r) = f/r(x) < \min(x)$. Now $(x, |x|^*k - r)$ and $(y, |y|^*k - r)$ are elements of $E^r$ with the same order type. Hence $f(x, |x|^*k - r) = f(y, |y|^*k - r) < \min(y)$. Therefore $f/r(x) = f/r(y) < \min(y)$ as required. □

We now introduce some notation for the next lemma. Let $q, t > 0$ and $y \in N^{qt}$. The coordinates of $y$ are obviously separated into $q$ consecutive blocks, each of length $t$. Let $1 \leq i \leq q$ and $1 \leq j \leq t$. We let $y\#(i, j)$ be the $qt$-tuple obtained by locating the coordinate of $y$ in the $j$-th position in the $i$-th block, and replacing all succeeding coordinates of $y$ with this coordinate.

LEMMA 5.5. *Let $p, q, r, t > 0$, $s \geq r$, $B \in \mathrm{BEF}(r+1, r)$, and $C_1, \ldots, C_q \in \mathrm{BEF}(s+t, r)$. Then there exist finite closed $A \subseteq N^{s+qt}$ and $f\colon A \to \mathrm{fld}(A)$ such that the following hold*:

(i) *For all $(x, |x|^* s+qt-r) \in A$, $f(x, |x|^* s+qt-r) = \min\{j \in \mathrm{fld}(A)\colon j = |x|$ or $B(x, j)$ holds in $f/r\}$*;

(ii) *Suppose that $x \in N^s$, $y \in N^{qt}$, $|x| < y_1 < \cdots < y_{qt}$, $1 \leq i \leq q$. Then $(f(x, y\#(i,1)), \ldots, f(x, y\#(i,t)))$ (if defined) is the lexicographically least $z \in \mathrm{fld}(A)^t$ such that $|z| < \min(x)$ and $C_i(x, z)$ holds in $f/r$ if it exists; $|x|^* t$ otherwise*;

(iii) *$f$ is regressively regular over some $E$ of cardinality $p$.*

*Proof.* Let $r$, $B$, $C_1, \ldots, C_r$ be given. We use Lemma 5.3 which asserts that any appropriate definition by recursion on the sup norm can be made on some finite closed $A$ such that the function is regressively regular over some set of cardinality $p$.

Accordingly, we can fix finite $A \subseteq N^{s+qt}$ and define $f\colon A \to \mathrm{fld}(A)$ by recursion on the sup norm as follows. Let $(x, v) \in A$, where $x \in N^s$ and $v \in N^{qt}$.

*Case* 1. $(x, v) = (x', |x'|^* s+qt-r)$. Define $f(x, v) = \min\{j \in \mathrm{fld}(A)\colon j = |x'|$ or $B(x, j)$ holds in $f/r\}$.

*Case* 2. $v = y\#(i,1)$ for some $1 \leq i \leq q$ and strictly increasing $y \in N^{qt}$, where $|x| < y_1$. If there exists $z \in \mathrm{fld}(A)^t$ such that $|z| < \min(x)$ and $C_i(x, z)$ holds in $f/r$, then define $f(x, v)$ to be the least $z_1$ among such $z$'s. Otherwise, define $f(x, v) = |x|$.

*Case* 3. $v = y\#(i,j)$ for some $1 \leq i \leq q$ and $1 < j \leq t$, and strictly increasing $y \in N^{q^t}$, where $|x| < y_1$. If there exists $z \in \mathrm{fld}(A)^t$ such that $|z| < \min(x)$ and $C_i(x, z)$ holds in $f/r$ and $z_1 = f(x, y\#(i,1))\,\&\,\ldots\,\&z_{j-1} = f(x, y\#(i, j-1))$, then define $f(x, v)$ to be the least $z_j$ among such $z$'s. Otherwise, define $f(x, v) = |x|$.

*Case* 4. Otherwise, define $f(x, v) = |x|$.

The quantifications over $\mathrm{fld}(A)^t$ are represented by $t$ quantifications over $\mathrm{fld}(A)$, which are in turn represented by a quantification over diagonal elements of $A$ (i.e., elements of $A$ all of whose coordinates are identical). If $A$ is closed, then these representations are valid.

By Lemma 5.3, fix finite closed $A \subseteq N^{s+qt}$ and $f\colon A \to \mathrm{fld}(A)$ defined as above, where $f$ is regressively regular over $E \in S_p(N)$. Then $f$, $E$ are as required. ∎

Lemma 5.6. *Let* $p, q, r, t > 0$, $s \geq r$, $B \in \mathrm{BEF}(r+1, r)$, *and* $C_1, \ldots, C_q \in \mathrm{BEF}(s+t, r)$. *Then there exists finite closed* $A \subseteq N^{s+qt}$, $f \colon A \to \mathrm{fld}(A)$, *and* $E \subseteq N$, $|E| = p$, *such that the following hold*:

(i) *For all* $(x, |x|^* s + qt - r) \in A$, $f(x, |x|^* s + qt - r) = \min\{j \in \mathrm{fld}(A) \colon j = |x|$ *or* $B(x, j)$ *holds in* $f/r\}$;

(ii) *Suppose that* $x \in E^s$, $u \in \mathrm{fld}(A)^t$, $|u| < \min(x)$, $1 \leq i \leq q$, *and* $C_i(x, u)$ *holds for* $f/r$. *Then there exists* $w \in \mathrm{fld}(A)^t$ *such that for all* $y \in E^s$ *of the same order type as* $x$, $C_i(y, w)$ *holds in* $f/r$;

(iii) $f$ *is regressively regular over* $E$.

*Proof.* Let $p$, $q$, $r$, $s$, $t$, $B$, $C_1, \ldots, C_q$ be given. Let $A \subseteq N^{s+qt}$ be finite closed, and $f \colon A \to \mathrm{fld}(A)$ obeys (i)–(iii) of Lemma 5.5, with $|E| = p + qt$. Let $E' = \{E_1, \ldots, E_p\}$. We now show that $f$, $E'$ is as required.

Suppose that $x \in E'^s$, $u \in \mathrm{fld}(A)^t$, $|u| < \min(x)$, and $C_i(x, u)$ holds in $f/r$. Let $|x| < y_1 < \cdots < y_{qt}$, where $y_1, \ldots, y_{qt} \in E$. Then $|w| < \min(x)$ and $C_i(x, w)$ holds for $f/r$, where

$$w = (f(x, y\#(i,1)), f(x, y\#(i,2)), \ldots, f(x, y\#(i,t)))\,.$$

Now by the regressive regularity of $f$ over $E$, this $w$ depends only on the order type of $x$, and not on $y_1, \ldots, y_{qt}$. This establishes (ii) as required. □

Lemma 5.7. *Let* $p, q, r, t > 0$, $s \geq r$, $B \in \mathrm{BEF}(r+1, r)$, *and* $C_1, \ldots, C_q \in \mathrm{BEF}(s+t, r)$. *Then there exist finite closed* $A \subseteq N^r$ *and* $E \in S_p(N)$ *such that the following hold*:

(i) *Suppose that* $x \in E^s$, $u \in \mathrm{fld}(A)^t$, $|u| < \min(x)$, *and* $C_i(x, u)$ *holds in* $Df(B; A)$. *Then there exists* $w \in \mathrm{fld}(A)^t$ *such that for all* $y \in E^s$ *of the same order type as* $x$, $C_i(y, w)$ *holds in* $Df(B; A)$;

(ii) $Df(B; A)$ *is regressively regular over* $E$.

*Proof.* Let $p$, $q$, $r$, $s$, $t$, $B$, $C_1, \ldots, C_q$ be as given. Let $A$, $f$, $E$ be as given by Lemma 5.6. By clause (i) of Lemma 5.6, we see that for all $x \in A/r$, $f/r(x) = \min\{j \in \mathrm{fld}(A) \colon j = |x|$ or $B(x, j)$ holds in $f/r\}$. By Lemma 5.4, $\mathrm{fld}(A) = \mathrm{fld}(A/r)$. Hence $f/r = Df(B; A)$. This also establishes (i) here. And (ii) follows by another application of Lemma 5.4. □

We can reformulate Lemma 5.7 in terms of a strengthening of regressive regularity.

We say that $f$ is $(t, r)$-regular over $E$ if and only if

(i) $t, r > 0$;

(ii) $f \colon A \to \mathrm{fld}(A)$, where $A \subseteq N^r$ is finite and closed;

(iii) $E^t \subseteq A$;

(iv) Let $C \in \mathrm{BEF}(2t, t, r)$, $x \in E^t$, $u \in \mathrm{fld}(A)^t$, $|u| < \min(x)$, and $C(x, u)$ hold in $f$. Then there exists $w \in \mathrm{fld}(A)^t$ such that for all $y \in E^t$ of the same order type as $x$, we have that $C(y, w)$ holds in $f$.

LEMMA 5.8. *Let $r, p > 0$, $t \geq r$, $B \in \mathrm{BEF}(r+1, r)$. Then some $Df(B; A)$ is $(t, r)$-regular over some $E$ of cardinality $p$.*

*Proof.* Let $r$, $t$, $p$, $B$ be given. We will apply Lemma 5.7 for suitably chosen parameters. Let $C_1, \ldots, C_q$ enumerate all elements of $\mathrm{BEF}(2t, t, r)$. The result follows by Lemma 5.7 since $s = t$. □

LEMMA 5.9. *Let $t \geq s > 0$ and $E \subseteq N$. Every $(t, r)$-regular function over $E$ is $(s, r)$-regular over $E$.*

*Proof.* Let $t$, $s$ be as given, and $f: A \to \mathrm{fld}(A)$ be $(t, r)$-regular. Let $C \in \mathrm{BEF}(2s, s, r)$, $x \in E^s$, $u \in \mathrm{fld}(A)^s$, $|u| < |x|$, and $C(x, u)$ hold in $f$. Let $C' \in \mathrm{BEF}(2t, t, r)$ be defined by $C'(z, w)$ if and only if $C(z_1, \ldots, z_s, w_1, \ldots, w_s)$. Then $C'(x, |x|^*t - s, u, |u|^*t - s)$ holds in $f$. Let $w \in \mathrm{fld}(A)^t$ be such that for all $y \in E^t$ of the same order type as $(x, |x|^*t - s)$, we have $C'(y, w)$. Then for all $y' \in E^s$ of the same order type as $x$, we have that $C'(y', |y|^*t - s, w)$ holds in $f$. Hence for all $y' \in E^s$ of the same order type as $x$, we have $C(y', w_1, \ldots, w_s)$ holds in $f$. □

We now prove an important consequence of $(2r, r)$-regularity. Let $x, y \in N^k$ and $j \in N$. We say that $x, y$ are $j$-related if and only if

(i) $x$, $y$ have the same order type;
(ii) for all $1 \leq i \leq k$, if $x_i \leq j$, then $y_i = x_i$;
(iii) for all $1 \leq i \leq k$, if $x_i > j$, then $y_i > |x|$.

LEMMA 5.10. *Let $r > 1$ and $p \geq 4r$ and $f$ be $(2r, r)$-regular over $E$, where $|E| = p$. Let $x, y \in E^r$ be $f(x)$-related. Then $f(x) = f(y)$.*

*Proof.* Let $r, p, f, x, y$ be as given and let $f : A \to \mathrm{fld}\,(A)$. Write $x = (E_{i_1}, \ldots, E_{i_r})$ and $y = (E_{j_1}, \ldots, E_{j_r})$, where $x, y$ are $f(x)$-related. Assume that $f(x) \neq f(y)$. We can push the $i$'s and $j$'s down to $a$'s and $b$'s such that $E(_{a_1}, \ldots, E_{a_r}, E_{b_1}, \ldots, E_{b_r})$ has the same order type as $(E_{a_1}, \ldots, E_{a_r}, E_{b_1}, \ldots, E_{b_r})$, and $\{a_1, \ldots, a_r, b_1, \ldots, b_r)$ is an initial segment of the positive natural numbers. Thus the largest of the $a$'s and $b$'s is $2r$.

Consider the true statements that $f(E_{i_1}, \ldots, E_{i_r}) \neq f(E_{j_1}, \ldots, E_{j_r})$ and $E_{i_1}, \ldots, E_{i_r}, E_{j_1}, \ldots, E_{j_r}$ are $f(E_{i_1}, \ldots, E_{i_r})$-related. From $(2r, r)$-regularity, we can derive an appropriate simple indiscernibility, thereby obtaining that $f(E_{a_1}, \ldots, E_{a_r}) \neq F(E_{b_1}, \ldots, E_{b_r})$, and $(E_{a_1}, \ldots, E_{a_r})$, $(E_{b_1}, \ldots, E_{b_r})$ are $f(E_{a_1}, \ldots, E_{a_r})$-related. See Lemma 5.11.

Next, note that the part of $E_{a_1}, \dots, E_{a_r}$ that lies strictly above $f(E_{a_1}, \dots, E_{a_r})$ corresponds to the part of $E_{b_1}, \dots, E_{b_r}$ that lies strictly above $f(E_{b_1}, \dots, E_{b_r})$. That is, the two parts lie in the same positions among the $a$'s and $b$'s. And the remaining parts are identical, i.e., are the same elements of $E$ in the same positions. This is clear again by simple indiscernibility.

Thus in order to get the required contradiction, it suffices to show that we can

(i) move the part of $E_{a_1}, \dots, E_{a_r}$ that lies strictly above $f(E_{a_1}, \dots, E_{a_r})$ — henceforth called the upper part — all the way to the top of $E$ (preserving order type) without changing the value of $f$; and

(ii) move the part of $E_{b_1}, \dots, E_{b_r}$ that lies strictly above $f(E_{b_1}, \dots, E_{b_r})$ all the way to the top of $E$ (preserving order type) without changing the value of $f$.

As the cases are symmetric, we may work only with $E_{a_1}, \dots, E_{a_r}$.

Since $p \geq 4r$, we can move the upper part of $E_{a_1}, \dots E_{a_r}$ so that the new $E$'s lie among $E_{2r+1}, \dots, E_{3r}$. We write this as $E_{c_1}, \dots, E_{c_r}$. Also we write $E_{d_1}, \dots, E_{d_r}$ for the second move where the new part is at the top of $E$, i.e., forms a tail in $E$. Hence the new part for the second move lies among $E_{3r+1}, \dots, E_p$.

We need to verify that $f(E_{a_1}, \dots, E_{a_r}) = f(E_{d_1}, \dots, E_{d_r})$. Assume this is false. Then by simple indiscernibility, we have

$$f(E_{a_1}, \dots, E_{a_r}) \neq f(E_{c_1}, \dots, E_{c_r}) \neq f(E_{d_1}, \dots, E_{d_r}) .$$

We can rewrite these three distinct values in the form

$$\begin{aligned} f'(E_{x_1}, \dots, E_{x_q}, E_{y_1}, \dots, E_{y_s}) &\neq f'(E_{x_1}, \dots, E_{x_q}, E_{z_1}, \dots, E_{z_s}) \\ &\neq f'(E_{x_1}, \dots, E_{x_q}, E_{w_1}, \dots, E_{w_s}) , \end{aligned}$$

where $E_{x_1}, \dots E_{x_q}$ is the common lower part, and the $y$'s, $z$'s, and $w$'s are the upper parts. Here

$$x_1 < \cdots < x_q < y_1 < \cdots < y_s < z_1 < \cdots z_s < w_1 < \cdots < w_s ,$$

and $f'$ results from $f$ by permuting arguments. Note that these values are $< E_{y_1}$.

Observe that there exist $u_1, \dots, u_q$, $v < E_{y_1}$ such that

$$f'(u_1, \dots, u_q, E_{y_1}, \dots, E_{y_s}) = v \neq f'(u_1, \dots, u_q, E_{z_1}, \dots, E_{z_s}) .$$

Hence by $(2r, r)$-regularity, there exist $u_1, \dots, u_q, v < E_1$ such that

$$f'(u_1, \dots, u_q, E_{y_1}, \dots, E_{y_s}) = v \neq f'(u_1, \dots, u_q, E_{z_1}, \dots, E_{z_s}) ,$$

and

$$f'(u_1, \ldots, u_q, E_{z_1}, \ldots, E_{z_s}) = v \neq f'(u_1, \ldots, u_q, E_{w_1}, \ldots, E_{w_s})\,.$$

But this is a contradiction. □

Here is an important weak consequence of $(t, r)$-regularity.

LEMMA 5.11. *Let* $t, r > 0$, $C \in \mathrm{BEF}(t, t, r)$, *and* $f$ *be* $(t, r)$*-regular over* $E$. *Let* $x, y \in E^t$ *be of the same order type, where* $\min(x), \min(y) > \min(\mathrm{fld}(f))$. *Then* $C(x)$ *holds if and only if* $C(y)$ *holds in* $f$.

*Proof.* Let $t$, $r$, $C$, $f$, $x$, $y$ be as given. Define $C' \in \mathrm{BEF}(2t, t, r)$ by $C'(x, y)$ if and only if $C(x)$ holds. Now let $C(x)$ hold in $f$. Then $C'(x, z)$ holds in $f$, where $z = \min(\mathrm{fld}(A))^*t$. Hence $C'(y, z)$ holds in $f$, and therefore also $C(y)$. □

We now introduce another concept of indiscernibility.

We say that $E$ is a $(t, r) -$ SOI (set of indiscernibles) for $f$ if and only if

(i) $t, r > 0$;

(ii) $f\colon A \to \mathrm{fld}(A)$, where $A \subseteq N^r$ is finite and closed;

(iii) $E^r \subseteq A$;

(iv) Let $C \in \mathrm{BEF}(3t, t, r)$, $x, y, z \in E^t$, $w \in f[E^r]^t$, $x$, $y$ of the same order type, and $|z, w| < \min(x, y)$. Then $B(x, z, w)$ if and only if $B(y, z, w)$.

Note that in the inequality in (iv), we are using a comma for concatenation.

We leave it to the reader to check that if $E$ is a $(t, r)$-SOI for $f$ and $t > s$, then $E$ is an $(s, r)$-SOI for $f$.

LEMMA 5.12. *Let* $r, t > 0$ *and* $p > 3r$ *and* $f$ *be* $(2t(r + 2), r)$*-regular over* $E$, *where* $|E| = p + r$. *Then* $\{E_2, \ldots, E_p\}$ *is a* $(t, r)$-SOI *for* $f$.

*Proof.* Let $p$, $t$, $r$, $f$, $E$ be as given, and let $f\colon A \to \mathrm{fld}(A)$. Without loss of generality, we can assume that $r, t, p > 1$. Let $E' = \{E_2, \ldots, E_p\}$. Let $C \in \mathrm{BEF}(3t, t, r)$, $x, y, z \in E'^t_{u_t}$, and $x$, $y$ have the same order type, and $|z| < \min(x, y)$. Also assume $u_1, \ldots, u_t \in E'^r$, where $|f(u_1), \ldots, f(u_t)| < \min(x, y)$. We need to prove that $C(x, z, f(u_1), \ldots, f(u_t))$ if and only if

$$C(y, z, f(u_1), \ldots, f(u_t))$$

holds in $f$.

Start with the assertion $C(x, z, f(u_1), \ldots, f(u_t))$ holds, which is treated as a statement about the tuples $x$, $z$ from $E'$, with the $t$ parameters shown from $\mathrm{fld}(A)$.

Now by Lemma 5.9, $f$ is $(2r, r)$-regular over $E$, and so by Lemma 5.10, we can write each $f(u_i)$ as $f(u_i^*) = f(u_i)$, where $u_i^* \in E^r$ is of the same order type as $u_i$, where $u_i$ and $u_i^*$ are $f(u_i)$-related, and every coordinate of $u_i^* > f(u_i)$ lies in $\{E_{p+1}, \ldots, E_{p+r}\}$. We thus have the equivalent assertion $C(x, z, f(u_i^*), \ldots f(u_i^*))$.

We now view $C(x, z, f(u_1^*), \ldots, f(u_t^*))$ as $C'(x, z, u_1^*, \ldots, u_t^*)$, where $C' \in \mathrm{BEF}(t(r+2), 2t, r)$. Note that $(x, z, u_1^*, \ldots, u_t^*)$ and $(y, z, u_1^*, \ldots, u_t^*)$ are of the same order type. Hence by $(t(r+2), r)$-regularity and Lemma 5.11, $C'(x, z, u_1^*, \ldots, u_t^*)$ holds if and only if $C'(y, z, u_1^*, \ldots, u_t^*)$ holds in $f$. Therefore, $C(x, z, f(u_1^*), \ldots, f(u_t^*))$ if and only if $C(y, z, f(u_1^*), \ldots, f(u_t^*))$, as required. □

LEMMA 5.13. *Let* $r, t, p > 0$ *and* $B \in \mathrm{BEF}(r+1, r)$. *Then there exists finite closed* $A \subseteq N^r$ *such that* $Df(B; A)$ *has a* $(t, r)$-SOI *of cardinality* $p$. *Furthermore we can require that* $\max(\mathrm{fld}(A)) = \max$ *of the* SOI.

*Proof.* Let $r$, $t$, $p$, $B$ be as given. We can assume without loss of generality that $p > 3r$. By Lemma 5.8, let $Df(B; A)$ be $(2t(r+2), r)$-regular over $E$ of cardinality $p + r$. By Lemma 5.12, $\{E_2, \ldots, E_{p+1}\}$ is a $(t, r)$-SOI for $Df(B; A)$ as required. To obtain the further requirement, simply cut down $A$ to the elements of norm at most the max of the SOI. □

5.2. *Existence of finite towers of finite sets with strong atomic indiscernibles, weak staggered comprehension, and a weak least element principle.* We now start a rather technical development which culminates in Lemma 5.33. We need to build appropriate finite towers of finite sets so that we can apply compactness in order to obtain workable linearly ordered structures.

Let $k, r, t > 0$. We say that $h: N^t \to (N \cup N^k)^r$ is basic if and only if there is a quantifier-free formula $D$ (with only $<$ and $=$ on $N$) which defines $h$ in the sense that $h(x) = y$ if and only if $D(x, y)$. No parameters are allowed. This involves mixed sorts and tuples of tuples, which are handled in the obvious way. We say that $x \in N^t$ is constant if all of its terms are equal.

Note that if $h: N^t \to (N \cup N^k)^r$ is basic and $E \subseteq N$, then $h: E \to (E \cup E^k)^r$.

Let $E \subseteq N$ and $n > 0$. We define $E[n] = \{E_1, \ldots, E_n\}$. For $x \in (N \cup N^k)^r$, $|x| = \max(|x_1|, \ldots, |x|_r)$.

LEMMA 5.14. *Let* $r, p > 0$. *Then there exists* $t > 0$, $0 < n_1 < \cdots < n_p \leq t$, *and basic* $h: N^t \to (N \cup N^{t+r})^r$ *such that the following hold*:

(i) *For all* $1 \leq i < p$ *and* $E \subseteq N$, $h$ *maps* $E[n_{i+1} - r]^t$ *onto* $(E[n_i + r] \cup E[n_i + r]^{t+r})^r$;
(ii) *For all* $j \in N$, $h(j^*t) = j^*r$;
(iii) *If* $x \in N^t$ *is not constant, then* $|h(x)| < |x|$.

*Proof.* Let $r, p > 0$. Set $n_1 = 8r$, and for all $1 \leq i \leq p$, $n_{i+1} = n_i^{8r}$. Set $t = n_{p+1}$. Then for all $1 \leq i < p$, $(n_i + r + (n_i + r)^{t+r})^r \leq (n_{i+1} - r)^{(t - n_{i+1})}$.

Note that $|([1,n_i+r] \cup [1,n_i+r]^{t+r})^r| = (n_1+r+(n_i+r)^{t+r})^r$. Furthermore, let $[1,n_{i+1}-r]^{t'}$ be the set of all elements of $[1,n_{i+1}-r]^t$ whose set of coordinates is exactly $[1,n_{i+1}-r]$. Note that $(n_{i+1}-r)^{(t-n_{i+1})} \leq |[1,n_{Pi+1}-r]^{t'}|$.

By the above inequality, for each $1 \leq i < p$, there exists a surjective map $h_i\colon [1,n_{i+1}-r]^{t'} \to ([1,n_i+r] \cup [1,n_i+r]^{t+r})^r$.

We now define $h$. Let $x \in N^t$.

*Case* 1. For some $i > 0$, the number of coordinates of $x$ is exactly $n_{i+1}-r$. Let $x'$ be the unique element of $[1,n_{i+1}-r]^t$ which is order isomorphic to $x$. Then define $h(x)$ to be the unique $y$ such that $(x,y)$ is order isomorphic to $(x', h_i(x'))$.

*Case* 2. Otherwise. Define $h(x) = x_1{}^*r$.

Clearly $h$ is basic. Now let $E \subseteq N$ and $y \in (E[n_i+r] \cup E[n_i+r]^{t+r})^r$, where $1 \leq i \leq p$. Then (i) follows by transference of $E$ order isomorphically onto an initial segment of $N$. Also (ii) and (iii) follow from Case 2. ☐

For $x \in N^{t+r}$, we define $\#x$ to be the greatest $i$ such that $x_1,\ldots,x_t < x_{t+1} < \cdots < x_{t+i} = \cdots = x_{t+r}$. If this does not exist, then $\#x$ is undefined. Note that if $\#x$ exists then $1 \leq \#x \leq r$. If $\#x = 1$, then the last $r$ coordinates of $x$ are equal, and if $\#x = r$, then the last $r$ coordinates of $x$ are strictly increasing. Below, we will be interested almost exclusively in only those $x \in N^{t+r}$ for which $\#x$ exists.

If $1 \leq i \leq r$, then define $x\#(i) = (x_1,\ldots,x_{t+i},\ldots,x_{t+i})$.

Let $f \in \mathrm{FPF}(N^{t+r})$. For $x \in N$ define $f(x) = x$. For $x \in (N \cup N^{t+r})^r$ define $f(x) = (f(x_1),\ldots,f(x_r))$.

Let $f \in \mathrm{FPF}(N^{t+r})$ and $h\colon N^t \to (N \cup N^{t+r})^r$. We now define $f^*[h] \in \mathrm{FPF}(N^{t+r}, N^r)$ as follows. Let $x \in N^{t+r}$.

*Case* 1. $\#x$ and $f(x)$ are defined. Define $f^*[h](x) = f(h(x_1,\ldots,x_t))$.

*Case* 2. If otherwise, then $f^*[h](x)$ is undefined.

Let $P$, $Q$ be expressions. It is useful to use the notation $P =' Q$ for $P = Q$ or ($P$ and $Q$ are undefined). Of course, $P = Q$ means that $P$ and $Q$ are defined and equal.

LEMMA 5.15. *Let $t,r > 0$ and $h\colon N^t \to (N \cup N^{t+r})^r$ be basic. Then there exists $B \in \mathrm{BEF}(t+2r, r+1, t+r)$ such that for all $f \in \mathrm{FPF}(N^{t+r})$, $x \in N^{t+r}$, and $y \in N^r$, $f^*[h](x) = y$ if and only if $B(x,y)$ is true in $f$.*

*Proof.* The proof is by inspection of the definition of $f^*[h]$. ☐

For $f \in \mathrm{FPF}(N^k)$ and $S \subseteq N$, we write $f[S]$ for $f[S^k]$.

LEMMA 5.16. *Let $r, p > 0$. Then there exists $t > 0$, $0 < n_1 < \cdots < n_p \le t$, and basic $h\colon N^t \to (N \cup N^{t+r})^r$, such that the following hold. Let $f \in \mathrm{FPF}(N^{t+r})$, $E^{t+r} \subseteq \mathrm{dom}\,(f)$, $|E| \ge n_p$, $1 \le i < p$, and $x, y \in N^{t+r}$.*

(i) $\mathrm{fld}(f^*[h][E[n_i + r]])^r \subseteq (E[n_i + r] \cup f[E[n_i + r]])^r \subseteq f^*[h][E[n_{i+1}]]$;

(ii) *If $\#x, \#y, f(x), f(y)$ exist and $x$, $y$ have the same first $t$ coordinates, then $f^*[h](x) =' f^*[h](y)$;*

(iii) *If $f^*[h](x)$ exists, then $f^*[h](x)$ is an $r$-tuple consisting of certain elements of $\{x_1, \ldots, x_t\}$ together with the value of $f$ at certain $y \subseteq (x_1, \ldots, x_t)$.*

*Proof.* Let $r, p > 0$. Let $t, n_1, \ldots, n_p, h$ be according to Lemma 5.14. Let $f$, $E$, $i$ be as given.

For (i), first observe that $\mathrm{fld}(f^*[h][E[n_i + r]]) \subseteq E[n_i + r] \cup f[E[n_i + r]]$ by inspection. Now let $x_1, \ldots, x_r \in E[n_i + r] \cup f[E[n_i + r]]$. Then there exists $y \in E[n_{i+1} - r]^t$ such that $f(h(y)) = (x_1, \ldots, x_r)$. Let $z \in E[n_{i+1}]^{t+r}$ extend $y$ by the elements of $E$ in positions $n_{i+1} - r + 1, \ldots, n_{i+1}$. Then $\#z$ and $f(z)$ are defined. So $f^*[h](z) = f(h(y)) = (x_1, \ldots, x_r)$, as required.

Claim (ii) is by inspection of the definitions. For (iii), suppose $f^*[h](x)$ exists. Then $f^*[h](x) = f(h(x))$. Also $\#x$ exists, and so $x$ is not constant. Hence $h(x)$ is a coordinate of $x$ or $h(x) \subseteq x = (x_1, \ldots, x_t)$. □

We now fix $r \ge 2$ and $p \ge 9$ until the proof of Lemma 5.33 is complete and we fix the corresponding $t, n_1, \ldots, n_p, h$ given by Lemma 5.16. We will abbreviate $f^*[h]$ by $f^*$.

We now present the crucial definition by induction on the sup norm. The equivalent definition given in Lemma 5.17 is clearer, but is not directly an induction on the sup norm. For this purpose, we fix an order-invariant $R \subseteq N^r \times N^r = N^{2r}$, such that for all $x, y \in N^r$, $R(x, y) \Rightarrow |x| > |y|$.

Let $A$ be a finite subset of $N^{t+r}$. We define $f\colon A \to \mathrm{fld}(A)$ as follows. Let $x \in A$.

*Case* 1. There exists $y$ such that

(i) $|y| < |x|$;

(ii) $f(y) = |y| < |f^*(x)|$;

(iii) $R(f^*(x), f^*(y))$;

(iv) For all $1 \le j < \#x$, $f^*(y)_j = f(x\#(j))$.

Define $f(x) = \min\{f^*(y)_{\#x}\colon y \text{ is as above}\}$.

*Case* 2. Otherwise. Define $f(x) = |x|$.

We now apply Lemmas 5.13 and 5.15. Note that the determination of whether Case 1 applies to $f(x)$ involves only applications of $f$ to tuples of lower sup norm than $|x|$. This follows from Lemma 5.16 (iii).

We summarize the data that will be fixed until the proof of Lemma 5.33 is complete:

(i) $t, r, p, n_1, \ldots, n_p, h, R$ are previously fixed;
(ii) $f\colon A \to \mathrm{fld}(A)$ obeying the preceding inductive definition, where $A$ is a finite closed subset of $N^{t+r}$;
(iii) $E \subseteq N$, $|E| = n_p$;
(iv) $E$ is a $(t+r, t+r)$-SOI for $f$;
(v) $\max(\mathrm{fld}(A)) = \max(E)$.

Clause (iv) will only be used in the proofs of Lemmas 5.26 and 5.31.

Lemma 5.17. *Let $x \in A$. Then $f(x) \leq |x|$. Also $f^*(x)$ is defined if and only if $\#x$ is defined. And if $\#x$ is defined, then $f^*(x) < |x|$. Also $E \subseteq f[E]$. Furthermore, $f$ is the unique function from $A$ into* fld$(A)$ *such that the following hold for all $x \in A$:*
*Case $1'$. There exists $y$ such that*
(i) $f(y) = |y| < |f^*(x)|$;
(ii) $R(f^*(x), f^*(y))$.
*Let $z =$ lex min$\{f^*(y)\colon y$ is as above$\}$. Then $f(x) = z_{\#x}$.*
*Case $2'$. Otherwise. Define $f(x) = |x|$.*

*Proof.* The first claim is by induction on $|x|$, by Lemma 5.16 (iii). The second claim follows from Lemma 5.16 (iii) and the fact that $A$ is closed. The third claim follows from the first two claims and Lemma 5.16(iii). For the fourth claim, note that since $E$ is a $(t+r, t+r)-$SOI for $f$, $E^{t+r} \subseteq A$. Hence for all $j \in E$, $j^*t + r \in A$ and $\#(j^*t + r)$ is undefined. Hence $f^*(j^*t + r)$ is undefined, and therefore $f(j^*t + r) = j \in f[E]$.

For the final claim, let $f$ be as given. Assume Case 1 applies to $f(x)$. Then Case 1 applies to $f(x\#(1))$, according to Lemma 5.16 (ii). And $f(x\#(1))$ is the least possible first coordinate of the associated $f^*(y)$. Also Case 1 applies to $f(x\#(2))$, using the $f^*(y)$'s whose first coordinate has been minimized, and $f(x\#(2))$ is the least possible second coordinate of these $f^*(y)$'s. We continue in this way till we get to $f(x)$. It is clear that we have merely made the first $\#(x)$ steps in the standard lexicographically least element construction on the set of $y$ such that Case $1'$ applies to $f(x\#(1))$ – which is of course the same as the set of $y$ such that Case $1'$ applies to $f(x\#(2))$ and higher, up through $f(x)$. This matches Cases $1'$ and $2'$ in the lemma, when we note that clause (i) in Case 1 follows from clause (i) in Case $1'$. □

From now on we will refer only to clauses $1'$ and $2'$.

We define

$$f\hat{}(x) = \text{lex min}\{f^*(y)\colon f(y) = |y| < |f^*(x)| \text{ and } R(f^*(x), f^*(y))\}$$

if this is nonempty; it is undefined otherwise.

Lemma 5.18. *For all $x \in A$, $f\hat{}(x)$ is defined if and only if $|f(x)| < |x|$. For all $x \in A$, if $f\hat{}(x)$ is defined, then $|f\hat{}(x)| < |f^*(x)|$. Also if $f^*(x) = f^*(y)$, then $f\hat{}(x) =' f\hat{}(y)$.*

*Proof.* This is by Lemma 5.17 and inspection. □

Lemma 5.19. *If $f\hat{}(x)$ is defined and $1 \leq i \leq \#x$, then $f\hat{}(x)_i = f(x\#(i))$. If $f\hat{}(x)$ is defined, then $f(x) = f\hat{}(x)_{\#x}$.*

*Proof.* Assume $f\hat{}(x)$ is defined and $1 \leq i \leq \#x$. Then Case $1'$ applies to $f(x)$, and the $z$ in Case $1'$ is $f\hat{}(x)$. Also Case $1'$ applies to $f(x\#(i))$, and the $z$ in Case $1'$ is again $f\hat{}(x)$. Hence $f(x\#(i)) = z_i$ because $\#(x\#(i)) = i$.

The second claim follows by the identity $x\#(\#x) = x$. □

Lemma 5.20. *For all $x \in A$, $f(x) = f\hat{}(x)_{\#x}$ if $f\hat{}(x)$ is defined; $|x|$ otherwise. If $f\hat{}(x)$ is defined and $\#x = r$, then*

$$f\hat{}(x) = (f(x\#(1)), \ldots, f(x\#(r)))\,.$$

*If $x \in E^{t+r}$ and $\#x$ is defined, then $f^*(x) \in f[E]^r$. If $x \in E^{t+r}$, $f\hat{}(x)$ is defined, and $\#x = r$, then $f\hat{}(x) \in f[E]^r$.*

*Proof.* Let $x \in A$. If Case $1'$ applies to $f(x)$, then $f\hat{}(x)$ is defined, and so by Lemma 5.19, $f\hat{}(x)_{\#x} = f(x)$. If Case $1'$ does not apply to $f(x)$, then $f(x) = |x|$. This establishes the first claim. The second claim is a special case of Lemma 5.19. For the third claim, let $x$ be as given. Since $E^{t+r} \subseteq \operatorname{dom}(f)$, we see that $f^*(x)$ is defined. Also by Lemma 5.16, $f^*(x)$ consists of elements of $E$ and $f[E]$. By Lemma 5.17, $E \subseteq f[E]$. The fourth claim follows from the second claim. □

Lemma 5.21. *If $f^*(x)$ is defined, then $|f^*(x)| < x_{t+1}$. If $f\hat{}(x)$ is defined, then $|f\hat{}(x)| < |f^*(x)|$.*

*Proof.* To establish the first claim, suppose $f^*(x)$ is defined. Then by Lemma 5.16, $f^*(x)$ is an $r$-tuple consisting of elements of $\{x_1, \ldots, x_t\}$ and values of $f$ at certain $y \subseteq (x_1, \ldots, x_t)$. Also $\#x$ is defined. Hence $|f^*(x)| < x_{t+1}$. By Lemma 5.17, $|f^*(x)| \leq \max(x_1, \ldots, x_t) < x_{t+1}$.

For the second claim, suppose $f\hat{}(x)$ is defined. Then $R(f^*(x), f\hat{}(x))$ by the definition of $f\hat{}(x)$. By the hypothesis on $R$, we have $|f\hat{}(x) < |f^*(x)|$. □

Lemma 5.22. *The following are equivalent for $x \in A$:*

(i) *$f\hat{}(x)$ is defined;*
(ii) *Case $1'$ applies to $f(x)$;*
(iii) *$f(x) < |x|$;*
(iv) *$f(x) < |f^*(x)|$.*

*Proof.* (i) and (ii) are equivalent by the definition of $f\hat{}$.

Assume (ii). By Lemma 5.21, the $f^*(y)$ occurring in the lex min in Case $1'$ applied to $f(x)$ must have $|f^*(y)| < |f^*(x)| < x_{t+1}$. Hence (iii) and (iv) hold.

Now obviously (iii) implies (ii). Hence (iii) implies (iv).

Finally, by Lemma 5.17, $|f^*(x)| < |x|$ if defined, and hence (iv) implies (iii). □

For $1 \leq i \leq p$, define $V(i) = f^*[E[n_i]]$. By Lemma 5.16 we have $(E[n_i + r] \cup \text{fld}(V(i)))^r \subseteq V(i+1)$, where $1 \leq i \leq p-1$. We also define the monadic predicate $P$ on $N^r$ by

$$P(x) \text{ if and only if there exists } y \text{ such that}$$
$$f^*(y) = x \text{ and } f(y) < |y|.$$

LEMMA 5.23. *If $P(f^*(x))$, then $f(x) < |x|$.*

*Proof.* Let $P(f^*(x))$. Let $f^*(y) = f^*(x)$ and $f(y) < |y|$. By Lemma 5.18, $f\hat{}(y) =' f\hat{}(x)$. By Lemma 5.22, $f\hat{}(y)$ is defined. Hence $f\hat{}(x)$ is defined, and so by Lemma 5.22, $f(x) < |x|$. □

LEMMA 5.24. *If $f(x) < |x|$, then $\neg P(f\hat{}(x))$.*

*Proof.* Let $f(x) < |x|$. By Lemma 5.22, $f\hat{}(x)$ is defined. By the definition of $f\hat{}$, let $f(y) = |y|$ and $f^*(y) = f\hat{}(x)$. Suppose $P(f\hat{}(x))$. Then by Lemma 5.23, $f(y) < |y|$, which is a contradiction. □

LEMMA 5.25. *Let $1 \leq i \leq p-1$, $x \in V(i)$, and $P(x)$. Then there exists $y \in V(i+1)$ such that $R(x, y)$ and $\neg P(y)$.*

*Proof.* Let $i$, $x$ be as given. Let $f^*(z) = x$, $z \in E[n_i]^{t+r}$. By Lemma 5.23, $f(z) < |z|$. By Lemma 5.24, $\neg P(f\hat{}(z))$. Now by the definition of $f\hat{}$, we have $R(x, f\hat{}(z))$. It remains to show that $f\hat{}(z) \in V(i+1)$. By Lemma 5.20, $f\hat{}(z) \in f[E[n_i]]^r$. By Lemma 5.16, $f[E[n_i]]^r \subseteq f^*[E[n_{i+1}]]$, as required. □

LEMMA 5.26. *Let $x \in V(p-1)$, $1 \leq i \leq n_p$, and $|x| < E_i$. Then there exists $|y| < E_i$ such that $f^*(y) = x$.*

*Proof.* Let $x$ be as given. Let $y' \in E[n_{p-1}]^{t+r}$, $f^*(y') = x$. The true sentence

$$\text{there exists } y \in A \text{ such that } |y| < \max(E) \text{ and } f^*(y) = x$$

can be viewed as a formula in $\text{BEF}(r+1, r+2, t+r)$ about $x$ and $\max(E)$. Now $x \in f[E]^r$. Since $E$ is a $(t+r, t+r)$-SOI for $f$, we can replace $\max(E)$ by $E_i$. □

Lemma 5.27. *Let $x \in \operatorname{rng}(f^*)$ and $|x| \in E$. Suppose there exists $y \in V(p-1)$ such that $R(x,y)$ and $\neg P(y)$. Then $P(x)$.*

*Proof.* Let $x$, $y$ be as given, and let $f^*(v) = x$. Since $R(x,y)$, we have $|y| < |x|$. Hence by Lemma 5.26, let $f^*(w) = y$ and $|w| < |x|$. Since $\neg P(y)$, we have $f(w) = |w|$. Combining this with $R(x,y)$, we see that Case 1′ applies to $f(v)$, and so $f(v) < |v|$. Therefore, $P(x)$. □

Lemma 5.28. *Let $1 \leq i \leq p-2$, $x \in V(i)$, and $|x| \in E$. Then $P(x)$ if and only if there exists $y \in V(i+1)$ such that $R(x,y)$ and $\neg P(y)$.*

*Proof.* The proof is from Lemmas 5.25 and 5.27. □

Lemma 5.29. *Let $1 \leq i \leq p-2$, $x \in V(i)$, $P(x)$, and $|x| \in E$. Then the lexicographically least $y \in V(p-1)$, such that $R(x,y)$ and $\neg P(y)$, exists, and lies in $V(i+1)$.*

*Proof.* Let $i$, $x$ be as given. By Lemma 5.25, this $y$ exists. Let $f^*(w) = x$, $w \in E[n_i]^{t+r}$, as in the proof of Lemma 5.25. There we showed that $R(x, f\hat{}(w))$ and $f\hat{}(w) \in V(i+1)$. It remains to show that $f\hat{}(w)$ is the least $y \in V(p-1)$ such that $R(x,y)$ and $\neg P(y)$. By definition, $f\hat{}(w)$ is the lexicographically least $f^*(u)$ such that $R(x, f^*(u))$ and $f(u) = |u| < |x|$.

Let $y \in V(p-1)$, $R(x,y)$, and $\neg P(y)$. By Lemma 5.26, since $|x| \in E$, let $f^*(u) = y$, $|u| < |x|$. Since $\neg P(y)$, we have $f(u) = |u|$. Hence by the definition of $f\hat{}$, we have $f\hat{}(w) \leq_{\text{lex}} y$, as required. □

Lemma 5.30. *There exists $B \in \text{BEF}(r+1, r+2, t+r)$ such that for all $x \in V(p-1)$, $B(x, \max(E))$ if and only if $P(x)$.*

*Proof.* Let $x \in V(p-1)$. We claim that $P(x)$ if and only if there exists $y \in \operatorname{dom}(f)$ such that $f^*(y) = x$ and $f(y) < |y| < \max(E)$. To see this, let $y' \in \operatorname{dom}(f)$, $f^*(y') = x$, and $f(y') < |y'|$. Now $|x| < \max(E)$. Now by Lemma 5.26, let $|y| < \max(E)$, $f^*(y) = f^*(y') = x$. By Lemma 5.18, $f\hat{}(y) = f\hat{}(y')$. Hence by Lemma 5.22, $f(y') < |y'|$ if and only if $f(y) < |y|$. Therefore $f(y) < |y| < \max(E)$ as required. Using Lemma 5.15, we see that this statement is of the form $B(x, \max(E))$, where $B \in \text{BEF}(r+1, r+2, t+r)$. □

We now give the form of indiscernibility used later.

Lemma 5.31. *Let $1 \leq i \leq r$, $x, y \in E[n_{p-2}]^i$, $z \in \text{fld}(V(p-2))^{r-i}$, where $x, y$ are of the same order type. Let $\pi$ be a coordinate permutation of all $r$-tuples (i.e., there are exactly $r!$ such $\pi$). If $|z| < \min(x,y)$, then $P(\pi(x,z))$ if and only if $P(\pi(y,z))$.*

*Proof.* By Lemma 5.30, $P(\pi(x,y))$ is equivalent to a formula in $\text{BEF}(r+1, r+3, t+r)$ about $x, y$, and $\max(E)$, for $(x,y) \in V(p-1)$. Now let $x, y, z, i$ be as given. Clearly $(x, z), (y, z) \in V(p-1)$. Thus we apply this representation together with the fact that $E$ is a $(t+r, t+r)$-SOI and $y, z \in f[E]^{r-i}$, to obtain the result. □

We now combine Lemmas 5.28, 5.29, and 5.31 into a single concise statement. Recall that we have previously fixed $r \geq 2$ and $p \geq 9$.

LEMMA 5.32. *For all* $1 \leq i \leq p$, *let* $A_i = \text{fld}(V(i))$.

(i) *For all* $1 \leq i < p$, $E[n_i] \subseteq A_{i+1}$;

(ii) *For all* $1 \leq i \leq p-4$ *and* $x \in (A_i)^r$ *with* $|x| \in E$, *we have* $P(x)$ *if and only if there exists* $y \in (A_{i+2})^r$ *such that* $R(x, y)$ *and* $\neg P(y)$;

(iii) *For all* $1 \leq i \leq p-4$ *and* $x \in (A_i)^r$ *with* $|x| \in E$, $P(x)$ *if and only if the lexicographically least* $y \in (A_{p-2})^r$ *such that* $R(x, y)$ *and* $\neg P(y)$ *exists and lies in* $(A_{i+2})^r$;

(iv) *For all* $1 \leq i \leq r$, *and* $x, y \in E[n_{p-2}]^i$ *of the same order type, and* $z \in (A_{p-3})^{r-i}$, *and coordinate permutations* $\pi$ *of* $r$*-tuples, if* $|z| < \min(x, y)$ *and* $|x, y| < \max(E)$, *then* $P(\pi(x, z))$ *if and only if* $P(\pi(y, z))$.

*Proof.* Recall that each $A_i = \text{fld}(V(i)) = \text{fld}(f^*[E[n_i]])$. By Lemma 5.16, for all $1 \leq i < p$, $(E[n_i + r] \cup \text{fld}(V(i))^r \subseteq V(i+1)$. Clearly (i) holds.

For (ii), let $i$, $x$ be as given. Then $x \in V(i+1)$. By Lemma 5.28, if $P(x)$, then there exists $y \in (A_{i+2})^r$ such that $R(x, y)$ and $\neg P(y)$. On the other hand, if there exists $y \in (A_{i+2})^r$ such that $R(x, y)$ and $\neg P(y)$, then there exists $y \in V(i+3)$ such that $R(x, y)$ and $\neg P(y)$, and hence by Lemma 5.28, $P(x)$.

For (iii), let $i$, $x$ be as given. Then $x \in V(i+1)$. Suppose $P(x)$. By Lemma 5.29, the lexicographically least $y \in V(p-1)$ such that $R(x, y)$ and $\neg P(y)$, lies in $V(i+2)$. Now $V(i+2) \subseteq (A_{i+2})^r \subseteq (A_{p-2})^r \subseteq V(p-1)$. Hence this $y$ is the lexicographically least $y \in (A_{p-2})^r$ such that $R(x, y)$ and $\neg P(y)$, and lies in $(A_{i+2})^r$.

And (iv) follows immediately from Lemma 5.31. □

We now streamline Lemma 5.32 as follows.

LEMMA 5.33. *Let* $2n \leq p-4$. *Then there exist nonempty finite sets* $B_1 \subseteq B_2 \subseteq \cdots \subseteq B_n \subseteq N$ *such that the following hold*:

(i) *For all* $1 \leq i \leq n$, $\max(B_i) = E_i$;

(ii) *For all* $1 \leq i < n$ *and* $x \in (B_i)^r$ *with* $|x| \in E$, $P(x)$ *if and only if there exists* $y \in (B_{i+1})^r$ *such that* $R(x, y)$ *and* $\neg P(y)$;

(iii) *For all* $1 \leq i < n$ *and* $x \in (B_i)^r$ *with* $|x| \in E$, $P(x)$ *if and only if the lexicographically least* $y \in (B_n)^r$, *such that* $R(x, y)$ *and* $\neg P(y)$, *exists and lies in* $(B_{i+1})^r$;

(iv) *For all* $1 \le i \le r$, *and* $x, y \in E[n]^i$ *of the same order type, and* $z \in (B_n)^{r-i}$, *and coordinate permutations* $\pi$ *of* $(B_n)^r$, *if* $|z| < \min(x, y)$, *then* $P(\pi(x, z))$ *if and only if* $P(\pi(y, z))$.

*Proof.* Let $2n \le p-4$. For all $1 \le i \le n$, set $B_i = A_{2i} \cap [0, E_i]$. Then (i) is evident. To verify (ii), let $1 \le i \le n$, $x \in (B_i)^r$, and $|x| \in E$. Hence $|x| \le E_i$. Now $(B_i)^r \subseteq (A_{2i})^r$, and so (ii) of Lemma 5.32 applies. So $P(x)$ if and only if there exists $y \in (A_{2i+2})^r$ such that $R(x, y)$ and $\neg P(y)$. Since $R(x, y)$ implies $|x| > |y|$, we can replace $(A_{2i+2})^r$ with $(B_{i+1})^r$.

To verify (iii), let $1 \le i < n$ and $x \in B_i^r$ with $|x| \in E$. Hence $|x| \le E_i$. Then we can apply (iii) of Lemma 5.32 to $2i$. Thus $P(x)$ if and only if the lexicographically least $y \in (A_{p-2})^r$, such that $R(x, y)$ and $\neg P(y)$, exists and lies in $(A_{2i+2})^r$. Again using $R(x, y)$ implies $|x| > |y|$, we see that $P(x)$ if and only if the lexicographically least $y \in (A_{p-2})^r$, such that $R(x, y)$ and $\neg P(y)$, exists and lies in $(B_{i+1})^r$.

Finally, (iv) is immediate from the (iv) of Lemma 5.32. ☐

We make use of the following construction. Let $R$ be an order invariant subset of $N^s$. Let $(A, <)$ be a linear ordering. Then $R$ can also be viewed as an order invariant subset of $A^s$ by taking the set of all elements of $A^s$ which are of the same order type as an element of $R$.

5.3. *Existence of infinite linearly ordered predicates with strong atomic indiscernibles, comprehension for bounded formulas of limited complexity, and a weak least element principle.* The next lemma is an application of the compactness theorem for predicate calculus, and brings us to the major milestone of working within an ordinary relational structure instead of finite towers of finite sets.

LEMMA 5.34. *Let* $k \ge 2$ *and* $R$ *be an order-invariant subset of* $N^k \times N^k$ *such that for all* $x, y \in N^k$, *if* $R(x, y)$ *then* $|x| > |y|$. *Then there exists* $(A, <, E, P)$ *obeying the following conditions*:

(i) $(A, <)$ *is a countable linear ordering*;

(ii) $E = \{e_1 < e_2 < \dots\}$ *is an unbounded subset of* $A$;

(iii) $P \subseteq A^k$;

(iv) *For all* $x \in A^k$ *with* $|x| \in E$, $P(x)$ *if and only if there exists* $y \in A^k$ *such that* $R(x, y)$ *and* $\neg P(y)$;

(v) *For all* $x \in A^k$ *with* $|x| \in E$, $P(x)$ *if and only if there exists a lexicographically least* $y \in A^k$ *such that* $R(x, y)$ *and* $\neg P(y)$;

(vi) *For all* $1 \le i \le k$, *and* $x, y \in E^i$ *of the same order type, and* $z \in A^{k-i}$, *and coordinate permutations* $\pi$ *of* $A^k$, *if* $|z| < \min(x, y)$, *then* $P(\pi(x, z))$ *if and only if* $P(\pi(y, z))$.

*Proof.* By compactness. Fix $k \geq 0$. We use the following symbols:

(a) the binary relation symbol $<$;
(b) monadic predicate symbols $A_1, A_2, \ldots$;
(c) the $k$-ary predicate symbol $P$;
(d) the constant symbols $e_1, e_2, \ldots$.

We consider the following set $T$ of axioms:

(1) $<$ is a linear ordering;
(2) for all $i$: $A_i(x) \Rightarrow A_{i+1}(x)$;
(3) for all $i$: $\max(A_i) = e_i$;
(4) for all $i < j$: $e_i < e_j$;
(5) for all $i, t$:

$$(A_i(x_1) \,\&\, \ldots \,\&\, A_i(x_k) \,\&\, \max(x_1, \ldots, x_k) = e_t)$$
$$\Rightarrow (P(x_1, \ldots, x_k) \Leftrightarrow (\exists\, y_1, \ldots y_k)$$
$$(A_{i+1}(y_1) \,\&\, \ldots \,\&\, A_{i+1}(y_k) \,\&\, R(x_1, \ldots, x_k, y_1, \ldots, y_k)$$
$$\&\, \neg P(y_1, \ldots, y_k)))\,;$$

(6) For all $i, t$:

$$(A_i(x_1) \,\&\, \ldots \,\&\, A_i(x_k) \,\&\, \max(x_1, \ldots, x_k) = e_t)$$
$$\Rightarrow (P(x_1, \ldots, x_k) \Leftrightarrow \text{there exists lexicographically least } (y_1, \ldots, y_k)$$
$$\text{such that } (R(x_1, \ldots, x_k, y_1, \ldots, y_k) \,\&\, \neg P(y_1, \ldots, y_k)))\,;$$

furthermore, this $(y_1, \ldots, y_k)$ has $A_{i+1}(y_1, \ldots, y_k)$;

(7) Let $1 \leq i \leq k$, and $1 \leq j_1, \ldots, j_i, p_1, \ldots, p_i$, where $(j_1, \ldots, j_i)$ and $(p_1, \ldots, p_i)$ have the same order type. Let $\pi$ be a coordinate permutation of $k$-tuples. For all $z_1, \ldots, z_{k-i} < \min(e_{j1}, \ldots, e_{ji}, e_{p1}, \ldots, e_{pi})$,
$P(\pi(e_{j1}, \ldots, e_{ji}, z_1, \ldots, z_{k-i})) \Leftrightarrow P(\pi(e_{p1}, \ldots, e_{pi}, z_1, \ldots, z_{k-i}))$.

We use Lemma 5.33 to show that every finite subset of $T$ has a (finite) model. Fix $T_0$ to be any finite subset of $T$ which uses at most $A_1, \ldots, A_n$ and $e_1, \ldots, e_n$, $n \geq 1$. In Lemma 5.33, set $r = k$, $p = 2n + 7$, and use the same $n$. Fix nonempty finite sets $B_1 \subseteq \cdots \subseteq B_n \subseteq N$ according to Lemma 5.33. Our model of $T_0$ will have universe $B_n$, and $<$ will be the standard ordering on $N$ restricted to $B_n$. The interpretation of $A_1, \ldots, A_n$ will be just $B_1, \ldots, B_n$. The interpretation of $P$ will be the $P$ of Section 5.2. The interpretation of $e_1, \ldots, e_n$ will be $E_1, \ldots, E_n$.

By the compactness theorem, $T$ is satisfiable, and so let

$$(X, <, P, e_1, e_2, \ldots; A_1, A_2, \ldots)$$

be a countably infinite model of $T$. Let $A$ be the union of the $A_i$'s. Note that $A$ may be a proper subset of $X$. Restrict this model down to the submodel $(A, <, P, e_1, e_2, \ldots; A_1, A_2, \ldots)$. Using axiom scheme (3), we see that this is still a model of $T$, and also condition (ii) holds. ∎

We now define $\mathrm{Fm}(j,n,k)$ for $j \geq 0$ and $k > 0$ by induction as follows. Here Fm abbreviates "Formula." Recall the notation $u^*m = u, u, \ldots, u$, where $m \geq 0$.

We let $\mathrm{Fm}(0,n,k)$ be the formulas of the following form:

(i) $y < z$, where $y$, $z$ are among the variables $x_1, \ldots, x_n$;

(ii) $\neg P(y_1, y_2, \ldots, y_{n-1}, y_n{}^*(k-n+1))$, where each $y_i$ is among $x_1, \ldots, x_n$.

Note that we are using $P$ as a $k$-any relation symbol.

In case $n \geq k$, (ii) will simply consist of all formulas of the form $\neg P(y_1, \ldots, y_k)$, where each $y_i$ is among $x_1, \ldots, x_n$.

Let $j \geq 0$. We define $\mathrm{Fm}^+(j,n,k)$ as $\mathrm{Fm}(j,n,k)$ together with the negations of elements of $\mathrm{Fm}(j,n,k)$. Let $D\mathrm{Fm}^+(j,n,k)$ be $\mathrm{Fm}^+(j,n,k)$ together with all disjunctions of distinct elements of $\mathrm{Fm}^+(j,n,k)$. Let $\mathrm{BFm}(j,n,k)$ be $D\mathrm{Fm}^+(j,n,k)$ together with all disjunctions of negations of distinct elements of $D\mathrm{Fm}^+(j,n,k)$. Note that every propositional combination of elements of $\mathrm{Fm}(j,n,k)$ is propositionally equivalent to an element of $\mathrm{BFm}(j,n,k)$, and vice versa.

Finally, we let $\mathrm{Fm}(j+1,n,k)$ be formulas obtained by putting a block of zero or more distinct existential quantifiers in front of elements of $\mathrm{BFm}(j,n,k)$. Note that $\mathrm{BFm}(j,n,k) \subseteq \mathrm{Fm}(j+1,n,k)$. These quantifiers are required to be among $x_i, \ldots, x_n$.

Note that each $\mathrm{Fm}(i,n,k)$ is finite. Furthermore, for all $k, k' \geq n$, $|\mathrm{Fm}(i,n,k)| = |\mathrm{Fm}(i,n,k')|$, and the same equality holds for $\mathrm{Fm}^+$, $D\mathrm{Fm}^+$, and BFm. Let $(A,<,P)$ be given, where $<$ is a linear ordering on the nonempty set $A$, and $P \subseteq A^k$, $k \geq 1$.

Let $Q \in \mathrm{Fm}(i,n,k)$, $y \in A$ and $s$ be a partial variable assignment. That is, $s$ assigns elements of $A$ to zero, some, or all of the free variables in $Q$.

We say that $Q$ holds $\leq y$ if and only if the formula is true under any interpretation of the unassigned free variables by elements of $\{z \in A\colon z \leq y\}$, and where all quantifiers range over $\{z \in A\colon z \leq y\}$. Note that this definition makes sense for any formula $Q$ in the language of $<$, $=$, $P$, and we will use this remark in Lemma 5.51 below.

Most often, we apply this definition to the case where there is no partial variable assignment – i.e., $s$ is empty. We also apply this definition to the case where $s$ is a total variable assignment (i.e., assigns elements of $A$ to every free variable). We will make sure that $s$ can always be inferred from context.

We fix the dimension $k > 0$ until the proof of Lemma 5.47 is complete.

For $p > 0$, let $Y(p)$ be the set of all $x \in N^k$ whose first $p+1$ coordinates are not all equal, and the first term is $|x|$. Note that $p < q$ implies $Y(p) \subseteq Y(q)$.

We say that $R$ is a roi if and only if $R \subseteq N^k \times N^k = N^{2k}$ is order invariant and for all $x, y \in N^k$, $R(x,y) \Rightarrow |x| > |y|$. (Here roi means regressive order invariant.)

We say that $R$ is a $p$-roi if and only if
(i) $R$ is a roi;
(ii) for all $x, y \in N^k$, if $R(x, y)$, then $x \in Y(p)$.

Let $R$ be a roi. We say that $(A, <, E, P)$ is $R$-good if and only if $(A, <, E, P)$ obeys conditions (i)–(vi) in Lemma 5.34, where of course $R$ is reinterpreted as the corresponding order invariant subset of $A^k \times A^k$ as discussed at the end of Section 5.2. Thus Lemma 5.34 asserts that for all roi $R$ there is an $R$-good $(A, <, E, P)$. We say that $(A, <, E, P)$ is $(R, p)$-good if and only if it obeys conditions (i)–(vi) in Lemma 5.34, where in (iv) and (v), we restrict $x$ to be in $Y(p)$.

In the next several lemmas we will build $p$-roi's for steadily increasing $p$. We want to use each roi to help define the next roi. Accordingly, let $R$ be a roi. We say that a roi $S$ is a $p$-extension of $R$ if and only if for all $x, y \in A^k$ with $x \in Y(p)$, we have $R(x, y)$ if and only if $S(x, y)$.

We say that $R$ is 1-empty if and only if for all $x, y \in A^k$, if $x \in Y(1)$ then $\neg R(x, y)$.

LEMMA 5.35. *Let $R$ be 1-empty, $(A, <, E, P)$ be $R$-good, $i > 0$, $x \in A^{k-1}$, and $|x| < e_i$. Then $\neg P(e_i, x)$.*

*Proof.* By (iv) in Lemma 5.34, $P(e_i, x)$ if and only if there exists $y \in A^k$ such that $R(e_i, x, y)$ and $\neg P(y)$. Note that if $|x| < e_i$, then $(e_i, x) \in Y(1)$. Hence if $|x| < e_i$, then $\neg P(e_i, x)$. □

LEMMA 5.36. *Let $n > 1$. There exists a 1-empty $n^2$-roi $R$ such that the following holds*:

*Let $Q$ be of the form $x_a < x_b$, where $1 \leq a, b \leq n$. Then there exists $2 \leq p \leq n^2$ such that for all $i > 1$ and $1 \leq m \leq k - p - n - 1$,*

$$Q \text{ iff } P(e_{i+m}{}^*p, e_{i+m-1}, \ldots, e_i, x_1, \ldots x_n, e_i{}^*)$$

*holds $\leq e_i$ in all $(R, n^2)$-good $(A, <, E, P)$.*

(*Note.* We write $e_i{}^*$ at the end of arguments for $P$ to indicate zero or more copies of $e_i$ in order to fill out the $k$ arguments needed for $P$.)

*Proof.* We first define $R(x, y)$ to be false for $x \in Y(1)$. Note that the $Q$'s in question do not mention $P$.

Now assign a unique index $2 \leq p \leq n^2$ to each $Q$ of the given form. (Normally we would need to use $n^2 + 1$ instead of $n^2$, but we only need to handle one formula of the form $x_a < x_a$.) We can obviously arrange $R$ to ensure that for all $2 \leq p \leq n^2$ and $x_1, \ldots, x_n \leq e_i$,

$$P(e_{i+m}{}^*p, e_{i+m-1}, \ldots, e_i, x_1, \ldots, x_n, e_i{}^*)$$

if and only if there exists $|y| < e_i$ such that $Q \,\&\, \neg P(e_i, y)$, where $Q$ is assigned index $p$. By Lemma 5.35 and $i > 1$, we see that the right-hand side of this equivalence is equivalent to $Q$. □

LEMMA 5.37. *Let $n > 1$. There is a 1-empty $|\mathrm{Fm}(0, n, k)|$-roi $R$ such that the following holds*:

*Let $Q \in \mathrm{Fm}(0, n, k)$. Then there exists $2 \le p \le |\mathrm{Fm}(0, n, k)|$ such that for all $i > 1$ and $1 \le m \le k - p - n - 1$,*

$$Q \text{ iff } P(e_{i+m}{}^*p, e_{i+m-1}, \ldots, e_i, x_1, \ldots, x_n, e_i{}^*)$$

*holds $\le e_i$ in all $(R, |\mathrm{Fm}(0, n, k)|)$-good $(A, <, E, P)$.*

*Proof.* Let $R$ be as in Lemma 5.36, which takes care of the formulas in $\mathrm{Fm}(0, n, k)$ that do not mention $P$. We now define $R'$ to be an $n^2$-extension of $R$ as follows. Assign to each remaining $Q$ in $\mathrm{Fm}(0, n, k)$ a unique index $p$ in the interval $(n^2, |\mathrm{Fm}(0, n, k)|]$. We can arrange for $R'$ to ensure that if $i > 0$, $1 \le m \le k - p - n - 1$, and $|x| \le e_i$, then

$$P(e_{i+m}{}^*p, e_{i+m-1}, \ldots, e_i, x_1, \ldots, x_n, e_i{}^*)$$

if and only if $Q$. Now $R'$ is as required. □

LEMMA 5.38. *Let $n > 1$. There is a 1-empty $|\mathrm{Fm}^+(0, n, k)|$–roi $R$ such that the following holds*:

*Let $Q \in \mathrm{Fm}^+(0, n, k)$. Then there exists $2 \le p \le |\mathrm{Fm}^+(0, n, k)|$ such that for all $i > 1$ and $2 \le m \le k - p - n - 1$,*

$$Q \textit{ if and only if } P(e_{i+m}{}^*p, e_{i+m-1}, \ldots, e_i, x_1, \ldots, x_n, e_i{}^*)$$

*holds $\le e_i$ in all $(R, |\mathrm{Fm}^+(0, n, k)|)$-good $(A, <, E, P)$.*

*Proof.* Let $R$ be as in Lemma 5.37. We define $R'$ to be an $|\mathrm{Fm}(0, n, k)|$-extension of $R$ as follows. We have to take care of $Q = \neg Q'$, where $Q' \in \mathrm{Fm}(0, n, k)$. Let $p$ be used for $Q'$ in Lemma 5.37 and let

$$p' \in (|\mathrm{Fm}(0, n, k)|, |\mathrm{Fm}^+(0, n, k)|]$$

be used for $Q$ here. We can arrange $R'$ to ensure that for all $|x| \le e_i$, $P(e_{i+m-1}{}^*p', e_{i+m}, \ldots, e_i, x, e_i{}^*)$ if and only if $\neg P(e_{i+m-1}{}^*p, e_{i+m-2}, \ldots, e_i, x, e_i{}^*)$ if and only if $\neg Q'$ if and only if $Q$. Here $x = x_1, \ldots, x_n$. □

LEMMA 5.39. *Let $n > 1$. There is a 1-empty $|D\mathrm{Fm}^+(0, n, k)|$–roi $R$ such that the following holds*:

*Let $Q \in D\mathrm{Fm}^+(0, n, k)$. Then there exists $2 \le p \le |D\mathrm{Fm}^+(0, n, k)|$ such that for all $i > 1$ and $3 \le m \le k - p - n - 1$,*

$$Q \text{ if and only if } P(e_{i+m}{}^*p, e_{i+m-1}, \ldots, e_i, x_1, \ldots, x_n, e_i{}^*)$$

*holds* $\le e_i$ *in all* $(R, |D\mathrm{Fm}^+(0, n, k)|)$*-good* $(A, <, E, P)$.

*Proof.* Let $R$ be as in Lemma 5.38. Define $R'$ to be an $|\mathrm{Fm}^+(0, n, k)|$–extension of $R$ as follows. We have to take care of $Q = (B_1$ or $\ldots$ or $B_t)$, where $B_1, \ldots, B_t \in \mathrm{Fm}^+(0, n, k)$. Let $p_1, \ldots, p_t$ be the indices used for $\neg B_1, \ldots, \neg B_t$ in Lemma 5.38. Let $p$ be used for $Q$ here. We can arrange $R'$ to ensure that for all $|x| \le e_i$, $P(e_{i+m}{}^*p, e_{i+m-1}, \ldots, e_i, x, e_i{}^*)$ if and only if there exists $1 \le j \le t$ such that $\neg P(e_{i+m-1}{}^*p_j, e_{i+m-2}, \ldots, e_i, x, e_i{}^*)$ if and only if $Q$. Here $x = x_1, \ldots, x_n$. □

LEMMA 5.40. *Let* $n > 1$. *There is a* 1*-empty* $|\mathrm{BFm}(0, n, k)|$-roi $R$ *such that the following holds*:

*Let* $Q \in \mathrm{BFm}(0, n, k)$. *Then there exists* $2 \le p \le |\mathrm{BFm}(0, n, k)|$ *such that for all* $i > 1$ *and* $4 \le m \le k - p - n - 1$,

$$Q \text{ if and only if } P(e_{i+m}{}^*p, e_{i+m-1}, \ldots, e_i, x_1, \ldots, x_n, e_i{}^*)$$

*holds* $\le e_i$ *in all* $(R, |\mathrm{BFm}(0, n, k)|)$*-good* $(A, <, E, P)$.

*Proof.* This is proved from Lemma 5.39 in the same way that Lemma 5.39 was proved from Lemma 5.38. □

LEMMA 5.41. *Let* $n > 1$. *There is a* 1*-empty* $|\mathrm{Fm}(1, n, k)|$-roi $R$ *such that the following holds*:

*Let* $Q \in \mathrm{Fm}(1, n, k)$. *Then there exists* $2 \le p \le |\mathrm{Fm}(1, n, k)|$ *such that for all* $i > 1$ *and* $5 \le m \le k - p - n - 1$,

$$Q \text{ if and only if } P(e_{i+m}{}^*p, e_{i+m-1}, \ldots, e_i, x_1, \ldots, x_n, e_i{}^*)$$

*holds* $\le e_i$ *in all* $(R, |\mathrm{Fm}(1, n, k)|)$*-good* $(A, <, E, P)$.

*Proof.* Let $R$ be as in Lemma 5.40. Define $R'$ to be a $|\mathrm{BFm}(0, n, k)|$-extension of $R$ as follows. Let $Q$ be $(\exists x_{q1}, \ldots, z_{qd})(Q')$, where

$$Q' \in \mathrm{BFm}(0, n, k)\,, \quad d > 0\,,$$

and $1 \le q_1, \ldots, q_d \le n$ are distinct. Let $p'$ be the index used for $\neg Q'$ (put in proper form) in Lemma 5.40, and $p$ be the index used for $Q$ here. We can arrange $R'$ to ensure that for all $|x| \le e_i$,

$$P(e_{i+m}{}^*p, e_{i+m-1}, \ldots, e_i, x_1, \ldots, x_n, e_i{}^*)$$

if and only if

$$\neg P(e_{i+m-1}{}^*p', e_{i+m-2}, \ldots, e_i, y_1, \ldots, y_n, e_i{}^*)$$

holds for some $y_1, \ldots, y_n \le e_i$ such that for all $j \ne q_1, \ldots, q_d$, $y_j = x_j$. The right-hand side of this equivalence is equivalent to $Q$. □

LEMMA 5.42. *Let $n > 1$ and $r \geq 0$. There is a 1-empty $|\mathrm{BFm}(r,n,k)|$-roi $R$ such that the following holds*:

*Let $Q \in \mathrm{BFm}(r,n,k)$. Then there exists $2 \leq p \leq |\mathrm{BFm}(r,n,k)|$ such that for all $i > 1$ and $4r+4 \leq m \leq k-p-n-1$,*

$$Q \text{ if and only if } P(e_{i+m}{}^{*}p, e_{i+m-1}, \ldots, e_i, x_1, \ldots, x_n, e_i{}^{*})$$

*holds $\leq e_i$ in all $(R, |\mathrm{BFm}(r,n,k)|)$-good $(A,<,E,P)$.*

*Proof.* The proof is by induction on $r$. Argue just as we argued cumulatively to get Lemma 5.41. □

5.4. *Existence of linearly ordered predicates with limit points, strong indiscernibles for bounded formulas of limited arity, and comprehension and the least element principle for bounded formulas.*

LEMMA 5.43. *Let $r, n, n' > 1$, $n' > n$, and $|\mathrm{BFm}(r,n',k)|+4r+n'+6 \leq k$. Then there exists a 1-empty $|\mathrm{BFm}(r,n',k)|$-roi $R$ such that the following hold*:

(i) *Let $Q \in \mathrm{BFm}(r,n,k)$. Then there exists $2 \leq p \leq |\mathrm{BFm}(r,n,k)|$ such that for all $i > 1$,*

$$Q \text{ if and only if } P(e_{i+4r+4}{}^{*}p, e_{i+4r+3}, \ldots, e_i, x_1, \ldots, x_n, e_i{}^{*})$$

*holds $\leq e_i$ in all $(R, |\mathrm{Fm}(r,n,k)|)$-good $(A,<,E,P)$;*

(ii) *Let $Q' \in \mathrm{BFm}(r,n',k)$. Then there exists $p' \leq |\mathrm{BFm}(r,n',k)|$ such that for all $i > 1$,*

$$Q' \text{ if and only if } P(e_{i+4r+4}{}^{*}p', e_{i+4r+3}, \ldots, e_i, x_1, \ldots, x_{n'}, e_i{}^{*})$$

*holds $\leq e_i$ in all $(R, |\mathrm{BFm}(r,n',k)|)$-good $(A,<,E,P)$.*

*Proof.* Let $R$ be a $|\mathrm{BFm}(r,n,k)|$-roi according to Lemma 5.42. Then obviously the above holds without considering $Q'$. We construct the required $|\mathrm{BFm}(r,n',k)|$-roi $R'$ to be a $|\mathrm{BFm}(r,n,k)|$-extension of $R$. Now $R'$ is obtained by repeating the proof of Lemma 5.42, starting with formulas in $\mathrm{BFm}(0,n',k) \backslash \mathrm{BFm}(r,n,k)$, using indices that are greater than $|\mathrm{BFm}(r,n,k)|$. All of $\mathrm{BFm}(r,n',k) \backslash \mathrm{BFm}(r,n,k)$ is handled in this manner. Finally, note that every element of $\mathrm{BFm}(r,n,k)$ is appropriately equivalent to an element of $\mathrm{BFm}(r,n',k) \backslash \mathrm{BFm}(r,n,k)$. □

If Lemma 5.43 holds for $r$, $n$, $n'$, $k$, $R$, then we say that $R$ handles $[\mathrm{BFm}(r,n,k), \mathrm{BFm}(r,n',k)]$.

LEMMA 5.44. *Let $r, n, n' > 0$, $n' > n$, and $|\mathrm{BFm}(r,n',k)|+4r+n'+6 \leq k$. Then there exists a $|\mathrm{BFm}(r,n',k)|+3$-roi $R$ which handles $[\mathrm{BFm}(r,n,k)$, $\mathrm{BFm}(r,n',k)]$ such that for all $(R, |\mathrm{BFm}(r,n',k)|+3)$-good $(A,<,E,P)$:*

(i) $(A,<)$ *is a countable linear ordering with a least element in which every element has a unique successor*;

(ii) $E = \{e_1 < e_2 < \cdots\}$ *is an unbounded set of limit points of* $(A,<)$.

*Proof.* Let $R$ handle $[\mathrm{BFm}(r,n,k), \mathrm{BFm}(r,n',k)]$. We can construct a $|\mathrm{BFm}(r,n',k)|$-extension $R'$ of $R$ such that for all $(R', |\mathrm{BFm}(r,n',k)|+2)$-good $(A,<,E,P)$:

(1) If $R'(x,y)$, then it is not the case that $x$ starts with exactly $|\mathrm{BFm}(r,n',k)|$ $+\ 1$ identical coordinates and no more;

(2) For all $n \geq 1$ and $|x| < e_n$, $\neg P(e_n{}^*|\mathrm{BFm}(r,n',k)|+1, x)$;

(3) For all $z < e_3$, $P(e_3{}^*|\mathrm{BFm}(r,n,',k)|+2, z^*)$ if and only if there exists a $k$-tuple $(z^*|\mathrm{BFm}(r,n',k)|+1, u^*)$ such that $\neg P(z^*|\mathrm{BFm}(r,n',k)|+1, u^*)$ and $u < z$ if and only if there exists a lexicographically least $k$-tuple $(z^*|\mathrm{BFm}(r,n',k)|+1, u^*)$ such that $\neg P(z^*|\mathrm{BFm}(r,n',k)|+2, u^*)$ and $u < z$;

(4) For all $i > 0$ and $u < e_i$, $P(e_{i+1}{}^* \mid \mathrm{BFm}(r,n',k)|+3, e_i, u^*)$ if and only if there exists $u < v < e_i$ such that $\neg P(e_i, v^*)$ if and only if there exists a lexicographically least $k$-tuple $(e_i, v^*)$ such that $\neg P(e_i, v^*)$ and $u < v < e_i$.

Note that (2) follows from (1).

By (2), $\neg P(e_2{}^*|\mathrm{BFm}(r,n',k)|+1, e_1{}^*)$. By setting $z = e_2$ and $u = e_1$ in (3), we see that $P(e_3{}^*|\mathrm{BFm}(r,n',k)|+2, z^*)$. Hence there is a lexicographically least $k$-tuple $(z^*|\mathrm{BFm}(r,n',k)|+1, u^*)$ such that $\neg P(z^*|\mathrm{BFm}(r,n',k)|+1, u^*)$. Hence the $u$ in this $k$-tuple must be the least element of $A$. Now since $P(e_3{}^*|\mathrm{BFm}(r,n',k)|+2, e_2{}^*)$, by the indiscernibility of $E$ (going back to Lemma 5.34 (vi)), we have $P(e_3{}^*|\mathrm{BFm}(r,n',k)|+2, e_1{}^*)$. Hence by (3), there exists $u < e_1$. Thus, the least element of $A$ is $< e_1$.

By 1-emptiness, we see that for all $v < e_i$, $\neg P(e_i, v^*)$. By (4), if $u < e_i$ and there exists $u < v < e_i$, then there is a least $v$ such that $u < v < e_i$. Also observe that if $u < e_i$ and there does not exist $u < v < e_i$, then $\neg P(e_{i+1}{}^* \mid \mathrm{BFm}(r,n,k) \mid +3, e_i, u^*)$. Fix such a $u < e_i$. By the indiscernibility of $E$, we have $\neg P(e_{i+2}{}^* \mid \mathrm{BFm}(r,n,k) \mid +3, e_{i+1}, u^*)$, which is a contradiction. Hence there is no greatest $u < e_i$. So for every $u < e_i$ there exists $u < v < e_i$, and in fact a least $v$ such that $u < v < e_i$. Now recall that $E$ is unbounded in $(A,<)$ by Lemma 5.34. This establishes the lemma. □

We let 0 be the least element of $A$. For any $x < e_i$, $i > 0$, we let $x+$ be the least element of $A$ that is $> x$. The result of iterating the $+$ operation starting with 0 yields elements $0 < 1 < \cdots$ below $e_1$ that form an initial segment of $A$ of order type $\omega$.

Henceforth we will view any element of $N$ as the corresponding element of the initial segment of $A$ of order type $\omega$. In general, we will say that $x$ is a finite element of $A$ if $x$ is in the initial segment of $A$ of order type $\omega$.

We say that $R$ sharply handles $[\mathrm{BFm}(r, n, k), \mathrm{BFm}(r, n', k)]$ if and only if $r, n, n', k, R$ are as in Lemma 5.44.

It is convenient to let $c(r, n) = |\mathrm{BFm}(r, n, k)| + 8r + 2n + 12$. Note that this expression does not depend on $k$ provided $n \leq k$, which we can always assume.

LEMMA 5.45. *Let $r, n > 1$ and $R$ sharply handle $[\mathrm{BFm}(r, n, k), \mathrm{BFm}(r, c(r, n), k)]$. Then there exists $q \leq |B\mathrm{Fm}(r, c(r, n), k)|$ such that the following holds*:

*Let $Q \in \mathrm{BFm}(r, n, k)$. Then there exists $p \leq c(r, n)$ such that for all $i > 1$,*

$$Q \text{ if and only if } P(e_{i+8r+9}{}^*q, e_{i+8r+8}, \ldots, e_i, x_1, \ldots, x_n, p, e_{i+4r+5}{}^*)$$

*holds $\leq e_i$ in all $(R, |\mathrm{BFm}(r, c(r, n))| + 3)$-good $(A, <, E, P)$.*

*Proof.* Let $r, n, R$ be as given. Let $Q'(x_1, \ldots, x_{4r+n+6})$ be the formal expansion of the following assertion:

(1) There exists $t \leq |\mathrm{BFm}(r, n, k)|$ such that $x_{4r+n+6} = t$ and

$$P(x_1{}^*t, x_2, \ldots, x_{4r+n+5}, x_{4r+5}{}^*)\,.$$

The formal expansion is obviously a disjunction (over $t$) of conjunctions of length 2, where the second conjunct of each conjunction is in

$$\mathrm{BFm}(0, c(r, n), k)\,.$$

In the expansion, we require that all bound variables be other than

$$x_1, \ldots, x_{4r+n+6}$$

in order to perform certain predicate calculus substitutions.

Observe that for $t \leq |\mathrm{BFm}(r, n, k)|$, the formula $x_{4r+n+6} > t$ expands into an element of $\mathrm{Fm}(1, c(r, n), k)$, when one uses a block of $t + 1$ existential quantifiers. Hence the formula $x_{4r+n+6} \leq t$ expands into an element of $\mathrm{BFm}(1, c(r, n), k)$. Also the formula $x_{4r+n+6} \geq t$ expands into an element of $\mathrm{Fm}(1, c(r, n), k)$. Hence the formula $x_{4r+n+6} = t$ expands into an element of $\mathrm{BFm}(1, c(r, n), k)$, up to tautological equivalence.

From this we see that the expansion of $Q'$ lies in $\mathrm{BFm}(1, c(r, n), k)$. From now on $Q'$ will indicate this expansion of $Q'$. Note that this expansion is valid in all $(R, |\mathrm{BFm}(r, c(r, n))| + 3)$–good $(A, <, E, P)$ by Lemma 5.44.

Observe that for all $i, j > 0$ and $x_1, \ldots, x_{4r+n+6}$,

(2) $Q'(x_1, \ldots, x_{4r+n+6})$ holds $\leq e_i$ if and only if $Q'(x_1, \ldots, x_{4r+n+6})$ holds $\leq e_j$.

By Lemma 5.43 (ii), choose $q \leq |\mathrm{BFm}(r, c(r, n), k)|$ such that for all $i > 1$,

(3) $Q'(x_1, \ldots, x_{4r+n+6})$ if and only if

$$P(e_{i+4r+4}{}^*q, e_{i+4r+3}, \ldots, e_i, x_1, \ldots, x_{c(r,n)}, e_i{}^*)$$

holds $\leq e_i$ in all $(R, |\mathrm{BFm}(r, c(r, n), k)| + 3)$-good $(A, <, E, P)$.

Hence by predicate calculus manipulations, for all $i > 1$,
(4) $Q'(x_1, \ldots, x_{4r+n+6})$ if and only if

$$P(e_{i+4r+4}{}^*q, e_{i+4r+3}, \ldots, e_i, x_1, \ldots, x_{4r+n+6}, e_i{}^*)$$

holds $\leq e_i$ in all $(R, |\mathrm{BFm}(r, c(r, n), k)| + 3)$-good $(A, <, E, P)$.

Therefore by substituting $i + 4r + 5$ for $i$, we see that for all $i > 1$,
(5) $Q'(x_1, \ldots, x_{4r+n+6})$ if and only if

$$P(e_{i+8r+9}{}^*q, e_{i+8r+8}, \ldots, e_{i+4r+5}, x_1, \ldots, x_{4r+n+6}, e_{i+4r+5}{}^*)$$

holds $\leq e_{i+4r+5}$ in all $(R, |\mathrm{BFm}(r, c(r, n), k)| + 3)$-good $(A, <, E, P)$.

Hence by substituting $e_{i+4r+4}, \ldots, e_i$ for $x_1, \ldots, x_{4r+5}$, we have that for all $i > 1$,
(6) $Q'(e_{i+4r+4}, \ldots, e_i, x_{4r+6}, \ldots, x_{4r+n+6})$ if and only if

$$P(e_{i+8r+9}{}^*q, e_{i+8r+8}, \ldots, e_{i+4r+5}, e_{i+4r+4}, \ldots, e_i, x_{4r+6}, \ldots, x_{4r+n+6}, e_{i+4r+5}{}^*)$$

holds $\leq e_{i+4r+5}$ in all $(R, |\mathrm{BFm}(r, c(r, n), k)| + 3)$-good $(A, <, E, P)$.

Therefore, by substituting $x_1, \ldots, x_{n+1}$ for $x_{4r+6}, \ldots, x_{4r+n+6}$, we see that for all $i > 1$,
(7) $Q'(e_{i+4r+4}, \ldots, e_i, x_1, \ldots, x_{n+1})$ if and only if

$$P(e_{i+8r+9}{}^*q, e_{i+8r+8}, \ldots, e_i, x_1, \ldots, x_{n+1}, e_{i+4r+5}{}^*)$$

holds $\leq e_{i+4r+5}$ in all $(R, |\mathrm{BFm}(r, c(r, n), k)| + 3)$-good $(A, <, E, P)$.

By (2), we have
(8) $Q'(e_{i+4r+4}, \ldots, e_i, x_1, \ldots x_{n+1})$ if and only if

$$P(e_{i+8r+9}{}^*q, e_{i+8r+8}, \ldots, e_i, x_1, \ldots, x_{n+1}, e_{i+4r+5}{}^*)$$

holds $\leq e_i$ in all $(R, |\mathrm{BFm}(r, c(r, n), k)| + 3)$-good $(A, <, E, P)$.

By Lemma 5.43, choose $p \leq |\mathrm{BFm}(r, n, k)| < c(r, n)$ so that for all $i > 1$,
(9) $Q$ if and only if $P(e_{i+4r+4}{}^*p, e_{i+4r+3}, \ldots, e_i, x_1, \ldots, x_n, e_i{}^*)$ holds $\leq e_i$ in all $(R, |\mathrm{BFm}(r, c(r, n), k)| + 3)$-good $(A, <, E, P)$.

Now by the definition of $Q'$, we see that for all $i > 1$,
(10) $Q'(e_{i+4r+4}, \ldots, e_i, x_1, \ldots, x_n, p)$ if and only if there exists $t \leq |\mathrm{BFm}(r, n, k)|$ such that $p = t$ and

$$P(e_{i+4r+4}{}^*t, e_{i+4r+3}, \ldots, e_i, x_1, \ldots, x_n, e_i{}^*)$$

holds $\leq e_i$ in all $(R, |\mathrm{BFm}(r, c(r, n), k)| + 3)$-good $(A, <, E, P)$.

Hence by (9) and (10), for all $i > 1$,
(11) $Q'(e_{i+4r+4}, \ldots, e_i, x_1, \ldots x_n, p)$ if and only if $Q$ if and only if

$$P(e_{i+4r+4}{}^*p, e_{i+4r+3}, \ldots, e_i, x_1, \ldots, x_n, e_i{}^*)$$

holds $\leq e_i$ in all $(R, |\mathrm{BFm}(r, c(r, n), k)| + 3)$-good $(A, <, E, P)$.

Also by (8), for all $i > 1$,

(12) $Q'(e_{i+4r+4}, \ldots, e_i, x_1, \ldots, x_n, p)$ if and only if

$$P(e_{i+8r+9}{}^*q, e_{i+8r+8}, \ldots, e_i, x_1, \ldots, x_n, p, e_{i+4r+5}{}^*)$$

holds $\leq e_i$ in all $(R, |\mathrm{BFm}(r, c(r,n), k)| + 3)$-good $(A, <, E, P)$.

Hence by (11) and (12), for all $i > 1$,

(13) $Q$ if and only if

$$P(e_{i+8r+9}{}^*q, e_{i+8r+8}, \ldots, e_i, x_1, \ldots x_n, p, e_{i+4r+5}{}^*)$$

holds $\leq e_i$ in all $(R, |\mathrm{BFm}(r, c(r,n), k)| + 3)$-good $(A, <, E, P)$, as required. □

We streamline the information contained in Lemma 5.45. To this end, we let $Z(p)$ for $p > 0$ be the set of $k$-tuples such that the first $p+1$ coordinates are not all equal, and the first coordinate is not the maximum coordinate.

Let $p, q > 0$. We say that $R$ is a $(p,q)$-roi if and only if

(i) $R \subseteq A^k \times A^k$;

(ii) $R$ is order invariant as a subset of $A^{2k}$;

(iii) $R(x,y)$ implies $x \in Y(p) \cup Z(q)$.

We say that $(A, <, E, P)$ is $(R, p, q)$-good if and only if it obeys conditions (i)–(vi) in Lemma 5.34, where in (iv) and (v), we restrict $x$ to be in $Y(p) \cup Z(q)$.

Lemma 5.46. *Let $r, n > 1$ and $|\mathrm{BFm}(r, c(r,n), k)| + 8r + 2n + 12 \leq k$. Then there exists a 1-empty $(|\mathrm{BFm}(r, c(r,n), k)| + 3, 2)$-roi $R$ such that*

(a) *For all $(R, |\mathrm{BFm}(r, c(r,n), k)| + 3, 2)$-good $(A, <, E, P)$, $(A, <)$ is a countable linear ordering with a least element, in which every element of $A$ has an immediate successor, and $E$ is an unbounded set of limit points in $A$ of order type $\omega$;*

(b) *Let $Q \in \mathrm{BFm}(r, n, k)$. Then there exists $p \leq c(r,n)$ such that for all $i > 1$,*

$$Q \text{ if and only if } P(e_i, \ldots, e_{i+8r+11}, p, x_1, \ldots, x_n, e_i^*)$$

*holds $\leq e_i$ in all $(R, |\mathrm{BFm}(r, c(r,n), k)| + 3, 2)$-good $(A, <, E, P)$.*

*Proof.* By Lemma 5.44, let $R$ be a $|\mathrm{BFm}(r, c(r,n), k)| + 3$-roi that sharply handles $[\mathrm{BFm}(r, n, k), \mathrm{BFm}(r, c(r,n), k)]$. We will construct a $|\mathrm{BFm}(r, c(r,n), k)| + 3$-extension $R'$ of $R$ as required. Since $R'$ will be a $|\mathrm{BFm}(r, c(r,n), k)| + 3$-extension of $R$, claim (a) is immediate from Lemma 5.44.

We have only to replace

$$P(e_{i+8r+9}{}^*q, e_{i+8r+8}, \ldots, e_i, x_1, \ldots, x_n, p, e_{i+4r+5}{}^*)$$

in Lemma 5.45 by $P(e_i, \ldots, e_{i+8r+11}, p, x_1, \ldots, x_n, e_i{}^*)$. But we can arrange $R'$ to ensure that for all $|x, p| \leq e_i$,

$$P(e_i{}^*2, e_{i+1}, \ldots, e_{i+8r+10}, p, x_1, \ldots, x_n, e_i{}^*)$$

if and only if $\neg P(e_{i+8r+9}{}^*q, e_{i+8r+8}, \ldots, e_i, x_1, \ldots, x_n, p, e_{i+4r+5}{}^*)$. And then we can also arrange $R'$ to ensure that if $|x, p| \leq e_i$, then

$$P(e_i, \ldots, e_{i+8r+11}, p, x_1, \ldots, x_n, e_i{}^*)$$

if and only if

$$\neg P(e_i{}^*2, e_{i+1}, \ldots, e_{i+8r+10}, p, x_1, \ldots, x_n, e_i{}^*)$$

if and only if

$$P(e_{i+8r+9}{}^*q, e_{i+8r+8}, \ldots, e_i, x_1, \ldots, x_n, p, e_{i+4r+5}{}^*)\,. \qquad \square$$

We now wish to sharpen Lemma 5.46 in order to include a transfinite induction principle.

LEMMA 5.47. *Let $r, n > 1$ and $|\mathrm{BFm}(r, c(r, 2n), k)| + 8r + 4n + 12 \leq k$. Then there exists a 1-empty $(|\mathrm{BFm}(r, c(r, 2n), k)| + 3, 3)$-roi $R$ such that*

(a) *For all $(R, |\mathrm{BFm}(r, c(r, 2n), k)| + 3, 3)$-good $(A, <, E, P)$, $(A, <)$ is a countable linear ordering with a least element, in which every element of $A$ has an immediate successor, and $E$ is an unbounded set of limit points in $A$ of order type $\omega$;*

(b) *Let $Q \in \mathrm{BFm}(r, 2n, k)$. Then there exists $p \leq c(r, 2n)$ such that for all $i > 1$,*

$$Q \text{ if and only if } P(e_i, \ldots, e_{i+8r+11}, p, x_1, \ldots, x_{2n}, e_i{}^*)$$

*holds $\leq e_i$ in all $(R, |\mathrm{BFm}(r, c(r, 2n), k)| + 3, 3)$-good $(A, <, E, P)$;*

(c) *The following holds in all $(R, |\mathrm{BFm}(r, c(r, 2n), k)| + 3, 3)$-good $(A, <, E, P)$: Let $Q \in \mathrm{BFm}(r, 2n, k)$, $i > 1$, $x \in A^n$, and $|x| \leq e_i$. If there exists $y \in A^n$, $|y| \leq e_i$, such that $Q(x, y)$ holds $\leq e_i$, then there is a lexicographically least $y \in A^n$ such that $|y| \leq e_i$ and $Q(x, y)$ holds $\leq e_i$.*

*Proof.* Let $R$ be according to Lemma 5.46 with $n$ replaced by $2n$. Then $R$ is a 1-empty $(|\mathrm{BFm}(r, c(r, 2n), k)| + 3, 2)$–roi. We will construct a $(|\mathrm{BFm}(r, c(r, 2n, k)| + 3, 2)$–extension $R'$ of $R$ as required.

According to Lemma 5.34 (iv), (v), we can arrange $R'$ such that the following holds. Let $(A, <, E, P)$ be $(R', |\mathrm{BFm}(r, c(r, 2n), k)| + 3, 3)$–good. For all $i > 1$ and $p, x_1, \ldots, x_n \leq e_i$,

(1) $P(e_i{}^*3, e_{i+1}, \ldots, e_{i+8r+12}, p, x_1, \ldots, x_n, e_i{}^*)$ if and only if there exists $x_{n+1}, \ldots, x_{2n} \in A$, $x_{n+1}, \ldots, x_{2n} \leq e_i$, such that $\neg P(e_i, \ldots, e_{i+8r+11}, p, x_1, \ldots, x_{2n}, e_i{}^*)$ if and only if there exists lexicographically least $x_{n+1}, \ldots, x_{2n} \in A$, $x_{n+1}, \ldots, x_{2n} \leq e_i$, such that $\neg P(e_i, \ldots, e_{i+8r+11}, p, x_1, \ldots, x_{2n}, e_i{}^*)$.

Now let $Q \in \mathrm{BFm}(r, 2n, k)$. By Lemma 5.46, let $p \leq e_1$ be such that the following holds. Let $i > 1$, $x_1, \ldots, x_{2n} \in A$, $x_1, \ldots, x_{2n} \leq e_i$. Then

(2) $\neg Q(x_1, \ldots, x_{2n})$ if and only if $P(e_i, \ldots, e_{i+8r+11}, p, x_1, \ldots, x_{2n}, e_i{}^*)$ holds $\leq e_i$ in all $(R', |\mathrm{BFm}(r, c(r, 2n), k)| + 3, 2)$–good $(A, <, E, P)$. Observe that (c) follows immediately from (1) and (2). □

We have arrived at the point where we can simplify matters by fixing $(A, <, E, P)$ with certain properties, and ignoring all that has gone on before. This simplification process ends with Lemma 5.50.

In the next lemma, we remove the $e_i{}^*$ at the end of the arguments for $P$ that are used in Lemma 5.47. In any linear ordering $(A, <)$, we say that $p$ is finite if and only if $p$ is in the initial segment of $A$ of order type $\omega$.

LEMMA 5.48. *Let $r, m > 1$ and $n > 8r+m+13$. There exists $(A, <, E, P)$ satisfying the following conditions*:

(i) *$(A, <)$ is a countable linear ordering with a least element in which every element has an immediate successor*;

(ii) *$E = \{e_1 < e_2 < \cdots\}$ is an unbounded set of limit points of $A$*;

(iii) *$P \subseteq A^{8r+m+13}$*;

(iv) *Let $Q \in \mathrm{BFm}(r, 2n, 8r + m + 13)$ have at most the free variables $x_1, \ldots, x_m$. Then there exists finite $p$ such that for all $i > 1$,*

$$Q \textit{ if and only if } P(e_i, \ldots, e_{i+8r+11}, p, x_1, \ldots x_m)$$

*holds $\leq e_i$*;

(v) *Let $Q \in \mathrm{BFm}(r, 2n, 8r + m + 13)$, $i > 1$, $x \in A^n$, and $|x| \leq e_i$. If there exists $y \in A^n$, $|y| \leq e_i$, such that $Q(x, y)$ holds $\leq e_i$, then there is a lexicographically least $y \in A^n$ such that $|y| \leq e_i$ and $Q(x, y)$ holds $\leq e_i$*;

(vi) *$E$ is a strong set of atomic indiscernibles for $P$. That is, for all $1 \leq i \leq 8r+m+13$, and $x, y \in E^i$ of the same order type, and $z \in A^{8r+m+13-i}$, and coordinate permutations $\pi$ of $A^{8r+m+13}$, if $|z| < \min(x, y)$, then $P(\pi(x, z))$ if and only if $P(\pi(y, z))$.*

*Proof.* Let $k = \max(8r + m + 13, |\mathrm{BFm}(r, c(r, 2n), k)| + 8r + 4n + 12) = \max(8r + m + 13, |\mathrm{BFm}(r, c(r, 2n), c(r, 2n))| + 8r + 4n + 12)$. By Lemma 5.47 and Lemma 5.34, let $(A, <, E, P)$, where $P \subseteq A^k$, be such that the following holds:

(a) $(A, <)$ is a countable linear ordering with a least element in which every element has an immediate successor, and $E$ is an unbounded set of limit points in $A$ of order type $\omega$;

(b) Let $Q \in \mathrm{BFm}(r, 2n, k)$ have at most the free variables $x_1, \ldots, x_m$. Then there exists finite $p$ such that for all $i > 1$,

$$Q \text{ if and only if } P(e_i, \ldots, e_{i+8r+11}, p, x_1, \ldots, x_{2n}, e_i^*)$$

holds $\leq e_i$ in $(A, <, E, P)$;

(c) Let $Q \in \mathrm{BFm}(r, 2n, k)$, $i > 1$, $x \in A^n$, and $|x| \leq e_i$. If there exists $y \in A^n$, $|y| \leq e_i$, such that $Q(x, y)$ holds $\leq e_i$, then there is a lexicographically least $y \in A^n$ such that $|y| \leq e_i$ and $Q(x, y)$ holds $\leq e_i$;

(d) $E$ is a strong set of atomic indiscernibles for $P$. That is, for all $1 \leq i \leq k$, and $x, y \in E^i$ of the same order type, and $z \in A^{k-i}$, and coordinate permutations $\pi$ of $E^k$, if $|z| < \min(x, y)$, then $P(\pi(x, z))$ if and only if $P(\pi(y, z))$.

Now let $P' \subseteq A^{8r+m+13}$ be defined by $P'(y_1, \ldots, y_{8r+m+13})$ if and only if $P(y_1, \ldots, y_{8r+m+13}, y_1{}^*)$.

Let $Q' \in \mathrm{BFm}(r, 2n, 8r + m + 13)$. We claim that there exists $Q \in \mathrm{BFm}(r, 2n, k)$ such that $Q$ has the same free variables, and for all $i > 1$ and all assignments $\leq e_i$ to the free variables,

$$Q \text{ holds } \leq e_i \text{ in } (A, <, P) \text{ if and only if}$$
$$Q' \text{ holds } \leq e_i \text{ in } (A, <, P').$$

To see this, just let $Q$ be the result of replacing each atomic subformula $P(z_1, \ldots, z_{8r+m+13})$ of $Q'$ by $P(z_1, \ldots, z_{8r+m+13}, z_1{}^*)$.

For (iv), let $Q' \in \mathrm{BFm}(r, 2n, 8r + m + 13)$ have at most the free variables $x_1, \ldots, x_m$.

By (b), there exists finite $p$ such that for all $i > 1$,

$$Q \text{ if and only if } P(e_i, \ldots, e_{i+8r+11}, p, x_1, x_2, \ldots, x_{2n}, e_i{}^*)$$

holds $\leq e_i$ in $(A, <, E, P)$.

Since the free variables in $Q$ are among $x_1, \ldots, x_m$, we see that by predicate calculus manipulations, there exists finite $p$ such that for all $i > 1$,

$$Q \text{ if and only if } P(e_i, \ldots, e_{i+8r+11}, p, x_1, \ldots, x_m, e_i{}^*)$$

holds $\leq e_i$ in $(A, <, E, P)$.

Therefore, by the displayed equivalence, there exists finite $p$ such that for all $i > 1$,

$$Q' \text{ if and only if } P'(e_i, \ldots, e_{i+8r+11}, p, x_1, \ldots, x_m)$$

holds $\leq e_i$ in $(A, <, E, P')$. This establishes (iv).

For (v), let $Q' \in \mathrm{BFm}(r, 2n, 8r + m + 13)$, $i > 1$, $x \in A^n$, and $|x| \leq e_i$. Construct $Q \in \mathrm{BFm}(r, 2n, k)$ as above. Then (c) holds of $Q$ in $(A, <, P)$. We need to check that (v) holds of $Q'$ in $(A, <, P')$. Suppose there exists $y \in A^n$, $|y| \leq e_i$, such that $Q'(x, y)$ holds $\leq e_i$ in $(A, <, P')$. Then there exists $y \in A^n$, $|y| \leq e_i$, such that $Q(x, y)$ holds $\leq e_i$ in $(A, <, P)$. Hence there is a lexicographically least $y \in A^n$ such that $|y| \leq e_i$ and $Q(x, y)$ holds $\leq e_i$ in $(A, <, P)$. This same $y$ must be the lexicographically least $y \in A^n$ such that $|y| \leq e_i$ and $Q'(x, y)$ holds $\leq e_i$ in $(A, <, P')$.

For (vi), note that $E$ is a strong set of atomic indiscernibles for $P$. Now let $1 \leq i \leq 8r+m+13$, and $x, y \in E^i$ of the same order type, and $z \in A^{8r+m+13-i}$, and $\pi$ be a coordinate permutation of $A^{8r+m+13}$, and $|z| < \min(x, y)$. It suffices to verify that $P(\pi(x, z), \pi(x, z)_1{}^*)$ if and only if $P(\pi(y, z), \pi(y, z)_1{}^*)$.

First assume that $\pi(x, z)_1$ is some $z_j$. Then $\pi(y, z)_1 = \pi(y, z)_1 = z_j$. Note that $P(\pi(x, z), z_j{}^*)$ if and only if $P(\pi(y, z), z_j{}^*)$ easily follows from the strong atomic indiscernibility of $E$ for $P$.

Finally assume that $\pi(x, z)_1$ is some $x_q$. Then $\pi(y, z)_1 = y_q$. But

$$P(\pi(x, z), x_q{}^*)$$

if and only if $P(\pi(y, z), y_q{}^*)$ also easily follows from the strong atomic indiscernibility of $E$ for $P$. □

We now apply compactness to make Lemma 5.48 more usable. We define $\mathrm{BFm}(r, \infty, k)$ as the union over $n$ of $\mathrm{BFm}(r, n, k)$.

Lemma 5.49. *Let $r, m > 1$. There exists $(A, <, E, P)$ satisfying the following conditions*:

(i) $(A, <)$ *is a linear ordering with a least element in which every element has an immediate successor*;

(ii) $E = \{e_1 < e_2 < \cdots\}$ *is an unbounded set of limit points in* $A$;

(iii) $P \subseteq A^{8r+m+13}$;

(iv) *Let* $Q \in \mathrm{BFm}(r, \infty, 8r + m + 13)$ *have at most the free variables* $x_1, \ldots, x_m$. *Then there exists* $u < e_1$ *such that for all* $i > 1$,

$$Q \text{ if and only if } P(e_i, \ldots, e_{i+8r+11}, u, x_1, \ldots, x_m)$$

*holds* $\leq e_i$;

(v) *Let* $i, n > 1$, $Q \in \mathrm{BFm}(r, \infty, 8r + m + 13)$ *with at most the free variables* $x_1, \ldots, x_{2n}$, $x \in A^n$, *and* $|x| \leq e_i$. *If there exists* $y \in A^n$, $|y| \leq e_i$, *such that* $Q(x, y)$ *holds* $\leq e_i$, *then there is a lexicographically least* $y \in A^n$ *such that* $|y| \leq e_i$ *and* $Q(x, y)$ *holds* $\leq e_i$;

(vi) $E$ *is a strong set of atomic indiscernibles for* $P$. *That is*, *for all* $1 \leq i \leq 8r + m + 13$, *and* $x, y \in E^i$ *of the same order type*, *and* $z \in A^{8r+m+13-i}$, *and* $x, y \in E^i$ *of the same order type*, *and* $z \in A^{8r+m+13-i}$, *and coordinate permutations* $\pi$ *of* $E^{8r+m+13}$, *if* $|z| < \min(x, y)$, *then* $P(\pi(x, z))$ *if and only if* $P(\pi(y, z))$.

*Proof.* Apply compactness to the integer parameter $n$ in 5.48. In order to make Lemma 5.48 first order, we weaken "finite $p$" in (vi) to "$p < e_1$." We obtain (i)–(vi) above without $E$ being unbounded in $A$. We then cut down to the initial segment determined by $E$. □

For ease of applications, we need a modified form of (v) and a sharper form of (vi).

Below, it is convenient to distinguish the parameters from the free variables in a formula. The parameters in a formula are just the free variables which have been assigned a particular object.

LEMMA 5.50. *Let $r, m > 1$. There exists $(A, <, E, P)$ satisfying the following conditions*:

(i) *$(A, <)$ is a linear ordering with a least element in which every element has an immediate successor*;

(ii) *$E = \{e_1 < e_2 < \cdots\}$ is an unbounded set of limit points in $A$*;

(iii) *$P \subseteq A^{8r+m+13}$*;

(iv) *Let $Q \in \mathrm{BFm}(r, \infty, 8r + m + 13)$ have at most the free variables $x_1, \ldots, x_m$. Then there exists $u < e_1$ such that for all $i > 0$,*

$$Q \text{ if and only if } P(e_i, \ldots, e_{i+8r+11}, u, x_1, \ldots, x_m)$$

*holds $\leq e_i$*;

(v) *Let $i > 0$ and $Q \in \mathrm{BFm}(r, \infty, 8r+m+13)$ have one free variable and any number of parameters. If there exists $x \leq e_i$, such that $Q(x)$ holds $\leq e_i$, then there is a least $x$ such that $Q(x)$ holds $\leq e_i$*;

(vi) *$E$ is a strong set of $m$-ary indiscernibles for $\mathrm{BFm}(r, \infty, 8r+m+13)$. That is, let $Q(x_1, \ldots, x_m) \in \mathrm{BFm}(r, \infty, 8r + m + 13)$, where all free variables are shown. Let $0 \leq j \leq m$, $i > 0$, $x, y \in E^j$ be of the same order type, $z \in A^{m-j}$, $|z| < \min(x, y)$, and $|x, y| < e_i$. Then $Q(x, z)$ holds $\leq e_i$ if and only if $Q(y, z)$ holds $\leq e_i$.*

*Proof.* Let $(A, <, E, P)$ be as in Lemma 5.49, except that we throw away $\min(E)$. This enables us to use $i > 0$ in (vi) and (v) instead of $i > 1$.

To obtain (v), let $Q(u_1, \ldots, u_n, x) \in \mathrm{BFm}(r, \infty, 8r + m + 13)$, where $u_1, \ldots, u_n \leq e_j$ are parameters, $x$ is a variable, and $i < j$. Assume that there exists $x \leq e_i$ such that $Q(u_1, \ldots, u_n, x)$ holds $\leq e_i$. We can find $Q' \in \mathrm{BFm}(r, \infty, 8r + m + 13)$ with $n + 2$ variables and no parameters such that for all $x \leq e_i$,

$$Q'(e_i, u_1, \ldots, u_n, x) \text{ holds} \leq e_j \text{ if and only if}$$
$$Q(u_1, \ldots, u_n, x) \text{ holds} \leq e_i.$$

By adding $n$ dummy variables to $x$, we see by Lemma 5.49 (v) that there is a least $x$ such that $Q'(e_i, u_1, \ldots, u_n, x)$ holds $\leq e_j$. Hence there is a least $x$ such that $Q(u_1, \ldots, u_n, x)$ holds $\leq e_i$.

Now let $Q$, $j$, $i$, $x$, $y$, $z$ be as given in (vi). By (iv), let $u < e_1$ be such that

$$Q(x_1, \ldots, x_m) \text{ if and only if}$$
$$P(e_i, \ldots, e_{i+8r+11}, u, x_1, \ldots, x_m)$$

holds $\leq e_i$.

Then

$$Q(x, z) \text{ holds} \leq e_i \text{ if and only if}$$
$$P(e_i, \ldots, e_{i+8r+11}, u, x, z)$$

and

$$Q(y, z) \text{ holds} \leq e_i \text{ if and only if}$$
$$P(e_i, \ldots, e_{i+8r+11}, u, y, z).$$

Hence by the strong atomic indiscernibility in Lemma 5.49 (vi), we have

$$Q(x, z) \text{ holds} \leq e_i \text{ if and only if}$$
$$Q(y, z) \text{ holds} \leq e_i,$$

which completes the proof. □

5.5. *Existence of pairing functions defined by bounded formulas.* We now introduce some new complexity classes which are more convenient than the BFm we have been using.

The formulas under consideration are all based on the $k$-ary predicate symbol $P$ and $<, =$. No constant symbols are allowed. We use the infinite collection of variables, $x_1, x_2, \ldots$. We inductively define the complexity classes $\mathrm{BF}(i, k)$ and $\mathrm{PBF}(i, k)$, for $i \geq 0$ and $k > 0$, as follows. (Here BF means "bounded formula" and PBF means "propositional combination of bounded formulas."

$\mathrm{BF}(0, k)$ consists of the atomic formulas in this language.

$\mathrm{BF}(i + 1, k)$ consists of all formulas of the form

$$(\exists\, y_1, \ldots y_p)(y_1 \leq z_1 \,\&\, \ldots \,\&\, y_p \leq z_p \,\&\, B)$$

where $p \geq 0$, $B \in \mathrm{PBF}(i, k)$, and $y_1, \ldots y_p$, $z_1, \ldots z_p$ are variables, and no $z$ is the same as any $y$.

$\mathrm{PBF}(i, k)$ consists of all propositional combinations of elements of $\mathrm{BF}(i, k)$.

Note that every element of $\mathrm{PBF}(i, k)$ is almost an element of $\mathrm{BFm}(i, \infty, k)$. The minor difference is that we allow all propositional combinations of elements of $\mathrm{BF}(i, k)$. For $B \in \mathrm{PBF}(i, k)$, we can normalize the propositional combinations used in the construction of $B$ so that the variant $B^*$ is literally in $\mathrm{BFm}(i, k)$. Note that $B$ and $B^*$ are logically equivalent in the sense of first order predicate calculus.

Let $B \in \mathrm{PBF}(i, k)$ and $(A, <, P)$ be given, $P \subseteq A^k$. Assume that all free variables of $B$ have been assigned elements of $A$ that are all $\leq x$. We say that $B$ holds $\leq x$ if and only if $B$ holds when all quantifiers are restricted to $\{y \in A \colon y \leq x\}$.

LEMMA 5.51. *Let $i \geq 0$, $k > 0$, $B \in \mathrm{PBF}(i, k)$. Let $(A, <, P)$ be such that $(A, <)$ is a linear ordering and $P \subseteq A^k$. Let $x \in A$ and $s$ be an assignment to the free variables of $B$, where all values of $s$ are $\leq x$. Then*

*$B$ holds in $(A,<,P)$ under $s$ if and only if*
*$B$ holds $\leq x$ in $(A,<,P)$ under $s$ if and only if*
*$B^*$ holds $\leq x$ in $(A,<,P)$ under $s$.*

*Proof.* By induction on $B$. We need only be concerned with the first two displayed statements since $B$ and $B^*$ are logically equivalent.

Assume this is true for $B$. We must show that it is true for $C=(\exists\, y_1,\dots y_p)(y_1 \leq z_1 \,\&\, \dots \,\&\, y_p \leq z_p \,\&\, B)$, where no $z$ is the same as any $y$. Let $s$ be an assignment to the free variables of $C$, where all values of $s$ are $\leq x$. Then

$C$ holds in $(A,<,P)$ under $s$ if and only if
$(\exists\, u_1,\dots u_p)(u_1 \leq s(z_1) \,\&\, \dots \,\&\, u_p \leq s(z_p) \,\&\, B$ holds in
$(A,<,E,P)$ under $s$ adjusted by $\{(y_1,u_1),\dots,(y_p,u_p)\})$.

The right hand side is equivalent to

$(\exists\, u_1 \dots u_p \leq x)(u_1 \leq s(z_1) \,\&\, \dots \,\&\, u_p \leq s(z_p) \,\&\, B$ holds $\leq x$
in $(A,<,P)$ under $s$ adjusted by $\{(y_1,u_1),\dots,(y_p,u_p)\})$.

This is equivalent to

$C$ holds $\leq x$ in $(A,<,P)$ under $s$. □

We now restate Lemma 5.50 for the new and more usable complexity classes.

LEMMA 5.52. *Let $r, m > 1$. There exists $(A,<,E,P)$ satisfying the following conditions*:

(i) *$(A,<)$ is a linear ordering with a least element in which every element has an immediate successor*;
(ii) *$E = \{e_1 < e_2 < \cdots\}$ is an unbounded set of limit points in $A$*;
(iii) *$P \subseteq A^{8r+m+13}$*;
(iv) *Let $Q \in \mathrm{PBF}(r, 8r+m+13)$ have at most the free variables $x_1,\dots,x_m$. Then there exists $u < e_1$ such that for all $i > 0$ and $x_1,\dots,x_m \leq e_i$,*

$$Q \text{ if and only if } P(e_i,\dots,e_{i+8r+11},u,x_1,\dots,x_m)$$

*holds in $(A,<,P)$*;
(v) *Let $Q \in \mathrm{PBF}(r, 8r+m+13)$ have one free variable and any number of parameters. If there exists $x$ such that $Q(x)$ holds in $(A,<,P)$, then there is a least $x$ such that $Q(x)$ holds in $(A,<,P)$*;
(vi) *$E$ is a strong set of $m$-ary indiscernibles for $\mathrm{PBF}(r, 8r+m+13)$. That is, let $Q(x_1,\dots,x_m) \in \mathrm{PBF}(r, 8r+m+13)$, where all free variables are shown. Let $0 \leq j \leq m$ and $x, y \in E_j$ be of the same order type, and $|z| < \min(x,y)$, $z \in A^{m-j}$. Then $Q(x,z)$ if and only if $Q(y,z)$ holds in $(A,<,P)$.*

*Proof.* Let $r, m > 1$. Let $(A, <, E, P)$ be according to Lemma 5.50. Let $Q$ be as given in (iv). By Lemma 5.50 (iv), let $u < e_1$ be such that for all $i > 0$,

$$Q^* \text{ if and only if } P(e_i, \dots, e_{i+8r+11}, u, x_1, \dots, x_m)$$

holds $\leq e_i$, where $Q^*$ is as defined just before Lemma 5.51. Then (iv) follows by Lemma 5.51.

Let $Q$ be as given in (v). Choose $j$ such that all parameters are $< e_j$ and there exists $x < e_j$ such that $Q(x)$. Then by Lemma 5.51, there exists $x < e_j$ such that $Q^*(x)$ holds $\leq e_j$ with these parameters. By Lemma 5.50 (v), there is a least $x$ such that $Q^*(x)$ holds $\leq e_j$ with these parameters. By Lemma 5.51, this least $x$ is also the least $x$ such that $Q(x)$ holds with these parameters.

Let $Q$, $j$, $x$, $y$, $z$ be as given in (vi). Let $|x, y| < e_q$. Then by Lemma 5.50 (vi), $Q^*(x, z)$ holds $\leq e_q$ if and only if $Q^*(y, z)$ holds $\leq e_q$. By Lemma 5.51, $Q(x, z)$ if and only if $Q(y, z)$. □

We now fix $r$, $m$, and $(A, <, E, P)$ according to Lemma 5.52. We will henceforth assume that $r, m > 7$. Recall that $P \subseteq A^{8r+m+13}$.

We now develop a pairing mechanism within $(A, <, E, P)$.

We define the linear ordering $<^*$ on $A^2$ by $(a, b) <^* (c, d)$ if and only if

(i) $|a, b| < |c, d|$; or

(ii) $|a, b| = |c, d|$ and $(a, b)$ is lexicographically earlier than $(c, d)$.

An initial pairing function (ipf) is a partial binary function $f$ from $A^2$ into $A$ satisfying the following conditions:

(i) $\operatorname{dom}(f)$ is of the form $\{(a, b) \colon (a, b) <^* (u, v)\}$, for some $u, v \in A$;

(ii) $\operatorname{rng}(f)$ is of the form $\{c \colon c < w\}$ for some $w \in A$;

(iii) $f$ is order-preserving in the sense that $(a, b) <^* (a', b')$ if and only if $f(a, b) < f(a', b')$.

Note that $u$, $v$, $w$ above are uniquely determined from $f$. Accordingly, we say that $f$ is an ipf of type $(u, v, w)$.

Note that if $f$ is an ipf of type $(u, v, w)$, then the extension $f^+ = f \cup \{(u, v, w)\}$ is an ipf of type $(u', v', w')$, where $(u', v')$ is the successor of $(u, v)$ in $<^*$ and $w'$ is the successor of $w$. That is, $f^+$ is the successor of the ipf $f$.

Lemma 5.53. *Let $f$, $g$ be* ipf*'s and $x \in \operatorname{dom}(f) \cap \operatorname{dom}(g)$. Suppose that for all $y <^* x$, $f(y) = g(y)$. Then $f(x) = g(x)$.*

*Proof.* Let $f$, $g$, $x$ be as given. Let $f$ be of type $(u, v, w)$ and $g$ be of type $(u', v', w')$. Then $f(x), g(x) < \min(w, w')$. In particular let $f(x) = g(z)$.

If $z <^* x$, then by order preservation, $f(z) < f(x) = g(z)$, contradicting the hypotheses.

Now suppose $x <^* z$. Then $g(x) < g(z) = f(x)$. Let $y$ be such that $f(y) = g(x)$. Since $f(y) < f(x)$ we have $y <^* x$. Hence $f(y) = g(y) = g(x)$, which is a contradiction.

We have thus shown that $z = x$, and so $f(x) = g(x)$, as required. □

Let $x \in A^{8r+m+10}$, and $i > 0$. We write $T(i,x)$ for the 3-ary relation on $A$ defined by $T(i,x) = \{(a,b,c)\colon a,b,c < e_i \,\&\, P(x_1,\dots,x_{8r+m+10},a,b,c)\}$. $T(i,x)$ is only defined for $x \in A^{8r+m+10}$.

LEMMA 5.54. *Every bounded subset of $A^3$ that is definable by a formula in* PBF$(r, 8r+m+13)$ *with $m-3$ parameters and three free variables is of the form $T(i,x)$, for some $i$, $x$. Furthermore, if the set is contained in $\{a\colon a < e_i\}^3$, then it is of the form $T(i,x)$, for some $x$.*

*Proof.* Let $Q = Q(b_1,\dots,b_{m-3},x_{m-2},x_{m-1},x_m)$ define the set

$$X \subseteq \{a\colon a < e_i\}^3\,,$$

where the $b$'s are the parameters and the $x$'s are the free variables, and where $Q \in \text{PBF}(r, 8r+m+13)$. Let $b_1,\dots,b_{m-3}$, $e_i < e_j$.

By Lemma 5.52 (iv), let $u < e_1$ be such that for all $x_1,\dots,x_m < e_j$,

$$Q \text{ if and only if } P(e_j,\dots,e_{j+8r+11},u,x_1,\dots,x_m).$$

Then for all $x_{m-2},x_{m-1},x_m < e_i$,

$$Q(b_1,\dots,b_{m-3},x_{m-2},x_{m-1},x_m) \text{ holds if and only if}$$
$$P(e_j,\dots,e_{j+8r+11},u,b_1,\dots,b_{m-3},x_{m-2},x_{m-1},x_m) \text{ if and only if}$$
$$T(i,y)(x_{m-2},x_{m-1},x_m),$$

where $y = (e_j,\dots,e_{j+8r+11},u,b_1,\dots,b_{m-3})$, as required. □

We are particularly interested in the condition that $T(i,x)$ is (the graph of) an ipf.

LEMMA 5.55. *Let $i > 0$. Any two* ipf*'s of the form $T(i,x)$ are comparable in the sense that one is included in the other.*

*Proof.* Suppose $T(i,x)$ and $T(i',x')$ are ipf's that disagree somewhere. By applying the least element principle from Lemma 5.52 (v) three times, we see that there is a lexicographically least triple $(a,b,c)$ such that

$$a = |b,c| < e_i, \text{ and } T(i,x)(b,c) \neq T(i',x')(b,c)\,.$$

Then $(b,c)$ is $<^*$-least such that $T(i,x)(b,c) \neq T(i',x')(b,c)$. But this contradicts Lemma 5.53. (We have used $r > 7$.) □

LEMMA 5.56. *Let $i > 0$ and $T(i,x)$ be an* ipf.
(i) *For all $(a,b) \in \text{dom}\,(T(i,x))$, $T(i,x)(a,b) \geq a,b$;*
(ii) *The type $(u,v,w)$ of $T(i,x)$ obeys $u,v \leq w$.*

*Proof.* To establish (i), it suffices to show that for all $(0,a) \in \text{dom}\,(T(i,x))$, $T(i,x)(0,a) \geq a$. By the least element principle in Lemma 5.52 (v), let $a$ be least such that $T(i,x)(0,a) < a$. Let $T(i,x)(0,a) = b$. Then $T(i,x)(0,b) < T(i,x)(0,a) = b$. This contradicts the minimality of $a$.

Let the type of $T(i,x)$ be $(u,v,w)$. Suppose $w < u$ or $w < v$. Then $(w,w) <^* (u,v)$, and so $(w,w)$ is in the domain of $T(i,x)$. From the type of $T(i,x)$ we have $T(i,x)(w,w) < w$. But this contradicts (i). □

Lemma 5.57. *There exists* $B \in \mathrm{PBF}(2, 8r+m+13)$ *with* $8r+m+14$ *free variables, such that the following holds. Let* $i > 0$, $x \in A^{8r+m+10}$, *and* $|x|, u, v, w < e_i$. *Then*

$$B(e_i, x, u, v, w) \text{ if and only if}$$
$$T(i,x) \text{ is an ipf of type } (u,v,w).$$

*Proof.* Recall conditions (i)–(iii) in the definition of ipf. $T(i,x)$ defines a function, and conditions (i) and (ii) can be put in $\mathrm{PBF}(2, 8r+m+13)$. Condition (iii) can be put in complexity class $\mathrm{PBF}(1, 8r+m+13)$. □

For $0 < i < j$, let $T^*(i,j)$ be the union of all ipf's of the form $T(i,x)$, where $|x| < e_j$.

Lemma 5.58. (i) *There exists* $B \in \mathrm{PBF}(3, 8r+m+13)$ *with five free variables such that for all* $0 < i < j$ *and* $a,b,c < e_i$,

$$B(e_i, e_j, a, b, c) \text{ if and only if}$$
$$T^*(i,j)(a,b) = c;$$

(ii) $T^*(i,j)$ *and* $T^*(i,j)^+$ *are* ipf'*s*;

(iii) *There exists* $C \in \mathrm{PBF}(5, 8r+m+11)$ *with five free variables such that for all* $0 < i < j$ *and* $a,b,c < e_i$,

$$C(e_i, e_j, a, b, c) \text{ if and only if}$$
$$T^*(i,j)^+(a,b) = c.$$

*Proof.* By Lemma 5.55, $T^*(i,j)$ is an order preserving function from an initial segment of $<^*$ onto an initial segment of $<$. We don't know yet that these initial segments are determined, respectively, by elements of $A^2$ and $A$. This is what is required to establish (ii).

Now (i) is easily verified from Lemma 5.57. The crucial point here is that a large block of existential quantifiers is used, which are treated as a single quantifier in terms of our measure of complexity.

To establish (ii), note that by the least element principle in Lemma 5.52 (v), there is a least element $w$ outside the range of $T^*(i,j)$. Now observe that $\mathrm{dom}\,(T^*(i,j)) \subseteq \{x\colon x < e_i\}$. So by three more uses of the least element principle, we can find a lexicographically least $(|u,v|, u, v)$ such that $(u,v)$ is outside the domain of $T^*(i,j)$. This $(u,v)$ is the $<^*$-least $(u,v)$ outside dom $(T^*(i,j))$. Hence $T^*(i,j)$ is an ipf of type $(u,v,w)$. Therefore $T^*(i,j)^+$ is an ipf.

For (iii), note that the type $(u,v,w)$ of the ipf $T^*(i,j)$ is the unique solution to the following:

(i) $|u,v,w| \leq e_i$;
(ii) For all $(u',v') <^* (u,v)$, $T^*(i,j)(u',v')$ is defined;
(iii) $T^*(i,j)(u,v)$ is not defined;
(iv) For all $w' < w$, $w' \in \operatorname{rng}(T^*(i,j))$;
(v) $w \notin \operatorname{rng}(T^*(i,j))$.

Thus we obtain (iii) from (i) as required, since $T^*(i,j)^+ = T^*(i,j) \cup \{(u,v,w)\}$. □

Lemma 5.59. *Let* $0 < i < j$.
(i) $T^*(i,j) = T^*(i,i+1)$;
(ii) $(\exists x)(T^*(i,i+1)^+ = T(i+1,x))$;
(iii) $T^*(i,i+1)^+ \subseteq T^*(i+1,i+2)$.

*Proof.* To establish (i), let $B$ be according to Lemma 5.58 (i). Now by Lemma 5.52 (vi), for all $a,b,c < e_i$,

$$B(e_i,e_j,a,b,c) \text{ if and only if}$$
$$B(e_i,e_{i+1},a,b,c).$$

Hence for all $a,b,c < e_i$, $T^*(i,j)(a,b) = c$ if and only if $T^*(i,i{+}1)(a,b) = c$, as required.

For (ii), first observe that by Lemma 5.58 (iii), $T^*(i,i+1)^+ \subseteq \{a\colon a < e_{i+1})^3\}$ is PBF$(5,8r+m+13)$-definable with 2 parameters. Hence (ii) follows by Lemma 5.54.

For (iii), fix $k > i+1$ so that there exists $|x| < e_k$ with $T^*(i,i+1)^+ = T(i+1,x)$. This can be regarded as a statement in PBF$(6,8r+m+13)$ about $e_i$, $e_{i+1}$, $e_k$. Hence by the indiscernibility in Lemma 5.52 (vi), we can replace $e_k$ by $e_{i+2}$.

We have thus shown that $T^*(i,i+1)^+$ is of the form $T(i+1,x)$, for some $x < e_{i+2}$. Therefore, (iii) follows immediately by the definition of $T^*(i+1, i+2)$. □

Lemma 5.60. *Let* $i > 0$. *Then* $T^*(i,i+1)$ *is an* ipf *of type* $(0,e_i,e_i)$.

*Proof.* Let $T^*(i,i+1)$ be of type $(u,v,w)$. Then $u,v,w \leq e_i$. Also $T^*(i, i+1)^+(u,v) = w$. By Lemma 5.56 (ii), $u,v \leq w \leq e_i$.

We first establish that $w = e_i$. Suppose $w < e_i$. Then $u,v,w < e_i$. Now clearly $(u,v,w) \in T^*(i,i+1)^+ \backslash T^*(i,i+1)$. Hence by Lemma 5.59 (iii), $(u,v,w) \in T^*(i+1,i+2) \backslash T^*(i,i+1)$. Let $B$ be as given in Lemma 5.58 (i). Then $B(e_{i+1}, e_{i+2}, u, v, w)$ and $\neg B(e_i, e_{i+1}, u, v, w)$. This contradicts the indiscernibility in Lemma 5.52 (vi). Hence $w = e_i$.

Now suppose $u, v < e_i$. As argued previously, $(u, v, e_i) \in T^*(i+1, i+2)$. By Lemma 5.58 (i), $B(e_{i+1}, e_{i+2}, u, v, e_i)$. So by the indiscernibility in Lemma 5.52 (vi), we have $B(e_{i+2}, e_{i+3}, u, v, e_{i+1})$. Hence by Lemma 5.58 (i), $(u, v, e_{i+1}) \in T^*(i+2, i+3)$. But $(u, v, e_i) \in T^*(i+1, i+2) \subseteq T^*(i+2, i+3)$, which contradicts $(u, v, e_{i+1}) \in T^*(i+2, i+3)$.

Hence $u = e_i$ or $v = e_i$. This implies that $\operatorname{dom}(T^*(i, i+1)) = \{a : a < e_i\}^2$. Hence $(u, v)$ is the least upper bound of $\{a : a < e_i\}^2$, which is $(0, e_i)$, as required. □

We now let $T^\#$ be the union over $i > 0$ of all $T^*(i, i+1)$.

Lemma 5.61. (i) *For all $i > 0$, $T^\#$ is a bijection from $\{a : a < e_i\}^2$ onto $\{a : a < e_i\}$;*

(ii) *There is a formula $B \in \mathrm{PBF}(3, 8r + m + 13)$ with five free variables such that for all $0 < i < j$ and $x, y, z < e_i$,*

$$T^\#(x, y) = z \text{ if and only if } B(e_i, e_j, x, y, z)\,;$$

(iii) *$T^\#$ is a bijection from $A^2$ onto $A$.*

*Proof.* By Lemmas 5.59 (iii) and 5.60, the $T^*(i, i+1)$ are comparable ifp's which are themselves bijections from $\{a : a < e_i\}^2$ onto $\{a : a < e_i\}$. Hence (i) follows.

Let $0 < i < j$. By Lemmas 5.59 (i) and 5.60, $T^\# | \{x : x < e_i\}^2 = T^*(i, i+1) = T^*(i, j)$. So (ii) follows from Lemma 5.58 (i).

Obviously (iii) follows from (i) since $E$ is unbounded in $A$. □

We extend the binary function $T^\#$ to a function that applies to any number (at least one) of arguments from $A$. Define $T^\#(x) = x$. For $k > 2$ we define $T^\#(x_1, \ldots, x_{k+1}) = T^\#(T^\#(x_1, \ldots, x_k), x_{k+1})$.

5.6. *Existence of linearly ordered binary relations with pairing function, limit points, strong indiscernibles, comprehension, and the least element principle.*

Lemma 5.62. (i) *For all $i, k > 0$, $T^\#$ is a bijection from $\{a : a < e_i\}^k$ onto $\{a : a < e_i\}$;*

(ii) *There is a formula $B \in \mathrm{PBF}(4, 8r + m + 13)$ with $k + 3$ free variables such that for all $k > 0$, $0 < i < j$, and $x_1, \ldots, x_k$, $y < e_i$,*

$$T^\#(x_1, \ldots, x_k) = y \text{ if and only if } B(e_i, e_j, x_1, \ldots, x_k, y)$$

*holds in $(A, <, P)$;*

(iii) *For all $k > 0$, $T^\#$ is a bijection from $A^k$ onto $A$.*

*Proof.* For (i), the case $k = 2$ is from Lemma 5.61. The case $k = 1$ is obvious. Assume true for $k, i > 0$. It suffices to show that it is true for $k+1$ and $i$. But this induction step is clear from the definition of $T^{\#}(x_1, \ldots, x_{k+1})$.

Clearly (ii) follows from Lemma 5.61 (ii).

Obviously (iii) follows from (i) since $E$ is unbounded in $A$. □

We define the important binary relation $S \subseteq A^2$ as follows. $S(x, y)$ if and only if there exists (necessarily unique) $z_1, \ldots, z_{8r+m+12}$ such that $P(z_1, \ldots, z_{8r+m+12}, y)$ and $T^{\#}(z_1, \ldots, z_{8r+m+12}) = x$.

LEMMA 5.63. *There is a formula $B \in \mathrm{PBF}(5, 8r+m+13)$ with four free variables such that the following holds in $(A, <, P)$. For all $0 < i < j$ and $x, y < e_i$,*

$$S(x, y) \text{ if and only if } B(e_i, e_j, x, y)\,.$$

*Proof.* Note that for all $x, y < e_i$, $S(x, y)$ if and only if there exists $z_1, \ldots, z_{8r+m+12} < e_i$ such that $P(z_1, \ldots, z_{8r+m+12}, y)$ and $T^{\#}(z_1, \ldots, z_{8r+m+12}) = x$. Hence the claim follows from Lemma 5.62 (ii). □

We say that $B \subseteq A^k$ is bounded if and only if there exists $x \in A$ such that for all $y \in B$, $|y| < x$.

LEMMA 5.64. *Let $n > 0$ and $X$ be a bounded subset of $A^n$ that is definable over $(A, <, P)$ by a formula in $\mathrm{PBF}(r-1, 8r+m+13)$ with any parameters and $n$ free variables. Then there exist $x, a \in A$ such that $X = \{y \in A^n \colon |y| < a \,\&\, S(x, T^{\#}(y))\}$.*

*Proof.* Let $i, t > 0$, $w \in A^t$, $|w| < e_i$, and $X = \{y \in A^n \colon |y| < e_i \,\&\, B(w, y)\}$, where $B \in \mathrm{PBF}(r-1, 8r+m+13)$ has $t+n$ free variables. The $t$ coordinates of $w$ serve as parameters.

Let $Y = T^{\#}[X]$ and $v = T^{\#}(w)$. Using Lemma 5.62 (ii), write $Y = \{z < e_i \colon C(e_i, e_{i+1}, v, z)\}$, where $C \in \mathrm{PBF}(r, 8r+m+13)$ has four free variables.

By Lemma 5.52 (iv), let $u < e_1$ be such that for all $z < e_i$,

$$C(e_i, e_{i+1}, v, z) \text{ if and only if}$$
$$P(e_{i+2}, \ldots, e_{i+8r+13}, u, e_i, e_{i+1}, v^*m-3, z).$$

Then $Y = \{z < e_i \colon S(x, z)\}$, where

$$x = T^{\#}(e_{i+2}, \ldots, e_{i+8r+13}, u, e_i, e_{i+1}, v^*m-3)\,.$$

Hence $X = \{y \in A^n \colon |y| < e_i \,\&\, S(x, T^{\#}(y))\}$ as required. □

Let $n > 0$ and $X \subseteq A^n$. We say that $X$ is boundedly definable over $(A, <, P)$ if and only if $X$ is bounded and $X$ can be defined by a formula in $\mathrm{PBF}(\infty, 8r+m+13)$ over $(A, <, P)$ with any parameters and $n$ free variables.

LEMMA 5.65. *Let $n>0$ and $X \subseteq A^n$ be boundedly definable over $(A,<,P)$. Then there exist $x, a \in A$ such that $X = \{y \in A^n\colon |y| < a \,\&\, S(x, T^{\#}(y))\}$.*

*Proof.* We show by induction on $i \geq 0$ that this is true for complexity class $\mathrm{PBF}(i, 8r+m+13)$, with parameters allowed. By Lemma 5.64 this is true for $\mathrm{PBF}(r-1, 8r+m+13)$. Suppose this is true for $\mathrm{BF}(i, 8r+m+13)$, $i \geq r-1$. We first show that this is true for $\mathrm{PBF}(i, 8r+m+13)$.

Let $X$ be a Boolean combination of sets of the form $(y \in A^n\colon U(y)\}$, where $X$ is bounded and $U$ is in $\mathrm{BF}(i, 8r+m+13)$ with parameters allowed. Then $X$ is also a Boolean combination of bounded sets of this form. By the induction hypothesis, $X$ is a Boolean combination of sets of the form $\{y \in A^n\colon |y| < a \,\&\, S(x, T^{\#}(y))\}$. By Lemmas 5.62 and 5.63, $X$ is a bounded set that is definable by a formula in $\mathrm{PBF}(6, 8r+m+13)$ with parameters. By Lemma 5.64, $X$ can be defined in the appropriate form.

We now show that this is true for $\mathrm{PB}(i+1, 8r+m+13)$. Accordingly, let $X = \{y \in A^n\colon |y| < e_j \,\&\, (\exists\, z)(z_1 \leq w_1 \,\&\, \ldots \,\&\, z_p \leq w_p \,\&\, B(z,y))\}$, where $B \in \mathrm{PBF}(i, 8r+m+13)$, each $w_i$ is either a parameter or some $y_j$, parameters are allowed in $B$, and all parameters used throughout are $< e_j$. By the induction hypothesis, let $x, a \in A$ be such that for all $z \in A^p$, $y \in A^n$, if $|z,y| < e_j$, then

$$B(z,y) \text{ if and only if } |z,y| < a \,\&\, S(x, T^{\#}(z,y))\,.$$

Hence for all $y \in A^n$, $y \in X$ if and only if

$$|y| < e_j \,\&\, (\exists\, z)(z_1 \leq w_1 \,\&\, \ldots \,\&\, z_p \leq w_p \,\&\, |z,y| < a \,\&\, S(x, T^{\#}(z,y)))\,.$$

By Lemmas 5.62 (ii) and 5.63, $X$ is definable by a formula in $\mathrm{PBF}(6, 8r+m+13)$ with parameters. The induction step, therefore, follows by Lemma 5.64, and we have established the claim. □

LEMMA 5.66. *Every nonempty subset of $A$ that is boundedly definable over $(A,<,P)$ has a least element.*

*Proof.* This follows immediately from Lemmas 5.52 (v), 5.62, 5.63, and 5.65. □

We now consider bounded formulas in the system $(A,<,S,T^{\#})$. Here the language is based on the binary relations $<$, $S$, and the binary function $T^{\#}$.

Specifically, we define $\mathrm{BF}^{\#}(0)$ as the atomic formulas in the language of $(A,<,S,T^{\#})$. (This will involve nested terms, as is usual in first order predicate calculus.)

$\mathrm{BF}^{\#}(i+1)$ consists of all formulas of the form

$$(\exists\, y_1, \ldots, y_p)(y_1 \leq t_1 \,\&\, \ldots \,\&\, y_p \leq t_p \,\&\, B)\,,$$

where $p \geq 0$, $B \in \mathrm{PBF}^{\#}(i)$, $y_1, \ldots, y_p$ are distinct variables, and $t_1, \ldots, t_p$ are terms in which no $y$ appears.

$\mathrm{PBF}^{\#}(i)$ consists of all propositional combinations of elements of $\mathrm{BF}^{\#}(i)$.

The bounded formulas for $(A, <, S, T^{\#})$ are the formulas in $\mathrm{PBF}^{\#}(\infty)$.

A subset of $A^n$ is boundedly definable in $(A, <, S, T^{\#})$ if and only if it is bounded and definable over $(A, <, S, T^{\#})$ by a bounded formula with any parameters and $n$ free variables.

LEMMA 5.67. *For all $n > 0$ and $X \subseteq A^n$, $X$ is boundedly definable over $(A, <, S, T^{\#})$ if and only if $X$ is boundedly definable over $(A, <, P)$.*

*Proof.* We are not claiming that every subset of $A^n$ that is definable over $(A, <, S, T^{\#})$ by a bounded formula is definable over $(A, <, P)$ by a bounded formula, or vice versa. We need to know that the subset is bounded.

Suppose that $X$ is boundedly definable over $(A, <, P)$. Then the representation of $X$ given by Lemma 5.65 shows that $X$ is boundedly definable over $(A, <, S, T^{\#})$.

Now suppose that $X$ is boundedly definable over $(A, <, S, T^{\#})$. Let $X \subseteq (a\colon a < e_i)^n$. Then using Lemmas 5.62 and 5.63, we can replace every atomic subformula appearing in the definition of $X$ by a bounded formula with parameters $e_i$, $e_{i+1}$, so that the definition is still valid for all arguments $< e_i$. This also relies on the fact that $T^{\#}$ maps each $\{a\colon a < e_i\}^m$ into $\{a\colon a < e_i\}$. □

We now need some more precise information. We need to refer back to $(A, <, P)$.

LEMMA 5.68. *Let $n > 0$, $j \geq 0$, and $Q(x_1, \dots, x_n) \in \mathrm{PBF}^{\#}(j)$ with only the free variables shown. Then there exists $Q^*(x_1, \dots, x_{n+2}) \in \mathrm{PBF}(j + 6, 8r + m + 13)$ with only the free variables shown such that the following holds. For all $i > 0$ and $x_1, \dots, x_n < e_i$,*

$$Q(x_1, \dots, x_n) \textit{ holds in } (A, <, S, T^{\#}) \textit{ if and only if}$$
$$Q^*(x_1, \dots, x_n, e_i, e_{i+1}) \textit{ holds in } (A, <, P).$$

*Proof.* The proof is by induction on $j$. We rely on the fact that $T^{\#}$ maps each $(a\colon a < e_i)^m$ into $\{a\colon a < e_i\}$.

Let $j = 0$, and $Q(x_1, \dots, x_n)$ be as given. Let $Q'(x_1, \dots, x_n)$ be the obvious predicate calculus equivalent of $Q$ obtained by putting a block of existential quantifiers on $Q$ and replacing all terms in $Q$ by unnested terms. Let $Q''(x_1, \dots, x_{n+1})$ result from $Q'$ by strictly bounding these existential quantifiers to $x_{n+1}$. Finally, let $Q^*(x_1, \dots, x_{n+2})$ result from $Q''$ by using the definitions of $S$, $T^{\#}$ in Lemmas 5.61 and 5.63, where the two parameters appearing there, $e_i$ and $e_j$, are replaced by $x_{n+1}$ and $x_{n+2}$. This results in the required $Q^* \in \mathrm{PBF}(6, 8r + m + 13)$.

Suppose this is true for $\mathrm{PBF}^{\#}(j)$. We handle $\mathrm{BF}^{\#}(j + 1)$ by changing the bound variables so that they are not any $x_1, \dots, x_{n+2}$, adding additional existential quantifiers to unravel the bounding terms, bounding all the relevant existential quantifiers to $x_{n+1}$, and using the definitions from Lemmas 5.61 and 5.63 as above. The resulting $Q^*$ is in $\mathrm{PB}(j+6, 8r+m+13)$. Finally $\mathrm{PBF}^{\#}(j+1)$ is handled by passing through the propositional operations. □

LEMMA 5.69. *Let $Q(x_1, \ldots, x_{m-2}) \in \mathrm{PBF}^{\#}(r-6)$ where all free variables are shown. Let $0 \leq j \leq m-2$ and $x, y \in E^j$ be of the same order type, and $|z| < \min(x, y)$, $z \in A^{m-2-j}$. Then $Q(x, z)$ if and only if $Q(y, z)$ holds in $(A, <, S, T^{\#})$.*

*Proof.* Let $Q$, $j$, $x$, $y$, $z$ be as given. Let $|x, y| < e_i$ and $Q^*(x_1, \ldots, x_m)$ be according to Lemma 5.68. Then $Q^* \in \mathrm{PBF}(r, 8r+m+13)$. By Lemma 5.52 (vi), we have

$$Q^*(x, z, e_i, e_{i+1}) \text{ if and only if } Q^*(y, z, e_i, e_{i+1})$$

holds in $(A, <, P)$. By Lemma 5.68 we see that $Q(x, z)$ if and only if $Q(y, z)$ holds in $(A, <, S, T^{\#})$ as required. □

We now transfer entirely over to $(A, <, E, S, T^{\#})$ and never consider $P$ again. From now on, we say that $X \subseteq A^n$ is boundedly definable if and only if $X$ is bounded and can be defined by a formula in $\mathrm{PBF}^{\#}(\infty)$ with any parameters allowed.

LEMMA 5.70. *The system $(A, <, E, S, T^{\#})$ obeys the following*:

(i) *$(A, <)$ is a linear ordering with a least element in which every element has an immediate successor*;
(ii) *$E = \{e_1 < e_2 < \cdots\}$ is an unbounded set of limit points in $A$*;
(iii) *$S \subseteq A^2$*;
(iv) *For all $i > 0$, $T^{\#}$ is a bijection from $\{x : x < e_i\}^2$ onto $\{x : x < e_i\}$*;
(v) *Let $n > 0$ and $X$ be a boundedly definable subset of $A^n$. Then there exists $x, a \in A$ such that $X = \{y \in A^n : |y| < a \,\&\, S(x, T^{\#}(y))\}$*;
(vi) *Every nonempty boundedly definable subset of $A$ has a least element*;
(vii) *Let $Q(x_1, \ldots, x_{m-2}) \in \mathrm{PBF}\#(r-6)$, where all free variables are shown. Let $0 \leq j \leq m-2$ and $x, y \in E^j$ be of the same order type, and $|z| < \min(x, y)$, $z \in A^{m-2-j}$. Then $Q(x, z)$ if and only if $Q(y, z)$.*

*Proof.* This is immediate from Lemmas 5.65, 5.66, 5.67, and 5.69. □

LEMMA 5.71. *There exists a system $(A, <, E, S, G)$ such that the following holds*:

(i) *$\{A, <\}$ is a linear ordering with a least element in which every element has an immediate successor*;
(ii) *$E = \{e_1 < e_2 < \ldots\}$ is an unbounded set of limit points in $A$*;
(iii) *$S \subseteq A^2$*;

(iv) *For all* $i > 0$, $G$ *is a bijection from* $\{x\colon x < e_i\}^2$ *onto* $\{x\colon x < e_i\}$;
(v) *Let* $n > 0$ *and* $X$ *be a boundedly definable subset of* $A^n$. *Then there exists* $x, a \in A$ *such that* $X = \{y \in A^n\colon |y| < a \,\&\, S(x, G(y))\}$;
(vi) *Every nonempty boundedly definable subset of* $A$ *has a least element*;
(vii) *Let* $n > 0$ *and* $Q(x_1, \ldots, x_n)$ *be a bounded formula, where all free variables are shown. Let* $0 \leq j \leq n$ *and* $x, y \in E^j$ *be of the same order type, and* $|z| < \min(x, y)$, $z \in A^{n-j}$. *Then* $Q(x, z)$ *if and only if* $Q(y, z)$.

(In (v), $G$ is extended to $A^n$, $n > 0$, by: for $n = 1$, $G(x) = x$; for $n > 2$, $G(x_1, \ldots, x_{k+1}) = G(G(x_1, \ldots, x_k), x_{k+1})$. Also boundedly definable refers to the system $(A, <, S, G)$, where $G$ is used instead of $T^\#$.)

*Proof.* We can try to obtain this result without the unboundedness of $E$ in clause (ii) from Lemma 5.70 by applying the compactness theorem for predicate calculus. We would like to cut down to the initial segment of $A$ determined by $E$. But this may ruin clause (v).

We modify this argument as follows. We treat $E$ as a monadic predicate symbol (instead of an infinite sequence of constants) and replace clause (ii) with "$E$ is unbounded above." We then obtain a system $(A, <, E, S, G)$ as desired, except that clause (ii) reads "$E$ is an unbounded set of limit points in $A$." We then make $E$ of order type $\omega$ by using any cofinal subset of $E$ of order type $\omega$. □

We now fix a system $(A, <, E, S, G)$ according to Lemma 5.71. Henceforth, definability and boundedly definable will refer to $(A, <, S, G)$.

LEMMA 5.72. *Let* $n \geq 0$, $i > 0$, *and* $B(x_1, \ldots, x_n)$ *be a prenex formula with free variables shown. Let* $B/e_i(x_1, \ldots, x_n)$ *be the result of bounding all of the quanifiers in* $B$ *to* $\{a\colon a < e_i\}$. *Then for all* $x_1, \ldots, x_n < e_i$, $B(x_1, \ldots, x_n)$ *if and only if* $B/e_i(x_1, \ldots, x_n)$.

*Proof.* By induction on the number $k \geq 0$ of quantifiers in $B$, the basis case being vacuous, we may assume that for all $n \geq 0$, $i > 0$, prenex $B(x_1, \ldots, x_n)$ with $k$ quantifiers, and $x_1, \ldots, x_n < e_i$, $B(x_1, \ldots, x_n)$ if and only if $B/e_i(x_1, \ldots, x_n)$. Let $B(x_1, \ldots, x_n)$ be prenex with $k + 1$ quantifiers and $i > 0$.

Write $B(x_1, \ldots, x_n) = (Qy)(C(x_1, \ldots, x_n, y))$ where $C$ has $k$ quantifiers. Let $x_1, \ldots, x_n < e_i$. Choose $i \leq j$ such that

$$(Qy < e_j)(C(x_1, \ldots, x_n, y)) \text{ if and only if } (Qy)(C(x_1, \ldots, x_n, y)).$$

By the induction hypothesis, for all $y < e_j$,

$$C(x_1, \ldots, x_n, y) \text{ if and only if } C/e_j(x_1, \ldots, x_n, y).$$

Hence

$$(Qy)(C(x_1,\dots,x_n,y)) \text{ if and only if}$$
$$(Qy < e_j)(C(x_1,\dots,x_n,y)) \text{ if and only if}$$
$$(Qy < e_j)(C/e_j(x_1,\dots,x_n,y)) \text{ if and only if}$$
$$B/e_j(x_1,\dots,x_n).$$

Now by the indiscernibility from Lemma 5.71 (vii),

$$B/e_j(x_1,\dots,x_n) \text{ if and only if}$$
$$B/e_i(x_1,\dots,x_n).$$

Hence

$$B(x_1,\dots,x_n) \text{ if and only if}$$
$$B/e_i(x_1,\dots,x_n),$$

as required. □

We say that a subset of $A^n$ is definable if and only if it can be defined by a formula (not necessarily bounded) with any parameters from $A$.

Lemma 5.73. *For all $n > 0$, every bounded definable subset of $A^n$ is boundedly definable.*

*Proof.* This is by Lemma 5.72, where $e_i$ is chosen to be greater than the parameters used and the bound on the subset of $A^n$. □

Lemma 5.74. *There exists a system $(A,<,E,S,G)$ such that the following holds*:

(i) *$(A,<)$ is a linear ordering with a least element in which every element has an immediate successor*;
(ii) *$E = \{e_1 < e_2 < \cdots\}$ is an unbounded set of limit points in $A$*;
(iii) *$S \subseteq A^2$*;
(iv) *For all $i > 0$, $G$ is a bijection from $\{a\colon a < e_i\}^2$ onto $\{a\colon a < e_i\}$*;
(v) *Let $n > 0$ and $X$ be a bounded definable subset of $A^n$. Then there exist $x, a \in A$ such that $X = \{y \in A^n\colon |y| < a \,\&\, S(x, G(y))\}$*;
(vi) *Every nonempty definable subset of $A$ has a least element*;
(vii) *Let $n > 0$ and $Q(x_1,\dots,x_n)$ be a formula, where all free variables are shown. Let $0 \le j \le n$ and $x, y \in E^j$ be of the same order type, and $|z| < \min(x,y)$, $z \in A^{n-j}$. Then $Q(x,z)$ if and only if $Q(y,z)$.*

*Proof.* This is immediate from Lemmas 5.71 and 5.73 except for (vii). Let $n$, $Q$, $j$, $x$, $y$, $z$ be as given. We can assume that $Q$ is prenex. Let $|x,y| < e_i$. By Lemma 5.71 (vii),

$$Q/e_i(x,z) \text{ if and only if } Q/e_i(y,z).$$

By Lemma 5.72,

$$Q(x,z) \text{ if and only if } Q/e_i(x,z), \text{ and}$$
$$Q(y,z) \text{ if and only if } Q/e_i(y,z).$$

Thus $Q(x,z)$ if and only if $Q(y,z)$, as required. □

5.7. *Construction of the cumulative hierarchy.*

LEMMA 5.75. *Let $n, i > 0$ and $X \subseteq \{b{:}\, b < e_i\}^n$ be definable. Then there exists $x$ and $a \leq e_i$ such that $X = \{y \in A^n{:}\, |y| < a \,\&\, S(x, G(y))\}$. Furthermore, if this equation holds for $x, a$, where $a \leq e_i$, then there exists $x < e_{i+1}$ such that this equation holds for $x, a$.*

*Proof.* The first claim is immediate from Lemma 5.74 (v); the bound on $a$ is obtained by the hypothesis on $X$.

Fix $n, i > 0$ and $a \leq e_i$. For all $x \in A$, define $f(x) = \{y \in A^n{:}\, |y| < a \,\&\, S(x, G(y))\}$. We say that $x$ is critical if and only if for all $x' < x$, $f(x') \neq f(x)$.

In order to complete the proof, it suffices to show that every critical $x$ lies strictly below $e_{i+1}$. Suppose this is false. Let $z_1$ be the least critical $x \geq e_{i+1}$. Then $z_1$ is definable from $a$ and $e_{i+1}$. Hence by the indiscernibility given in Lemma 5.74 (vii), $e_{i+1} \leq z_1 < e_{i+2}$.

By indiscernibility, some critical $x$ is $\geq e_{i+2}$. Hence let $z_2$ be the least critical $x \geq e_{i+2}$. Then by indiscernibility, $e_{i+2} \leq z_2 < e_{i+3}$.

By the same argument, let $z_3$ be the least critical $x \geq e_{i+3}$. Hence also $e_{i+3} \leq z_3 < e_{i+4}$.

Now we wish to compare the sets $f(z_1)$, $f(z_2)$, and $f(z_3)$. Since $z_1$, $z_2$, and $z_3$ are critical, these sets differ. Let $y_1$, $y_2$, and $y_3$ be, respectively, the lexicographically least element of $f(z_1)\Delta f(z_2)$, $f(z_2)\Delta f(z_3)$, and $f(z_1)\Delta f(z_3)$. Note that $|y_1|, |y_2|, |y_3| < a \leq e_i$. Hence by indiscernibility, $y_1 = y_2 = y_3$. But this is a contradiction. □

We say that a partial function from $A^n$ into $A$ is definable if and only if its graph is a definable subset of $A^{n+1}$. Of course, any parameters are allowed.

LEMMA 5.76. *Let $n > 0$. Every $n$-ary definable function with bounded domain has bounded range. For all $i > 0$, the field of every $n$-ary function with bounded domain that is definable with parameters $< e_i$ is contained in some $\{a{:}\, a < x\}$, where $x < e_i$.*

*Proof.* Clearly the second claim implies the first claim. Let $i > 0$ and $F$ be as given. Let $i < j$ be such that

$$(\exists\, y < e_j)(\operatorname{dom}(F) \subseteq \{a{:}\, a < y\}^n)\,.$$

By Lemma 5.74 (vii),

$$(\exists\, y < e_i)(\operatorname{dom}(F) \subseteq \{a{:}\, a < y\}^n)\,.$$

Fix such a $y$. By Lemma 5.74 (vii), clearly $\operatorname{rng}(F) \subseteq \{a : a < e_i\}$. Hence

$$(\exists\, z < e_{i+1})(\operatorname{rng}(F) \subseteq \{a : a < z\}).$$

By Lemma 5.74 (vii),

$$(\exists\, z < e_i)(\operatorname{rng}(F) \subseteq \{a : a < z\}).$$

Fix such a $z$. Then $x = \max(y, z)$ is as required. □

Lemma 5.77. *Let $B \subseteq A$ be definable and bounded, and $x \in A$. Then there are unique $y, f$ such that $f$ is a definable order-preserving map from $B$ one-one onto $[x, y)$.*

*Proof.* Fix $B$ as given. Uniqueness is immediate by the least element principle in Lemma 5.74 (vi).

Let $K(z, w, a, b)$ hold if and only if

$$\{u \in A^2 : |u| < a \,\&\, S(b, G(u))\} \text{ is an order-preserving function } f_{a,b} \text{ from } B \cap [0, z) \text{ onto } [x, w).$$

By the least element principle, if $K(z, w, a, b)$ and $K(z, w', a', b')$, then $w = w'$ and $f_{a,b} = f_{a',b'}$. Thus if $K(z, w, a, b)$, then we can write $f_z$ for this unique $f_{a,b}$, which is the unique definable order-preserving function from $B \cap [0, z)$ onto $[x, w)$, since it depends only on $z$.

We prove that for all $z$ there exists $w$, $a$, $b$ such that $K(z, w, a, b)$. Suppose this is false for some $z$. By the least element principle, assume $z$ is the least such that this is false.

Now for all $z' < z$, $f_{z'}$ is the unique definable order-preserving function from $B \cap [0, z')$ onto some $[x, w')$. By the least element principle, for all $z' < z'' < z$, $f_{z'} \subseteq f_{z''}$.

Let $g$ be the union of the $f_{z'}$, $z' < z$. Then $g$ is an order-preserving function from $B \cap [0, z)$ onto an initial segment $I$ of $[x, \infty)$. But also note that using the definition of the $f_{z'}$, $z' < z$, given in terms of the definition of $K$, we see that $g$ is definable. Hence by Lemma 5.76, the range of $g$ is bounded. So by the least element principle, $I$ must be the form $[x, w)$. Therefore, by Lemma 5.74 (v), there exist $a, b$ such that $K(z, w, a, b)$. This is the required contradiction. □

We now develop the cumulative hierarchy within $(A, <, E, S, G)$. An initial cumulative hierarchy (ich) as a system $(x, y, R, \mathrm{rk})$ satisfying the following conditions:

(i) $x, y \in A$;
(ii) $R \subseteq \{a : a < y\}^2$;
(iii) Let $b, c < y$. If $(\forall\, a < y)(R(a, b)$ if and only if $R(a, c))$, then $b = c$;
(iv) $\mathrm{rk}\colon \{a : a < y\}$ onto $\{a : a < x\}$;

(v) Let $a < y$. Then $\mathrm{rk}(a) = \sup\{\mathrm{rk}(b) : R(b,a)\}$, where this is the strict sup;
(vi) Let $b < x$ and $Y$ be a definable subset of $\{z : \mathrm{rk}(z) < b\}$. Then there exists $w < y$ such that $Y = \{z : R(z,w)\}$;
(vii) If $a \le b < y$, then $\mathrm{rk}(a) \le \mathrm{rk}(b)$.

From the point of view of set theory, this intuitively means that
(1) $(x, y, R, \mathrm{rk})$ codes the cumulative hierarchy of height $x$;
(2) the domain of the cumulative hierarchy structure is $\{a : a < y\}$;
(3) rk is the rank function from $\{a : a < y\}$ onto $\{a : a < x\}$;
(4) for all $b < x$, every set of domain elements of rank $< b$ is coded in the cumulative hierarchy;
(5) for all $a, b$ in the domain, if $\mathrm{rk}(a) < \mathrm{rk}(b)$, then $a < b$.

We will be using the following terminology:
(a) $x$ is the height of the ich;
(b) $y$ is the domain point of the ich;
(c) $R$ is the epsilon relation of the ich;
(d) rk is the rank function of the ich.

We say that an ich is definable if and only if its epsilon relation and rank function are definable (with any parameters allowed).

LEMMA 5.78. *Let* $X = (x, y, R, \mathrm{rk})$ *be a definable* ich. *Then there is a definable* ich $X^{\#} = (x^+, y', R', \mathrm{rk}')$ *such that*
(i) $y < y'$ *and* $\mathrm{rk} \subseteq \mathrm{rk}'$;
(ii) $R = R' \cap \{a : a < y\}^2$.

*Proof.* Let $X$ be as given, and let $y < e_i$. We say that $(b, c) \in A^2$ is minimal over $X$ if and only if
(a) $\{a < b : S(c,a)\}$ is different from every $\{a : R(a,d)\}$, $d < y$;
(b) $\{a < b : S(c,a)\} \subseteq [0, y)$;
(c) for all $(b', c') <_{\mathrm{lex}} (b, c)$, $\{a < b' : S(c', a)\} \ne \{a < b : S(c,a)\}$.

In other words, $(b, c)$ codes a subset of $[0, y)$ which is not coded inside $X$, and not coded by any lexicographically earlier $(b', c')$.

Let $T = \{G(b,c) : (b,c)$ is minimal over $X\}$. Clearly $T$ is definable. Also observe that by Lemma 5.75, $T$ is strictly bounded by $e_{i+1}$. By Lemma 5.77, let $f, y'$ be unique such that $f$ is a definable order-preserving function from $T$ onto $[y, y')$.

We define $R'$ by: $R'(z,w)$ if and only if either $(w < y \,\&\, R(z,w))$ or $(\exists\, b,c)(f(G(b,c)) = w \,\&\, z < b \,\&\, S(c,z))$. We extend rk to rk$'$ by setting all $\mathrm{rk}'(f(u)) = x$.

We now define $X^{\#} = (x^+, y', R', \mathrm{rk}')$. $X^{\#}$ obviously obeys clauses (i)–(iii), (vii) in the definition of ich.

We claim that for all $G(b,c) \in T$, $x = \sup\{\mathrm{rk}(z) : R'(z, f(G(b,c)))\}$. Note that $R'(z, f(G(b,c)))$ if and only if $(S(c,z) \,\&\, z < b)$.

Suppose that this is false for some $G(b,c) \in T$. Let $x' < x$ be strictly greater than all $\mathrm{rk}(z)$ for which $S(c,z) \,\&\, z < b$. Then $\{z < b\colon S(c,z)\}$ is a definable subset of $\{z\colon \mathrm{rk}(z) < x'\}$. So by clause (vi) in the definition of ich, we see that $\{z < b\colon S(c,z)\}$ is of the form $\{z\colon R(z,w)\}$, $w < y$. But this contradicts that $G(b,c) \in T$. This establishes the claim.

Hence clause (v) holds for $X^{\#}$. To establish (vi) for $X^{\#}$, let $u < x^{+}$ and $Y$ be a definable subset of $\{z\colon \mathrm{rk}'(z) < u\}$. Then $Y$ is a definable subset of $[0,y)$. By Lemma 5.74 (v), $Y$ is of the form $\{a < b\colon S(c,a)\}$. By the construction of $X^{\#}$ and the least element principle, we see that $Y$ is of the form $\{z\colon R'(z,w)\}$, $w < y'$.

It remains to establish clause (iv) in the definition of ich. It suffices to prove that there is a definable subset of $[0,y)$ that differs from every $\{a{:}R(a,d)\}$, $d < y$. But $\{d < y\colon \neg R(d,d)\}$ is such a definable subset of $[0,y)$. □

Lemma 5.79. *Let* $X = (x,y,R,\mathrm{rk})$ *be a definable* ich *and* $a \leq x$. *Then* $\{z\colon \mathrm{rk}(z) < a\}$ *is of the form* $[0,b)$, $b \leq y$. *Also* $X \mid a = (a,b,R \cap [0,b)^2, \mathrm{rk} \cap [0,b)^2)$ *is a definable* ich.

*Proof.* Let $X$, $a$ be as given. By clause (vii) for $X$, $\{z\colon \mathrm{rk}(z) < a\}$ is an initial segment of $A$. By the least element principle, let $b$ be such that $\{z\colon \mathrm{rk}(z) < a\} = [0,b)$. All the clauses in the definition of ich clearly hold for $(a,b,R \cap [0,b\}^2, \mathrm{rk} \cap [0,b)^2)$ except for (vi). Let $c < a$ and $Y$ be a definable subset of $\{z\colon \mathrm{rk}(z) < c\}$. By (vi) for $X$, let $w < y$, $Y = \{z\colon R(z,w)\}$. Then $\mathrm{rk}(w) \leq c < a$, and hence $w < b$ as required. □

For any definable ich $X$, we let $X^{\#}$ be the definable ich given by the construction in the proof of Lemma 5.78. Notice that $X^{\#}$ depends on $X$ and not on the definition of $X$ (although the existence of $X^{\#}$ depends on the existence of a definition of $X$).

We say that $X = (x,y,R,\mathrm{rk})$ is a preferred ich if and only if

(i) $X$ is a definable ich;

(ii) for all $a^{+} \leq x$, $X \mid a^{+} = (X \mid a)^{\#}$.

To make this more precise, we can use the representation in Lemma 5.74 (v) of bounded definable subsets of $A^n$ to introduce codes for definable ich's as certain quadruples from $A$. Using the proof of Lemma 5.78, we then define a binary relation $B$ such that for all such codes $\alpha$, $\beta$, $B(\alpha,\beta)$ if and only if $\beta$ codes the # of the ich coded by $\alpha$. Then the definition of preferred ich can be given as follows: $X = (x,y,R,\mathrm{rk})$ is a preferred ich if and only if

(i) $X$ is a definable ich;

(ii) for all $a^{+} \leq x$, we have $B(\alpha,\beta)$, where $\alpha$ is a code of the ich $X|a$ and $\beta$ is a code of the ich $X|a^{+}$.

Lemma 5.80. *Let $X$, $X'$ be preferred ich's of heights $x$, $x'$, with rank functions* rk, rk$'$. *For all $a \leq x$, $X \mid a$ is a preferred* ich. *If $x \leq x'$ then $X = X' \mid x$. Any two preferred* ich*'s of the same height are equal.*

*Proof.* Let $X$, $X'$, $x$, $x'$ be as given. Assume $a \leq x$. By Lemma 5.79, $X \mid a$ is a definable ich. Let $b+ < a$. Then $X \mid a \mid b^+ = X \mid b^+ = (X \mid b)^\# = (X \mid a \mid b)^\#$.

For the second claim, assume $x \leq x'$. We claim that for all $b \leq x$, $X \mid b = X' \mid b$. By way of contradiction, let $b$ be least such that $X \mid b \neq X' \mid b$. By clause ii) in the definition of preferred ich, clearly $b$ is not a successor in $A$. Hence $b$ is a limit.

Clearly the heights of $X \mid b$ and $X' \mid b$ are both $b$. Also the domain, epsilon relation, and rank function of $X \mid b$ is the union of the domain, epsilon relation, and rank function of $X \mid c$, $c < b$. The same holds for $X'$. Hence $X \mid b = X' \mid b$, which is the required contradiction.

The third claim follows from the second claim. □

We now introduce the predicate $M(x, y, a, b, c, d)$ which holds if and only if $(x, y, \{z \in A^2\colon |z| < a \,\&\, S(b, G(z))\},\ \{w \in A^2\colon |w| < c \,\&\, S(d, G(w))\})$ is a preferred ich.

Lemma 5.81. *$M$ is definable.*

*Proof.* Clearly $M(x, y, a, b, c, d)$ is equivalent to the fact that $(x, y, \{z \in A^2\colon |z| < a \,\&\, S(b, G(z))\},\ \{w \in A^2\colon |w| < c \,\&\, S(d, G(w))\})$ is an ich $X$ such that for all $a^+ \leq x$, $X \mid a^+ = (X \mid a)^\#$. It remains to show that the assertion $(x, y, \{z \in A^2\colon |z| < a \,\&\, S(b, G(z))\},\ \{w \in A^2\colon |w| < c \,\&\, S(d, G(w)\})$ is an ich is definable.

The only clause in the definition of ich that is problematic is clause (vi). However, according to Lemma 5.74 (v), the sets $Y$ appearing there are the sets of the form $\{y < z\colon S(w, y)\}$, for some $z$, $w$, which removes the difficulty. □

Lemma 5.82. *For all $x$ there exist $y$, $a$, $b$, $c$, $d$ such that $M(x, y, a, b, c, d)$. There is a unique preferred* ich *of each height.*

*Proof.* By Lemma 5.80, there is at most one preferred ich of any given height. We now prove the first claim by transfinite induction, which follows from the least element principle and Lemma 5.81.

Fix $x$ such that for all $x' < x$, there exists $y'$, $a'$, $b'$, $c'$, $d'$ such that $M(x', y', a', c', c', d')$. If $x = 0$, then there obviously exists $y$, $a$, $b$, $c$, $d$ such that $M(x, y, a, b, c, d)$. If $x = x'^+$, then let $M(x', y', a', b', c', d')$. Let $X$ be the preferred ich of height $x'$ coded by $x'$, $y'$, $a'$, $b'$, $c'$, $d'$. We can choose $x$, $y$, $a$, $b$, $c$, $d$ which appropriately codes the preferred ich $X^\#$ of height $x$.

Finally, suppose $x$ is a limit point in $A$. By the least element principle and Lemma 5.81, there is a definable map $h\colon [0,x) \to A^6$ such that for all $x' < x$, $M(h(x'))$. By Lemma 5.76, the range of $h$ is bounded.

Now for $x' < x$, $h(x')$ codes the unique preferred ich of height $x'$, and these preferred ich's are comparable according to Lemma 5.80. Let $(x, y, R, \mathrm{rk})$ be such that $y$ is the sup of the domain points of the $h(x')$, $R$ is the union of the epsilon relations of the $h(x')$, and rk is the union of the rank functions of the $h(x')$. Then $(x, y, R, \mathrm{rk})$ is a preferred ich of height $x$. By Lemma 5.74 (v), choose $a, b, c, d$ so that $(x, y, a, b, c, d)$ appropriately codes up $X$. Hence $M(x, y, a, b, c, d)$ as required. This establishes the first claim. The second claim follows. □

Let $R$ be the union of the epsilon relations of the preferred ich's, and let rk be the union of the rank functions of the preferred ich's.

LEMMA 5.83. $(A, R, \mathrm{rk})$ *is a system that is definable without parameters in* $(A, <, S, G)$, *and has the following properties*:
(i) $(A, R)$ *is extensional*;
(ii) rk *is a rank function for* $(A, R)$ *which is onto* $A$;
(iii) *The inverse image of every bounded subset of* $A$ *under* rk *is bounded*;
(iv) *Let* $Y \subseteq A$ *be bounded and definable in* $(A, <, S, G)$.
*Then there exists* $a \in A$ *such that* $Y = \{x\colon R(x, a)\}$.

*Proof.* Observe that by a straightforward transfinite induction argument, if $(x, y, R, \mathrm{rk})$ is a definable ich, then $x \leq y$. Clauses (i), (ii), (iv) follow from Lemmas 5.80 and 5.82. For clause (iii), note that $\mathrm{rk}(u) < x$ is equivalent to $u$ being in the domain of the preferred ich of height $x$. □

5.8. *Existence of a model of set theory with large cardinals, and application of the second incompleteness theorem.*

LEMMA 5.84. $(A, R)$ *satisfies* ZFC. $(A, R)$ *satisfies strong separation in the sense that for all bounded* $Y \subseteq A$ *that are definable in* $(A, <, S, G)$, *there exists* $x \in A$ *such that* $Y = \{a\colon R(a, x)\}$.

*Proof.* Extensionality is from clause (iii) in the definition of ich. Pairing is from clause (vi) in the definition of ich.

Separation is from clause (vi) in the definition of ich and the fact that $R$ itself is definable in $(A, <, S, G)$. In fact, this yields the strong form of separation, as stated in the lemma.

For union, first note that $R(a, b)\ \&\, R(b, c)$ implies $\mathrm{rk}(a) < \mathrm{rk}(c)$. The union axiom follows from Lemma 5.83 (iii) and strong separation.

For power set, note that by clause (v) in the definition of ich, if $(\forall x)(R(x,a)$ implies $R(x,b))$, then $\mathrm{rk}(a) \leq \mathrm{rk}(b)$. Power set then follows from Lemma 5.83 (iii) and strong separation.

For replacement, first use the fact that every definable function with bounded domain has bounded range from Lemma 5.76. Then apply Lemma 5.83 (iii) and strong separation. (In the presence of separation, replacement is taken in the form $(\forall x \in a)(\exists !y) \Rightarrow (\exists z)(\forall x \in a)(\exists y \in z)$.)

For foundation, we simply use the rank function, since $R(x,y)$ implies $\mathrm{rk}(x) < \mathrm{rk}(y)$.

We now come to the axiom of infinity: there exists a set $x$ such that the empty set is in $x$, and for all sets $a \in x$, $a \cup \{x\} \in x$. Let $Y = \{a\colon \mathrm{rk}(a) < e_1\}$. Apply Lemma 5.83 (iii) and strong separation, together with the fact that $e_1$ is a limit point, to obtain the required $x$.

We have now shown that $(A,R)$ satisfies ZF, which allows us a certain amount of liberty. To verify the axiom of choice, we start with what is regarded within $(A,R)$ as a pairwise disjoint collection of nonempty sets. This translates externally (i.e., in $(A,<,S,G)$) into an appropriate binary relation that is bounded and definable. We can obviously find a choice set that is definable in $(A,<,S,G)$ by taking least elements. Now apply Lemma 5.83 (iii) and strong separation. □

LEMMA 5.85. *Proposition* A *for #-decreasing implies the consistency of* ZFC, *provably within* ZFC. *If* ZFC *is consistent, then Proposition* A *for #-decreasing cannot be proved in* ZFC.

*Proof.* The first claim is by Lemma 5.84, and the second claim is by the Gödel second incompleteness theorem. □

We now want to show that $(A,R)$ satisfies the existence of subtle cardinals of every standard finite order, thereby showing that Proposition A for #-decreasing cannot be proved in ZFC + {there exists a subtle cardinal of order $k\}_k$, assuming this augmented theory is consistent. Recall the definition of $k$-subtle cardinals from [Ba75] discussed in Part 2: $\lambda$ is $k$-subtle if and only if

(i) $\lambda$ is an infinite cardinal;

(ii) for all closed unbounded $C \subseteq \lambda$ and regressive $f\colon S_k(\lambda) \to S(\lambda)$, there exists an $f$-homogeneous $B \in S_{k+1}(C)$.

LEMMA 5.86. *Let* $i > 0$, $x < e_i$, *and* $F\colon [0,x) \to [0,e_i)$ *be definable in* $(A,<,S,G)$. *Then* $\mathrm{rng}\,(F)$ *is bounded in* $[0,e_i)$.

*Proof.* Let $i > 0$ and suppose that this is false for some $x$ and $F$. Using the representation of bounded definable sets in Lemma 5.74 (v), we can find a counterexample $x$, $F$, where $x$, $F$ are definable in $(A,<,S,G)$ using only the parameter $e_i$. Now make the same definitions with the parameter $e_i$ replaced by the parameter $e_{i+1}$. By Lemma 5.74 (vii), this results in the same $x$, $F$. Now consider the true statement

$$\mathrm{rng}\,(F) \text{ is bounded in } [0, e_{i+1})$$

formulated with the definition of $F$ using only the parameter $e_{i+1}$. By Lemma 5.74 (vii), the statement

$$\mathrm{rng}\,(F) \text{ is bounded in } [0, e_i)$$

must hold, where the definition of $F$ using only the parameter $e_i$ is used. This is the required contradiction. □

We say that $C$ is a closed unbounded subset of $[0, e_i)$ if and only if
(i) $C \subseteq [0, e_i)$;
(ii) for all $x < e_i$ there exists $x < y \in C$;
(iii) for all $0 < x < e_i$, if $(\forall\, y < x)(\exists\, z \in C)(y < z < x)$, then $x \in C$.

Lemma 5.87. *Let $C$ be a closed unbounded subset of $[0, e_i)$ and $F$: $S_k([0, e_i)) \to S([0, e_i))$ be regressive. Assume that $C$ is definable in $(A, <, S, G)$. Also assume that the "characteristic relation" $x \in F(y)$ of $F$ is definable in $(A, <, S, G)$. Then there exists an $F$-homogeneous $B \in S_{k+1}(C)$.*

*Proof.* Let $i > 0$, and suppose that there is a definable counterexample $C$, $F$. Using Lemma 5.74 (v), we can find a counterexample $C$, $F$, where $C$, $F$ are definable in $(A, <, S, G)$ using only the parameter $e_i$. Now make the same definitions with $e_i$ replaced by $e_{i+k+1}$. By Lemma 5.74 (vii), this results in $C'$, $F'$, where
(1) $C'$ is a closed unbounded subset of $[0, e_{i+k+1})$;
(2) $F'\colon S_k([0, e_{i+k+1})) \to S([0, e_{i+k+1}))$ is regressive;
(3) $C = C' \cap [0, e_i)$;
(4) the characteristic relation of $F$ is the restriction of the characteristic relation of $F'$ to $[0, e_i)$.

Hence $e_i \in C'$. So again by Lemma 5.74 (vii), $e_i, e_{i+1}, \ldots, e_{i+k} \in C'$. It is also clear by Lemma 5.74 (vii) that $\{e_i, \ldots, e_{i+k}\}$ is $F'$-homogeneous. Hence

$$\text{there exists an } F'\text{-homogeneous } B \in S_{k+1}(C')$$

holds by the definitions with only the parameter $e_{i+k+1}$. Therefore, by Lemma 5.74 (ii),

$$\text{there exists an } F\text{-homogeneous } B \in S_{k+1}(C)$$

holds using the definitions with only the parameter $e_i$. This is the required contradiction. □

Lemma 5.88. *For all $a \in A$ there exists a unique $a^* \in A$ such that* $\mathrm{rk}(a^*) = a$ *and "$a^*$ is epsilon-connected and transitive" holds in $(A, R)$. (The* *statement in quotes is also written "$a^*$ is an ordinal.") The $^*$-operation is an order-preserving bijection from $A$ onto the set of ordinals of $(A, R)$ (with the ordinals of $(A, R)$ ordered by $R$), and is definable without parameters in $(A, <, S, G)$.*

*Proof.* Since $(A, R)$ satisfies ZFC, it satisfies the classical theorem about comparability of ordinals – for any two distinct ordinals, the first is an element of the second or the second is an element of the first. Thus their ranks must be different. This establishes uniqueness.

We prove existence by transfinite induction in $(A, <, S, G)$, which is valid in this context because by Lemma 5.83, $(A, R, \mathrm{rk})$ is definable in $(A, <, S, G)$. Let $a \in A$, and assume that for all $b < a$, there exists (unique) $b^*$ such that $\mathrm{rk}(b^*) = b$ and $b^*$ is an ordinal in $(A, R)$. If $a = 0$, then we can take the empty set in $(A, R)$, and if $a = b^+$, we can take $b^* \cup \{b^*\}$ in $(A, R)$.

Now assume that $a$ is a limit. Let $Y = \{b^*: b < a\}$. By the strong separation proved in Lemma 5.84, let $a'$ be such that $Y = \{c: R(c, a')\}$. Then $\mathrm{rk}(a') = a$. We also claim that $(A, R)$ satisfies that $Y$ is transitive. To see this, let $R(x, b^*)$, $b < a$. Since $b^*$ is an ordinal in $(A, R)$, so is $x$. Let $\mathrm{rk}(x) = c < \mathrm{rk}(b^*) = b$. By the uniqueness in the previous paragraph, $x = c^*$, and so $x \in Y$. Obviously $a'$ is a transitive set of ordinals in $(A, R)$, and so $a'$ is an ordinal in $(A, R)$. Also $\mathrm{rk}(a')$ is the strict sup of the $\mathrm{rk}(b^*)$, $b < a$, which is clearly $a$.

For the second claim, the surjectivity is a consequence of the fact that every $x \in A$ has a rank $\mathrm{rk}(x)$. To check order preservation, suppose $x < y$. Then $R(x^*, y^*)$ or $R(y^*, x^*)$. The latter is impossible since otherwise $\mathrm{rk}(y^*) = y < \mathrm{rk}(x^*) = x$. The definability follows from Lemma 5.83. □

We are not claiming that $*$ is continuous.

LEMMA 5.89. *For all $i > 0$, "$e_i^*$ is a regular cardinal" holds in $(A, R)$.*

*Proof.* It suffices to prove that in $(A, R)$, the range of every function from an element of $e_i^*$ into $e_i^*$ has a strict upper bound $< e_i^*$. Let $f$, $x$ be such that in $(A, R)$, "$x < e_i^*$ and $f$ maps $x$ into $e_i^*$." By Lemma 5.88, let $x = a^*$, $a < e_i$. By Lemma 5.88, when viewed externally, $f$ defines a function $F: [0, a) \to e_i$ which is definable in $(A, <, S, G)$, where for all $b < a$, $f(b^*) = F(b)^*$ holds in $(A, R)$. The lemma then follows from Lemmas 5.86 and 5.88. □

LEMMA 5.90. *For all $i, k > 0$, "$e_i^*$ is a $k$-subtle cardinal" holds in $(A, R)$.*

*Proof.* In $(A, R)$, let $C$ be a closed unbounded subset of $e_i^*$, and $f: S_k(e_i^*) \to S(e_i^*)$. Let $C' = \{a: a^* \in C\}$. Then by Lemma 5.88, $C'$ is a closed unbounded subset of $[0, e_i^*)$ that is definable in $(A, <, S, G)$. Let $F: S_k([0, e_i)) \to S([0, e_i))$ be given by $F(x) = \{a: a^* \in f(x)\}$. Then the hypotheses of Lemma 5.87 apply to $C'$, $F$, and $e_i$. Let $B \in S_{k+1}(C')$ be $F$-homogeneous. Then in $(A, R)$, $\{b^*: b \in B\}$ is an $f$-homogeneous element of $S_{k+1}(C)$ as required. □

THEOREM 5.91. *Proposition* A *for #-decreasing implies the consistency of* ZFC + {*there exists a k-subtle cardinal*}$_k$, *provably within* ZFC. *If* ZFC + {*there exists a k-subtle cardinal*}$_k$ *is consistent, then Proposition* A *for #-decreasing cannot be proved in* ZFC+{*there exists a k-subtle cardinal*}$_k$. *In fact, Proposition* A *for #-decreasing cannot be proved from any consistent subset of the theorems of* ZFC + {*there exists a k-subtle cardinal*}$_k$ *from which* ZFC *can be derived. The same results hold for all forms of Propositions* A–D.

*Proof.* The first claim is proved by Lemmas 5.84 and 5.90. The second claim is proved by Gödel's second incompleteness theorem.

For the third claim, we use $A$ for "Proposition A for #-decreasing." We fix a finite fragment ZFC′ of ZFC which suffices to prove the implication "$A$ implies the consistency of ZFC + {there exists a $k$-subtle cardinal}$_k$." Let $T$ be any consistent subset of the theorems of ZFC + {there exists a $k$-subtle cardinal}$_k$ that derives ZFC + $A$. Let $T'$ be a finite subset of $T$ that derives ZFC′ + $A$. Then ZFC + $T'$ is consistent since its axioms are included in the theorems of $T$. Also ZFC + $T'$ proves the consistency of ZFC + {there exists a $k$-subtle cardinal}$_k$. Now $T'$ is a finite set of theorems of ZFC + {there exists a $k$-subtle cardinal}$_k$, and hence also ZFC + $T'$ proves this fact. Hence ZFC + $T'$ proves its own consistency. This contradicts the consistency of ZFC + $T'$ by the Gödel second incompleteness theorem.

The fourth claim follows immediately from Theorem 3.14. □

The following result gives a precise sense to the word "necessary" in the title of this paper.

COROLLARY 1. *Any extension of* ZFC *that suffices to prove Proposition* A *for #-decreasing is an extension of* ZFC *in which* ZFC + {*there exists a k-subtle cardinal*}$_k$ *is interpretable. The same result holds for all forms of Propositions* A–D.

*Proof.* This follows immediately from the first claim of Theorem 5.91 and the usual proof of Gödel's completeness theorem for first order predicate calculus with equality. □

THE OHIO STATE UNIVERSITY, COLUMBUS, OH
*E-mail address*: friedman@math.ohio-state.edu